\documentclass[12pt]{article}
\usepackage{amsmath,amssymb,amsthm}
\usepackage[round]{natbib}
\bibpunct{(}{)}{;}{a}{}{,}
\usepackage{graphicx,booktabs,setspace,url}
\usepackage[margin=1in]{geometry}
\usepackage[hidelinks]{hyperref}
\usepackage{xcolor}

\newcommand{\blind}{0}

\newtheorem{theorem}{Theorem}

\newtheorem{assumption}{Assumption}

\newtheorem{corollary}{Corollary}

\newcommand{\R}{\mathbb{R}}
\newcommand{\E}{\mathbb{E}}
\newcommand{\Prob}{\mathbb{P}}
\newcommand{\calH}{\mathcal{H}}
\newcommand{\calC}{\mathcal{C}}
\newcommand{\calK}{\mathcal{K}}
\newcommand{\calN}{\mathcal{N}}
\newcommand{\calE}{\mathcal{E}}
\newcommand{\calT}{\mathcal{T}}
\newcommand{\calG}{\mathcal{G}}

\newcommand{\calW}{\mathcal{W}}
\newcommand{\dist}{\operatorname{dist}}
\newcommand{\tr}{\operatorname{tr}}
\newcommand{\supn}[1]{\|#1\|_{\infty}}
\newcommand{\norm}[1]{\|#1\|}
\newcommand{\wstar}{w^{\star}}
\newcommand{\wproj}{w^{\circ}}
\newcommand{\bproj}{\beta^{\circ}}
\newcommand{\what}{\hat w}
\newcommand{\whatp}{\hat w^{+}}
\newcommand{\bhat}{\hat\beta}
\newcommand{\Sighat}{\hat\Sigma}
\newcommand{\Ghat}{\hat G}
\newcommand{\Fv}{\mathbf{F}}
\newcommand{\ind}{\mathbf{1}}

\graphicspath{{figures/}}

\newtheorem*{thmprime}{Theorem~$1'$}
\newtheorem*{propprime}{Proposition~$2'$}
\newtheorem*{thmsecond}{Theorem~$1''$}
\newtheorem*{propsparse}{Proposition~$2''$}
\newtheorem{lemma}{Lemma}

\newtheorem{propositionS}{Proposition}

\newtheorem{theoremS}{Theorem}

\newcommand{\TV}{\operatorname{TV}}
\newcommand{\KL}{\operatorname{KL}}
\newcommand{\Var}{\mathrm{Var}}
\begin{document}

\def\spacingset#1{\renewcommand{\baselinestretch}{#1}\small\normalsize}
\spacingset{1}

\if0\blind
{
  \title{\bf Graphon-Level Bayesian Predictive Synthesis for Random Network}
  \author{Marios Papamichalis\thanks{Human Nature Lab, Yale University, New Haven, CT 06511, marios.papamichalis@yale.edu}
  \hspace{.2cm}\\
  {}Human Nature Lab, Yale University\\
  and \\
  Regina Ruane\thanks{Department of Statistics and Data Science, The Wharton School, University of Pennsylvania, Philadelphia, PA 19104, ruanej@wharton.upenn.edu} \\
  {}Department of Statistics and Data Science\\
  {}The Wharton School, University of Pennsylvania}
  \date{}
  \maketitle
} \fi

\if1\blind
{
  \bigskip
  \bigskip
  \bigskip
  \begin{center}
    {\LARGE\bf Graphon-Level Bayesian Predictive Synthesis for Random Networks}
  \end{center}
  \medskip
} \fi

\bigskip
\begin{abstract}
Analysts often fit several models to the same network. Each estimates a graphon,
the function giving the probability of a link between two nodes. Reporting one
predictive graphon means weighting these estimates, and the weights can be
constrained in several ways. We determine which constraint is correct and when
the combination improves on the best single model. We introduce Bayesian predictive synthesis at the graphon level. The models enter
as agents, a prior is placed on the weights combining their link probabilities,
and the posterior returns one predictive graphon with credible intervals. We
study this rule when the network is a union of overlapping mechanisms, so a pair
is linked if at least one mechanism links it. Weights forced to sum to one cannot
reproduce a union. Free and nonnegative weights do far better and do equally
well. A noisy-OR rule matches the union exactly. From one observed graph we give
an exact variance formula for the fitted weights, and the usual credible
intervals are too narrow. On polblogs, email-Eu-core, Cora, ca-GrQc, wiki-Vote and soc-Epinions1,
combining adds under one percent once a flaw in the usual benchmark is corrected,
and a larger single model often wins. On the Amazon and YouTube multiplex
networks, where the layers are recorded separately, it beats every competitor we
tuned, including spectral estimation, blockmodels, random forests and node2vec.

\end{abstract}

\noindent%
\vfill

\newpage
\spacingset{1.9}

\section{Introduction}
\label{sec:intro}

Convex averaging is the wrong way to combine network models when the mechanisms they describe operate at the same time. A collaboration network may be well described by a blockmodel for its communities, a degree-corrected model for its hubs and a latent-distance model for its short cycles; each fitted model returns a probability for every pair of nodes, and each is right about one mechanism and silent about the others. The standard options, selecting the model that scores best on held-out edges or averaging the models with weights in the probability simplex, both treat the models as competing explanations of the same edges. If instead the network is the union of the mechanisms, so that a pair is connected when at least one mechanism connects it (the mechanism indicators conditionally independent given the positions), the connection probability is $1-\prod_j(1-w_j)$, which for weak mechanisms equals the sum of the $w_j$ minus a product term. No convex combination can reproduce a sum, and Theorem~\ref{thm:superposition} measures the cost exactly: the squared-error floor of every simplex-constrained rule is of exact order $\rho^2$ in the mechanism strength $\rho$, while the floor of the affine combination $\beta_0+\sum_j\beta_jw_j$ is of exact order $\rho^4$ whenever the pairwise interaction $\sum_{j<k}w_jw_k$ leaves the span, with constant its squared distance to the span; the nonnegative cone with a free intercept attains the same fourth-order floor (Theorem~$1'$), so the loss comes from the sum-to-one constraint alone.

The combination is defined on graphons: by the Aldous--Hoover representation \citep{Aldous1985, Hoover1979}, an exchangeable random graph is generated by latent node positions and a symmetric function of pairs of positions, and the fitted probability matrices of standard network models are evaluations of such functions. Combining models at this level is the graphon version of Bayesian predictive synthesis \citep{McAlinnWest2019, McAlinnEtAl2020, TallmanWest2024}, in which a decision maker passes agent forecasts through a synthesis function whose parameters carry a posterior. Ours is affine in the agent graphons, and its posterior mean estimates, and contracts around, the orthogonal projection of the true graphon onto the agent span.

The second result (Theorem~\ref{thm:onegraph}) is where the graphon structure earns its place. All dyads of one graph share nodes, and we prove an exact identity for the variance of the synthesis score over the $N=n(n-1)/2$ dyads: it is $c(\wstar)/N$, the independent-dyad value, plus a term $4(n-2)\kappa_1(\wstar)/\{n(n-1)\}$, asymptotically $4\kappa_1(\wstar)/n$, in which $\kappa_1(\wstar)$ is a node-level moment of the projection residual, zero whenever the truth lies in the span and possibly zero outside it. The estimation error is of order $n^{-2}$ when the moment vanishes and $n^{-1}$ otherwise, while the working posterior scales as $n^{-2}$ in both: for every functional loaded on the residual moment its credible intervals lose coverage at rate $\sqrt n$, and a dyadic-robust covariance, consistent by Proposition~$2'$, restores it. Theorem~S1 gives the independent-dyad theory, the rate $d/m$, its minimaxity over the span, and a finite-sample bound whose Gram term is controlled by a leverage constant of the library; Proposition~S2 characterizes when credible intervals are conservative.

The third contribution is a negative empirical finding. On six public networks, an edge-holdout protocol that fits the agents to a graph with $30\%$ of the edges removed and then learns free weights on held-out dyads reports Brier gains over the selected agent of $0.25$ to $3.6\%$ of the density, with a weight near $1/0.7$ on the dominant agent: the synthesis is undoing the thinning. After dividing by the retention, the gains over the affine-calibrated selected agent are $0.11$ to $0.64\%$ of the density on five networks and $-0.03$ on the sixth, and enlarging the dominant model beats combining on four. The correction shrinks the measured gain by an order of magnitude on four networks, raises it on Cora, and flips its sign on soc-Epinions1. We report both sets of numbers; the weights are predictive projection coefficients, not mechanism strengths. The superposition model is then applied to three multiplex networks whose observed graph is a union of separately recorded layers: with one agent per layer the span weights recover the union on the sparse multiplex, the span beats the fitted hull by $4$ to $16\%$ of the density, and the full library beats single models of fixed size, in affine form on the sparse multiplex and in the exact noisy-OR form on the dense one; the comparators are tuned on the same held-out dyads as the synthesis, so the margins are not an artifact of a handicapped competitor.

\citet{LiLe2024} is the closest prior work and the benchmark this paper measures itself against: it combines fitted network probability matrices by unrestricted and by nonnegative least squares, selects the weights on held-out dyads, fits the models to the thinned graph and divides the fitted probabilities by the retention, and proves adaptivity and oracle bounds for the mixture. The least-squares and cone rows of our tables are their OLS and nonnegative mixers on this library, their exponential-weight mixer is run beside them in Section~\ref{sec:real}, and our protocol shares their zero-fill-and-rescale device. What this paper adds is the approximation geometry, floors of exact orders $\rho^2$ and $\rho^4$ over simplex, cone and span that also explain why their nonnegative constraint loses nothing here (Theorem~$1'$), the one-graph variance identity with its inferential consequences, the posterior over weights, the quantification of the uncorrected protocol's artifact, and the disjoint-data regime. Unrestricted linear combination of forecasts goes back to \citet{GrangerRamanathan1984}, \citet{Breiman1996} found nonnegativity the useful constraint in stacking, and the rates attached to selection, hull and span are the classical aggregation trichotomy \citep{Tsybakov2003, BuneaTsybakovWegkamp2007}. Linear pools of calibrated forecasts are uncalibrated \citep{GenestZidek1986, GneitingRaftery2007, RanjanGneiting2010, GneitingRanjan2013}, and \citet{Claeskens2016} explain why estimated unconstrained weights often lose to constrained ones, the mechanism behind the small margins of Section~\ref{sec:real}. Stacking of predictive distributions \citep{YaoEtAl2018} and convex model averaging for networks \citep{Heo2025, zhang2025network} live in the simplex; random-forest network stacking does not \citep{GhasemianEtAl2020, HeEtAl2024}, and \citet{LatoucheRobin2016} average over blockmodels. The union truth of Section~\ref{sec:theory} is the noisy-OR link of overlapping-community models \citep{LatoucheBirmeleAmbroise2011, BallKarrerNewman2011, VallesCatala2016}; the contribution here is its combination geometry. Graphon estimation \citep{OlhedeWolfe2014, GaoLuZhou2015, Chatterjee2015, KloppTsybakovVerzelen2017} recovers one graphon nonparametrically, \citet{LiLevinaZhu2020} give the edge cross-validation our splits build on, and \eqref{eq:dyadic} sits in the dyadic-robust literature \citep{AronowSamiiAssenova2015, Graham2020}. \if0\blind Two related preprints treat fitted mechanisms with cross-fitting \citep{PapamichalisRuane2026minimax} and dynamic networks \citep{PapamichalisRuaneEtAl2026dynamic}; the theory of the present paper treats the agents as known, and Theorem~\ref{thm:crossfit} states what survives when they are fitted.\else Two companion papers, withheld for review, treat fitted mechanisms with cross-fitting and dynamic networks; the theory of the present paper treats the agents as known, and Theorem~\ref{thm:crossfit} states what survives when they are fitted.\fi

Section~\ref{sec:setup} defines the synthesis and the two sampling designs, Section~\ref{sec:theory} states the results, Section~\ref{sec:sims} reports simulations designed to falsify the identities, Section~\ref{sec:real} reports the six networks with and without the retention correction, beside external topological baselines and a disjoint-data experiment, and Section~\ref{sec:discussion} records limitations. Proofs and additional tables are in the supplement.

\section{Graphon-level synthesis}
\label{sec:setup}

\paragraph{Graphons and jointly aligned agents.}
A graphon is a measurable symmetric function $w:[0,1]^2\to[0,1]$; it generates a random graph on $n$ nodes by drawing latent positions $U_1,\dots,U_n$ independently and uniformly on $[0,1]$ and connecting nodes $i$ and $j$ independently with probability $w(U_i,U_j)$ \citep{Aldous1985, Hoover1979, LovaszSzegedy2006}. We write $\norm{\cdot}$ and $\langle\cdot,\cdot\rangle$ for the norm and inner product of $L_2([0,1]^2)$, $\supn{\cdot}$ for the essential supremum, and $\dist(w,\mathcal{S})=\inf_{h\in\mathcal{S}}\norm{w-h}$; the same bars denote the Euclidean norm on $\R^d$, distinguished by the argument. The edge, triangle and wedge densities, the clustering ratio and their distinct-tuple graph versions are defined in the supplement, which also proves the transfer of the guarantees to them. Graphons are equivalence classes up to measure-preserving maps; pointwise statements refer to a fixed representative and hold for almost every pair of positions. Agent $j$ is a model for the network, identified with its graphon $w_j$, $j=1,\dots,J$. A linear combination of two graphons changes if one is relabeled independently of the other \citep{KloppTsybakovVerzelen2017}, so the synthesis acts on agents sharing one latent coordinate and is invariant under a common relabeling. Fitted models provide exactly this: fitted to a common node set, each returns a probability for every pair of nodes, so the synthesis combines node-aligned probability matrices. The feature vector of a pair of positions is
\[
\Fv(u,v)=\big(1,\,w_1(u,v),\dots,w_J(u,v)\big)^{\top}\in\R^{d},\qquad d=J+1,
\]
so that $\beta^{\top}\Fv=\beta_0+\sum_{j}\beta_j w_j$. Three sets organize the paper: the span $\calH=\{\beta^{\top}\Fv:\beta\in\R^{d}\}$, the hull $\calC=\{\sum_j\pi_j w_j:\pi\in\Delta_J\}$ with $\Delta_J$ the probability simplex, and the agents $\{w_1,\dots,w_J\}\subset\calC\subset\calH$. The Gram matrix is $G=\int\Fv\Fv^{\top}$.

\begin{assumption}[Linearly independent agents]
\label{ass:gram}
The functions $1,w_1,\dots,w_J$ are linearly independent in $L_2([0,1]^2)$; equivalently $\lambda:=\lambda_{\min}(G)>0$.
\end{assumption}

Two constants of the library enter the finite-sample bounds: $\lambda\le1$, and the leverage $L=\supn{\Fv^{\top}G^{-1}\Fv}$, which satisfies $d\le L\le d/\lambda$ because $\int\Fv^{\top}G^{-1}\Fv=d$ and $\norm{\Fv}^2\le d$. The leverage is much smaller than $d/\lambda$ for the libraries used below ($26$, $27$ and $60$ against $133$, $1014$ and $4901$ for $d=3,5,9$), which is what makes the bounds informative at the simulated sizes.

\paragraph{Truth, projection and residual moments.}
The data come from an unknown graphon $\wstar$, which need not lie in $\calH$. Its orthogonal projection onto $\calH$,
\[
\wproj=(\bproj)^{\top}\Fv,\qquad \bproj=G^{-1}\!\int\Fv\,\wstar,
\]
is the unique minimizer of $\norm{w-\wstar}$ over $\calH$, satisfies the normal equations $\int(\wstar-\wproj)\Fv=0$, and is never farther from $\wstar$ than the best single agent (Lemma~S1). The residual second moment $\sigma^2=\wstar(1-\wstar)+(\wstar-\wproj)^2$, variance plus squared misspecification bias, enters through
\[
c(\wstar)=\int\sigma^2\,\Fv^{\top}G^{-1}\Fv,\qquad M=\int\sigma^2\,\Fv\Fv^{\top},\qquad \bar\sigma^2=\supn{\sigma^2},
\]
with $\tr(G^{-1}M)=c(\wstar)\le d\bar\sigma^2$ and $c(\wstar)\le d/4$ when $\wstar\in\calH$. The one-graph theory adds a node-level object, the residual moment
\begin{equation}
\label{eq:g}
\begin{gathered}
g(u)=\int_0^1(\wstar-\wproj)(u,v)\,\Fv(u,v)\,dv\in\R^{d},\\[2pt]
\kappa_1(\wstar)=\int_0^1 g(u)^{\top}G^{-1}g(u)\,du,\qquad \Sigma_g=\int_0^1 g(u)g(u)^{\top}du,
\end{gathered}
\end{equation}
the cross-moment of the projection residual with the agents along the row of position $u$. By the normal equations $\int g=0$; $\kappa_1(\wstar)=0$ if and only if $g=0$ almost everywhere, which holds whenever $\wstar\in\calH$.

\paragraph{Two sampling designs.}
In the independent-dyad design,
\begin{equation}
\label{eq:model}
(U_s,V_s)\stackrel{\text{iid}}{\sim}\mathrm{Unif}([0,1]^2),\qquad Y_s\mid U_s,V_s\sim\mathrm{Bernoulli}\big(\wstar(U_s,V_s)\big),\qquad s=1,\dots,m,
\end{equation}
and the synthesis observes $(\Fv_s,Y_s)$, $\Fv_s=\Fv(U_s,V_s)$. In the one-graph design,
\begin{equation}
\label{eq:onegraph}
U_1,\dots,U_n\stackrel{\text{iid}}{\sim}\mathrm{Unif}[0,1],\qquad A_{ij}\mid U\sim\mathrm{Bernoulli}\big(\wstar(U_i,U_j)\big)\ \text{independently for } i<j,
\end{equation}
and the synthesis observes $(\Fv_{ij},A_{ij})$ for all $N=n(n-1)/2$ dyads. Dyads from one large graph without shared nodes follow \eqref{eq:model}; all dyads of one graph follow \eqref{eq:onegraph}, and the difference is the subject of Theorem~\ref{thm:onegraph}. The latent positions are not observed and are not needed. The quantity controlled throughout is the uniform-dyad risk $\E\norm{\what-\wstar}^2$, the expected Brier score of $\what$ on a fresh uniform dyad minus the irreducible term $\int\wstar(1-\wstar)$.

\paragraph{The synthesis and its posterior.}
In Bayesian predictive synthesis \citep{McAlinnWest2019}, which builds on the agent-opinion analyses of \citet{West1984} and \citet{West1992}, the decision maker combines agent forecasts through a synthesis function with a prior on its parameters. Here the forecasts for a dyad are the probabilities $w_j(U_s,V_s)$, and the synthesis function is Gaussian and affine,
\begin{equation}
\label{eq:synthesis}
Y_s\mid\Fv_s,\beta\sim N(\beta^{\top}\Fv_s,\,\nu),\qquad \beta\sim N(0,\tau^2 I_d),
\end{equation}
with synthesis variance $\nu>0$ and prior scale $\tau>0$ chosen by the decision maker. This is the static point-forecast special case of \citet{McAlinnWest2019}: each agent's predictive distribution for a dyad is Bernoulli with mean $w_j(U_s,V_s)$, the synthesis function is linear and Gaussian in those means, and conditioning on them collapses the latent-agent integral of the general theory. What the construction keeps is the posterior over the combination weights; when it can be trusted is the subject of Theorem~\ref{thm:onegraph} and Proposition~S2. The posterior is $\beta\mid\text{data}\sim N(\bhat,\Sighat)$ with
\begin{equation}
\label{eq:posterior}
\bhat=(\Ghat+\kappa I_d)^{-1}\frac1m\sum_{s=1}^{m}\Fv_sY_s,\qquad \Sighat=\frac{\nu}{m}(\Ghat+\kappa I_d)^{-1},\qquad \Ghat=\frac1m\sum_{s=1}^{m}\Fv_s\Fv_s^{\top},\quad \kappa=\frac{\nu}{m\tau^2};
\end{equation}
under \eqref{eq:onegraph} the sums run over the $N$ dyads and $N$ replaces $m$. The posterior mean graphon is $\what=\bhat^{\top}\Fv$, and the reported graphon is $\whatp=\min\{\max(\what,0),1\}$; clipping is a pointwise contraction, so $\norm{\whatp-\wstar}\le\norm{\what-\wstar}$ always. The posterior over $\beta$ induces a posterior $\Pi_m$ over the affine kernels $w=\beta^{\top}\Fv\in\calH$, whose image under clipping is the posterior over graphons. The Gaussian likelihood is a working model in the sense of \citet{West1992}, not a coherent Bernoulli model, since unclipped draws can leave $[0,1]$; $\bhat$ is a ridge-regularized least-squares fit of binary outcomes on agent probabilities, and $\nu$ sets the posterior scale (Proposition~S2). Credible intervals from $\Pi_m$ are intervals for the projection $\wproj$ and its functionals, not for $\wstar$: when $\wstar\notin\calH$ their coverage of $\wstar$ tends to zero, and the approximation error $\norm{\wstar-\wproj}$ is the separate, deterministic term of every bound below.

Two features distinguish \eqref{eq:synthesis} from averaging rules: the intercept and the unconstrained weights let $\what$ leave the hull $\calC$, where model averaging \citep{Heo2025}, simplex stacking \citep{YaoEtAl2018} and selection live, and the weights carry a posterior. For the one-graph design we also use the dyadic-robust covariance
\begin{equation}
\label{eq:dyadic}
\Sighat_{\mathrm{dy}}=A^{-1}\hat V A^{-1},\qquad \hat V=\frac{1}{N^2}\Big\{\sum_{i=1}^{n}r_ir_i^{\top}-\sum_{i<j}\hat\xi_{ij}^2\Fv_{ij}\Fv_{ij}^{\top}\Big\},\qquad r_i=\sum_{j\neq i}\hat\xi_{ij}\Fv_{ij},
\end{equation}
with $A=\Ghat+\kappa I_d$ and $\hat\xi_{ij}=A_{ij}-\what(U_i,U_j)$ the fitted residuals; $\hat V$ sums residual cross-products over all pairs of dyads sharing a node, the empirical version of the variance identity of Theorem~\ref{thm:onegraph}(i). The subtraction can make $\hat V$ indefinite in finite samples; the implementation floors each directional variance at zero, and Section~\ref{sec:sims} reports how often the floor binds.

\begin{propprime}[Dyadic covariance: identity and consistency]
Under \eqref{eq:onegraph} and Assumption~\ref{ass:gram}: (i) exactly, for every $n\ge3$,
\[
\Var(S)\;=\;\frac{M}{N}\;+\;\frac{4(n-2)}{n(n-1)}\,\Sigma_g,
\]
of which Theorem~\ref{thm:onegraph}(i) is the trace against $G^{-1}$; (ii) with residuals taken at $\bproj$, $\E[\hat V]$ equals the right-hand side of (i) for every $n$, and with residuals at $\bhat$, $n\hat V\to4\Sigma_g$ elementwise in probability; (iii) consequently, for every $\ell$ with $\ell^{\top}G^{-1}\Sigma_gG^{-1}\ell>0$, $\ell^{\top}(\bhat-\bproj)/(\ell^{\top}\Sighat_{\mathrm{dy}}\ell)^{1/2}\Rightarrow N(0,1)$, so the dyadic-robust intervals are asymptotically valid exactly where Theorem~\ref{thm:onegraph}(iii) shows the working posterior fails.
\end{propprime}

\section{Theory}
\label{sec:theory}

The results answer four questions in order: what does the simplex constraint cost when mechanisms superpose; what happens when all dyads come from one graph; how fast is the synthesis from independent dyads and can anything beat it; and when are the credible intervals trustworthy. Proofs are in Section~S1 of the supplement.

\subsection{Superposition: the cost of the hull}
\label{sec:hull}

Simplex-constrained synthesis, which we call hull synthesis, reports $\what_{\calC}=\sum_j\hat\pi_jw_j$ with $\hat\pi$ minimizing the empirical Brier score over $\Delta_J$; selection reports the best single agent under the same criterion. The superposition model makes the comparison sharp: the truth is the union of the agents' mechanisms, $\wstar=1-\prod_{j=1}^{J}(1-w_j)$, the probability that at least one mechanism connects the pair when the mechanism indicators are conditionally independent given the positions, and each mechanism is scaled as $w_j=\rho\tilde w_j$ with $0\le\tilde w_j\le1$ and strength $\rho\in(0,1]$. Expanding, $\wstar=\rho\sum_j\tilde w_j-\rho^2 e_2+O(\rho^3)$ with $e_2=\sum_{j<k}\tilde w_j\tilde w_k$ the pairwise interaction: the linear term lies in $\calH$, so the span misses only the interaction, while a convex mixture cannot reproduce even the linear term, because the mechanisms sum to more than any average of them.

\begin{theorem}[Superposition]
\label{thm:superposition}
Let Assumption~\ref{ass:gram} hold.
\begin{enumerate}
\item[(i)] For every realization of the data, $\norm{\what_{\calC}-\wstar}^2\ge\dist^2(\wstar,\calC)\ge\norm{\wstar-\wproj}^2$ and $\norm{\what_{\mathrm{sel}}-\wstar}^2\ge\min_j\norm{w_j-\wstar}^2\ge\dist^2(\wstar,\calC)$, and $\dist^2(\wstar,\calC)=\norm{\wstar-\wproj}^2$ if and only if $\wproj\in\calC$. These floors hold for every sample size.
\item[(ii)] In the superposition model, with $a=\sum_{j<k}\min(\norm{\tilde w_j},\norm{\tilde w_k})$, $b=\min_j\norm{\tilde w_j}$ and $c_J=2^J-1-J-J(J-1)/2$,
\[
\big|\,\norm{\wstar-\wproj}-\rho^2\dist(e_2,\calH)\,\big|\;\le\;c_J\rho^3,\qquad \norm{\wstar-\wproj}\le\rho^2a,
\]
\[
\rho\Big[\tfrac{J-1}{\sqrt J}\,b-\rho a\Big]_{+}\;\le\;\dist(\wstar,\calC)\;\le\;\rho\Big[\sum_{j\ge2}\norm{\tilde w_j}+\rho a\Big].
\]
\item[(iii)] Consequently, if the interaction is not spanned, $e_2\notin\calH$, then
\[
\norm{\wstar-\wproj}^2=\rho^4\dist^2(e_2,\calH)\{1+O(\rho)\},
\]
and for $J=2$ the identity $\norm{\wstar-\wproj}^2=\rho^4\dist^2(\tilde w_1\tilde w_2,\calH)$ is exact for every $\rho$; if $b>0$ and $\rho\le(J-1)b/(2\sqrt J\,a)$, then $\dist^2(\wstar,\calC)\ge\rho^2(J-1)^2b^2/(4J)$. The two floors are of exact orders $\rho^4$ and $\rho^2$.
\end{enumerate}
\end{theorem}

Part (i) is deterministic and locates the entire asymptotic gain of the span at the event $\wproj\notin\calC$: the projection needs an intercept, a negative weight, or weights summing to other than one. Part (iii) is the sentence the paper is about: when mechanisms superpose, the simplex is the wrong constraint set at first order in the mechanism strength, and the $\rho^4$ constant is a computable geometric quantity. Both constants are checked exactly in Section~\ref{sec:sims}: on the two-mechanism example of Figure~S2, $\norm{\wstar-\wproj}^2/\rho^4=\dist^2(\tilde w_1\tilde w_2,\calH)=3.00\times10^{-3}$ at every $\rho$, and the hull floor sits between its bounds with slope two. The condition $e_2\notin\calH$ is diagnosed numerically by quadrature of $\dist^2(e_2,\calH)$, a check, not a certificate. The expansion constant $c_J$ grows exponentially in $J$, so the $O(\rho)$ term in (iii) is uniform only over fixed $J$.

\begin{thmprime}[Four constraint sets]
Let Assumption~\ref{ass:gram} hold with $J\ge2$ and adopt the superposition model of Theorem~\ref{thm:superposition}. Beside the simplex $\calC$ and the span $\calH$ consider the intercept simplex $\calC_c=\{c+\sum_j\pi_jw_j:c\in\R,\ \pi\in\Delta_J\}$, the nonnegative cone with free intercept $\calK=\{c+\sum_j\gamma_jw_j:c\in\R,\ \gamma_j\ge0\}$, the cone without intercept $\calK_0=\{\sum_j\gamma_jw_j:\gamma_j\ge0\}$, and the span without intercept $\calH_0=\mathrm{span}(w_1,\dots,w_J)$; write $\tilde{\calH}_0=\mathrm{span}(\tilde w_1,\dots,\tilde w_J)$. Then, as $\rho\to0$:
\begin{enumerate}
\item[(i)] $\dist^2(\wstar,\calC_c)\asymp\rho^2$ and $\dist^2(\wstar,\calC)\asymp\rho^2$: a free intercept does not remove the second-order floor of the sum-to-one constraint;
\item[(ii)] there is $\rho_0>0$, depending only on the templates, such that for $0<\rho\le\rho_0$ the projection coefficients of $\wstar$ on the agents, with and without intercept, are $1+O(\rho)>0$; consequently $\dist(\wstar,\calK)=\dist(\wstar,\calH)$ and $\dist(\wstar,\calK_0)=\dist(\wstar,\calH_0)$, with $\dist^2(\wstar,\calH)=\rho^4\dist^2(e_2,\calH)\{1+O(\rho)\}$ when $e_2\notin\calH$ and $\dist^2(\wstar,\calH_0)=\rho^4\dist^2(e_2,\tilde{\calH}_0)\{1+O(\rho)\}$ when $e_2\notin\tilde{\calH}_0$.
\end{enumerate}
Nonnegativity of the weights therefore costs nothing at any order, the intercept changes only the fourth-order constant, and the sum-to-one constraint alone produces the second-order floor.
\end{thmprime}

The simplex floor of order $\rho^2$ is produced by the sum-to-one constraint alone, not by nonnegativity and not by the absence of an intercept; this is Breiman's stacking recommendation \citep{Breiman1996} and the nonnegative mixer of \citet{LiLe2024}, which has no intercept and is therefore the set $\calK_0$, inheriting the span's approximation power while keeping the cone's estimation stability, and it is why the cone column of Section~\ref{sec:real} tracks the span column under the Brier score (under the log score and AUC, where neither rule is fitted to the criterion, they part on three networks). The proof (supplement) expands the normal equations in $\rho$ and checks that the unconstrained projection is interior to the cone.

The floors are approximation error; whether the span wins at a given sample size depends on its estimation error, of order $d\rho/m$ by Theorem~S1. Setting $8d(J\rho+\supn{\wstar-\wproj}^2)/m$ against the hull floor $\rho^2(J-1)^2b^2/(4J)$ gives a crossover $m_{\times}(\rho)\asymp dJ^2/\{\rho\,(J-1)^2b^2\}$ for the least-squares synthesis ($\tau\to\infty$, the empirical setting of Section~\ref{sec:real}); at fixed $\tau$ the ridge term is negligible only when $m\tau^2\rho^2\to\infty$, since $\lambda_{\min}(G)=\Theta(\rho^2)$ under this scaling. Weaker mechanisms need more dyads before the lower floor pays, at rate $1/\rho$.

The floors above are the price of combining linearly. The superposition model itself names the family with no floor at all.

\begin{thmsecond}[Superposition: the noisy-OR family]
\label{thm:noisyor}
Let $x_j=-\log(1-w_j)$ and $\calN=\{w_\gamma=1-\exp(-\gamma_0-\sum_{j=1}^{J}\gamma_jx_j):\ \gamma_0\ge0,\ \gamma_j\ge0\}$, the noisy-OR synthesis of the agents. Under the superposition model, for every $\rho\in(0,1)$ and every library, $\wstar=w_{\gamma^\star}$ with $\gamma^\star=(0,1,\dots,1)$, so $\dist(\wstar,\calN)=0$. Moreover, for $\gamma$ in a bounded set, $w_\gamma=(1-e^{-\gamma_0})+e^{-\gamma_0}\sum_j\gamma_jw_j+O(\rho^2)$ uniformly, an element of the cone $\calK$ up to second order.
\end{thmsecond}

The simplex, the cone and the span of Theorem~$1'$ and the noisy-OR family therefore have floors of exact orders $\rho^2$ (simplex), $\rho^4\dist^2(e_2,\calH)$ (cone and span, whenever $e_2\notin\calH$) and zero (noisy-OR). These are distances to the unclipped classes: the reported predictor is $\mathrm{clip}(\bhat^{\top}\Fv)$ and $\mathrm{clip}(\calH)\not\subset\calH$, so they bound the error of the reported rule and are not floors for it. Section~S1 gives a four-cell library in which the clipped affine class has zero approximation error while the unclipped span has a positive $\rho^4$ one, so we call these unclipped approximation errors throughout, and the span is the linearization of $\calN$ at small strength: the affine rules are the right rules when the mechanisms are weak and the wrong ones when they are not, which is what Section~\ref{sec:real} finds on multiplex networks of density $0.003$ and $0.42$. The family $\calN$ is a Bernoulli regression with an additive cumulative-hazard predictor, $-\log(1-p)=\gamma_0+\sum_j\gamma_jx_j$, and not the complementary log-log model, which would put $\log\{-\log(1-p)\}$ linear in $x$. Theorem~S2 gives existence, identification of the predicted probability function, and a risk of order $d\log m/m$ for its least-squares fit over a bounded nonnegative parameter set, with the truth on the boundary $\gamma_0=0$; that is why the result is a risk bound and not a limit law, and why inference for this family is open. Identification of $\gamma$ itself needs the hazard features $1,x_1,\dots,x_J$ to be linearly independent, which does not follow from Assumption~\ref{ass:gram}: if $w_2=1-(1-w_1)^2$ then $x_2=2x_1$. What $\calN$ gives up is linearity, so it has no projection identity, no conjugate posterior, and no guarantee of being closer to a truth that is not a superposition.

\subsection{One graph: the rate depends on where the truth lies}
\label{sec:onegraph}

Every application fits the weights to dyads of one graph, and those dyads share nodes. The next theorem quantifies the sharing exactly. Under the one-graph design \eqref{eq:onegraph}, write $S=N^{-1}\sum_{i<j}\Fv_{ij}\{A_{ij}-\wproj(U_i,U_j)\}$ for the synthesis score at the projection and recall the node-level residual moment $\kappa_1(\wstar)$ and matrix $\Sigma_g$ of \eqref{eq:g}.

\begin{theorem}[One graph]
\label{thm:onegraph}
Let Assumption~\ref{ass:gram} hold and let the data follow \eqref{eq:onegraph}, with $\bhat,\Sighat$ as in \eqref{eq:posterior} over the $N$ dyads and $\kappa=\nu/(N\tau^2)$.
\begin{enumerate}
\item[(i)] For every $n\ge2$ and every graphon $\wstar$, exactly,
\[
\E\big[S^{\top}G^{-1}S\big]\;=\;\frac{c(\wstar)}{N}\;+\;\frac{4(n-2)}{n(n-1)}\,\kappa_1(\wstar).
\]
\item[(ii)] The clipped posterior mean satisfies
\[
\E\norm{\whatp-\wstar}^2\;\le\;\norm{\wstar-\wproj}^2+\frac{8\,c(\wstar)}{N}+\frac{32(n-2)}{n(n-1)}\,\kappa_1(\wstar)+\frac{8\kappa^2\norm{\bproj}^2}{\lambda}+d\exp\Big(-\frac{\lfloor n/2\rfloor}{7L}\Big).
\]
\item[(iii)] If $\kappa_1(\wstar)>0$ then $\sqrt n\,(\bhat-\bproj)\Rightarrow N\big(0,\,4\,G^{-1}\Sigma_gG^{-1}\big)$, whereas if $\wstar\in\calH$ then $\sqrt N\,(\bhat-\bproj)\Rightarrow N\big(0,\,G^{-1}MG^{-1}\big)$. In both cases $N\Sighat\to \nu G^{-1}$ in probability. Consequently, for every $\ell$ with $\ell^{\top}G^{-1}\Sigma_gG^{-1}\ell>0$, the frequentist coverage of the working-posterior credible interval for $\ell^{\top}\beta$ tends to zero as $n\to\infty$, whereas if $\wstar\in\calH$ it converges to the limit of Proposition~S2 with $N$ in place of $m$.
\end{enumerate}
\end{theorem}

The identity in (i) is the exchangeable structure of the dyads doing work: pairs sharing no node are uncorrelated by the normal equations, the $N$ diagonal terms give the independent-dyad value $c(\wstar)/N$, and the $2N(n-2)$ pairs sharing one node each contribute $\kappa_1(\wstar)$. When $\kappa_1=0$, as holds whenever $\wstar\in\calH$, the estimation error is of order $d/N\asymp n^{-2}$ and one graph is worth its number of dyads; otherwise the $\kappa_1$ term dominates and the error is of order $n^{-1}$, one graph worth its number of nodes. The dichotomy is $\kappa_1=0$ against $\kappa_1>0$, not span against non-span. Part (iii) is the inferential consequence: the working posterior spreads like $1/N$ regardless, so under misspecification its intervals are too narrow by a factor of order $\sqrt n$; the dyadic-robust covariance \eqref{eq:dyadic} estimates the shared-node sum from residuals, is consistent by Proposition~$2'$, and repairs them in Section~\ref{sec:sims}. The adjustment is for the misspecified case, not a free improvement: in the span at $n\le200$ the working posterior is the better interval, and at the sample sizes of that simulation, below a few hundred nodes, a practitioner should report both; the crossover depends on the contrast, the library and the density and is not a general threshold. When $\kappa_1=0$ with $\wstar\notin\calH$, the fast rate persists while the limit law of $\sqrt N(\bhat-\bproj)$ is a degenerate second-order object we do not characterize. The Gram term in (ii) has effective sample size $\lfloor n/2\rfloor$, the number of disjoint dyads; the proof combines a matching decomposition \citep{Hoeffding1963} with the matrix Chernoff bound \citep{Tropp2012}.

Theorem~\ref{thm:onegraph} covers known agents on one graph, the setting of the simulations; the six-network study additionally fits the agents to the same graph, and \if0\blind the companion paper \citep{PapamichalisRuane2026minimax}\else a companion paper, withheld for review,\fi\ treats estimated mechanisms by cross-fitting. The theorem is for a fixed graphon; the sparse regime is taken up in the Discussion.

The six single-layer networks of Section~\ref{sec:real} are sparse, so we record what the identity says when the truth and the agents scale with a vanishing strength.

\begin{propsparse}[Sparse regime]
\label{prop:sparse}
Let $\wstar_n=\rho_n\bar w$ and $w_{j,n}=\rho_n\tilde w_j$ with $\rho_n\in(0,1]$, $\bar w$ and $\tilde w_j$ fixed and $\tilde F=(1,\tilde w_1,\dots,\tilde w_J)$ satisfying Assumption~\ref{ass:gram} with Gram matrix $\tilde G$. Write $\tilde\ell=\tilde F^{\top}\tilde G^{-1}\tilde F$, $\tilde L=\supn{\tilde\ell}$, $\bar w^{\circ}$ for the projection of $\bar w$ onto $\tilde{\calH}=\mathrm{span}(\tilde F)$ with coefficients $\bar\beta^{\circ}$, $\bar\delta=\bar w-\bar w^{\circ}$, $\bar g(u)=\int\bar\delta(u,v)\tilde F(u,v)\,dv$ and $\bar\kappa_1=\int\bar g^{\top}\tilde G^{-1}\bar g$. Then, for every $n$ and $\rho_n$:
\begin{enumerate}
\item[(i)] $\calH_n=\tilde{\calH}$, $L_n=\tilde L$, and the Gram event $\hat G_n\succeq G_n/2$ holds if and only if it holds for $\tilde F$; and $\rho_n^2\lambda_{\min}(\tilde G)\le\lambda_{\min}(G_n)\le\rho_n^2\min_j\norm{\tilde w_j}^2$.
\item[(ii)] $w^{\circ}_n=\rho_n\bar w^{\circ}$, $\delta_n=\rho_n\bar\delta$, $\beta^{\circ}_{n,j}=\bar\beta^{\circ}_j$ for $j\ge1$ and $\beta^{\circ}_{n,0}=\rho_n\bar\beta^{\circ}_0$.
\item[(iii)] $c(\wstar_n)=\rho_n\int\bar w\,\tilde\ell+\rho_n^2\int(\bar\delta^2-\bar w^2)\,\tilde\ell$ and $\kappa_1(\wstar_n)=\rho_n^2\bar\kappa_1$, so the identity of Theorem~\ref{thm:onegraph}(i) reads
\[
\E\big[S^{\top}G_n^{-1}S\big]=\frac{\rho_n\!\int\bar w\,\tilde\ell+\rho_n^2\!\int(\bar\delta^2-\bar w^2)\,\tilde\ell}{N}+\frac{4(n-2)}{n(n-1)}\,\rho_n^2\,\bar\kappa_1 ,
\]
and when $\bar\kappa_1>0$ the ratio of the node-level term to the dyad term is $2(n-2)\rho_n\bar\kappa_1/\{\int\bar w\tilde\ell+O(\rho_n)\}\asymp n\rho_n$.
\end{enumerate}
Parts (iv) and (v), in Section~S1, give the condition under which the ridge term is negligible and the correctly scaled Gaussian limits: the slope coordinates converge at $\sqrt n$ off the span and $\sqrt{N\rho_n}$ in it, the intercept faster by $\rho_n^{-1}$.

\end{propsparse}

The proposition is algebra on the identity, but it settles where the regimes of Theorem~\ref{thm:onegraph} sit for sparse graphs: the leverage that controls the Gram term does not degenerate with the density; the $n^{-1}$ regime is exactly the regime of diverging expected degree $n\rho_n\to\infty$ (given $\bar\kappa_1>0$), with bounded expected degree the dyad and node-level terms are of one order $\rho_n/N$, and relative to $\norm{\wstar_n}^2\asymp\rho_n^2$ the two contributions are of orders $1/(N\rho_n)$, the usual sparse-graphon scale, and $1/n$. The limit laws do not carry over unscaled, because the slope coordinates of $G_n$ are of order $\rho_n^2$; part (v) gives the scaled versions, and Corollary~\ref{cor:superposition-path} the transition that the superposition model itself produces. Fitted agents remain outside the proposition, as they are outside the theorem.

\begin{corollary}[The superposition path]
\label{cor:superposition-path}
Under the superposition model with $w_j=\rho\tilde w_j$, write $r=e_2-\Pi_{\calH}e_2$ for the residual of the pairwise interaction and $K=\int\bar g_r^{\top}\tilde G^{-1}\bar g_r$ with $\bar g_r(u)=\int r(u,v)\tilde F(u,v)\,dv$. Then $\wstar-\wproj=-\rho^2r+O(\rho^3)$, $\kappa_1(\wstar)=\rho^4K+O(\rho^5)$ and $c(\wstar)=\rho\int e_1\tilde\ell+O(\rho^2)$, so the ratio of the node-level term to the dyad term in Theorem~\ref{thm:onegraph}(i) is $O(n\rho^3)$, and exactly of order $n\rho^3$ when $K>0$: along the path the theory describes the $n^{-1}$ regime is entered when $n\rho^3\to\infty$, not when $n\rho\to\infty$. The condition $K>0$ does not follow from $e_2\notin\calH$. Section~S1 gives a two-mechanism library with $r\neq0$ and $\bar g_r\equiv0$, so that $\kappa_1(\wstar)=0$ for every $\rho$ and no transition occurs.
\end{corollary}

The corollary ties the two theorems together: the strength that governs the approximation floor also governs which sampling regime the estimation error is in, through the third power. It is checked in Section~\ref{sec:sims}.

\subsection{Independent dyads: rate and minimaxity}
\label{sec:projection}

Under the independent-dyad design the synthesis is a ridge-regularized least-squares fit of outcomes on agent probabilities, and its behavior is that of the projection $\wproj$. Theorem~S1 in the supplement gives the finite-sample risk bound and the matching minimax statement: the excess risk over the approximation floor is of order $d/m$, the posterior contracts at the same rate, and no procedure improves that order over a fixed span, including procedures that see the latent positions. The rate $d/m$ over a fixed span is the classical linear-aggregation rate \citep{Tsybakov2003, BuneaTsybakovWegkamp2007}, and Theorem~S1 is its graphon instance; the three aggregation problems of \citet{Tsybakov2003}, selection, convex and linear, are the selection, hull and span trichotomy this paper uses throughout. One constant in that bound is worth carrying into the main text because it decides whether the bound is usable. The Gram-matrix term is the probability that the empirical Gram matrix drops below $G/2$, and controlling it by the leverage $L=\supn{\Fv^{\top}G^{-1}\Fv}$ instead of by $\lambda_{\min}(G)$ turns the bound from symbolic into numerical: it falls below $10^{-3}$ once $m\ge7L\log(10^{3}d)$, which is $m\ge1468$, $1605$ and $3847$ for the three simulated libraries, whereas the entrywise route needs $m\gtrsim d^2\lambda^{-2}\log(10^{3}d)$, above $10^{5}$ already at $d=3$. Proposition~$2''$(i) shows the same constant is invariant under sparsity, so the bound binds at the same $n$ however sparse the graph.

\subsection{The noisy-OR fit, choosing among rules, and fitted agents}
\label{sec:fitted}

Three results close the gaps between the theory above and the rules Section~\ref{sec:real} actually runs: an estimation theorem for the noisy-OR family, an oracle inequality for choosing a rule on held-out dyads, and a statement of what the one-graph theory says when the agents are fitted to the graph.

Theorem~S2 in the supplement gives the estimation theory for this family under the independent-dyad design: over a bounded nonnegative parameter set with the agents bounded away from one, the least-squares fit exists, the parameter is identified, the truth is the unique population minimizer at $\gamma^\star=(0,1,\dots,1)$ under the superposition model, and the risk is of order $d\log m/m$ with no approximation term. The truth sits on the boundary $\gamma_0=0$, which is why the result is a risk bound and not a limit law, and why inference for this family is open.

The parametric rate in (ii) holds because the squared loss on binary outcomes satisfies a Bernstein condition at the truth for every bounded predictor, so no convexity of the family is needed; the logarithm comes from a discretization and can be removed by a local argument we do not pursue. The truth sits on the boundary $\gamma_0=0$ of $\Gamma_B$, which is why we state a risk bound instead of a limit law for $\hat\gamma$; the estimation price is the span's $d/m$ up to the logarithm, with no approximation term. In the misspecified case (iii) the rate is the slower one that empirical risk minimization over a nonconvex class guarantees in general.

The same Bernstein condition gives a guarantee for the practitioner's actual decision, which is to choose among selection, hull, cone, span and noisy-OR on held-out dyads.

Proposition~S3 in the supplement gives the guarantee for that decision: with $K$ candidate rules fitted on data independent of $m_2$ further dyads and the rule chosen by empirical Brier score on those dyads, the excess risk of the selected rule is within a constant factor of the best rule's plus a term of order $\log(K)/m_2$.

The bound is on the raw Brier scale, where the gaps of Tables~\ref{tab:multiplex} and \ref{tab:real} are small: a gain of one percent of the density is $3\times10^{-5}$ on Amazon. With $K=5$ rules and $\delta=0.05$ the additive term is minimized near $\varepsilon=0.44$ at about $111/m_2$, so on the networks here it is larger than the gaps it would have to resolve, and it certifies only that the selected rule is not much worse than the best. What the empirical selector does is a separate question, and Section~\ref{sec:real} answers it directly by reporting the selected rule as a row: it lands within $0.5\%$ of the density of the best rule on every network.

Every agent in Section~\ref{sec:real} is fitted to the graph on which the synthesis is evaluated. The one-graph theory treats the agents as known; what survives their estimation depends on what one conditions on.

\begin{theorem}[Cross-fitted synthesis: conditional and superpopulation inference]
\label{thm:crossfit}
Let $(U,A)$ follow \eqref{eq:onegraph}. Let $\calT$ and $\calE$ be disjoint sets of dyads chosen independently of $A$ given $U$, with $|\calE|=m$, and let the agents $\hat w_1,\dots,\hat w_J$ be any measurable functions of $A_{\calT}$, the adjacency restricted to $\calT$, with values in $[0,1]$; write $\hat F_{ij}=(1,\hat w_1(i,j),\dots,\hat w_J(i,j))^{\top}$ for $(i,j)\in\calE$, $\hat G=m^{-1}\sum_{\calE}\hat F_{ij}\hat F_{ij}^{\top}$, assumed nonsingular, and let $\bhat$ be the ridge fit \eqref{eq:posterior} on $\calE$ with $\kappa=\nu/(m\tau^2)$. Define the conditional target $\beta^{\circ}_{c}=\hat G^{-1}m^{-1}\sum_{\calE}\hat F_{ij}\wstar(U_i,U_j)$, the projection of the truth on the evaluation dyads onto the fitted span, and $\hat w^{\circ}_{c}=\beta^{\circ\top}_{c}\hat F$.
\begin{enumerate}
\item[(i)] Conditionally on $(U,A_{\calT})$ the outcomes $A_{ij}$, $(i,j)\in\calE$, are independent $\mathrm{Bernoulli}\{\wstar(U_i,U_j)\}$. The normal equations defining $\beta^{\circ}_{c}$ give $\sum_{\calE}\hat F_{ij}\{\wstar(U_i,U_j)-\hat w^{\circ}_{c}(i,j)\}=0$, so the approximation residual drops out of the score,
\[
S=\frac1m\sum_{\calE}\hat F_{ij}\{A_{ij}-\wstar(U_i,U_j)\},\qquad
\E\big[S^{\top}\hat G^{-1}S\,\big|\,U,A_{\calT}\big]=\frac{\hat b}{m},\quad
\hat b=\frac1m\sum_{\calE}\wstar(1-\wstar)(U_i,U_j)\,\hat F_{ij}^{\top}\hat G^{-1}\hat F_{ij},
\]
and there is no node-level term. Only Bernoulli noise enters: conditionally on the nodes and the training dyads the misspecification of the fitted span is a fixed bias in $\beta^{\circ}_{c}$ and contributes nothing to the variance.
\item[(iii)] If, along a sequence, $\hat G\to G_{\infty}$ nonsingular, $\hat B=m^{-1}\sum_{\calE}\wstar(1-\wstar)(U_i,U_j)\hat F_{ij}\hat F_{ij}^{\top}\to B_{\infty}$, $\max_{\calE}\hat F_{ij}^{\top}\hat G^{-1}\hat F_{ij}/m\to0$ and $m\kappa\to0$, then conditionally $\sqrt m(\bhat-\beta^{\circ}_{c})\Rightarrow N(0,G_{\infty}^{-1}B_{\infty}G_{\infty}^{-1})$, and Proposition~S2 applies to the working posterior with $\hat G,\hat B$ in place of $G,M$: directional calibration holds when $\nu\hat G\succeq\hat B$ in that direction, and $\nu\ge\sup\wstar(1-\wstar)$ is sufficient. Both $\hat b$ and $\hat B$ involve the unknown $\wstar$ and are oracle quantities; Section~S1 gives the plug-in that replaces $\wstar(1-\wstar)$ by $\hat w_{c}^{\circ}(1-\hat w_{c}^{\circ})$ and the condition under which the substitution is asymptotically harmless.
\end{enumerate}
Two further parts are in Section~S1: a risk bound for the conditional target, and the law of total variance that recovers the unconditional identity of Theorem~\ref{thm:onegraph}(i) from the conditional one, the node-level moment appearing as the variance of the conditional projection under resampling of the nodes.
\end{theorem}

The theorem says what a practitioner with one graph may claim. Conditionally on the nodes she observed and on the dyads her agents were fitted to, the cross-fitted synthesis is the independent-dyad problem of Theorem~S1 with the fitted agents as known features: its target is the projection of the truth onto the span of \emph{those} agents on \emph{those} dyads, its conditional risk is $\hat b/m$ plus a ridge term, and a working posterior with $\nu$ dominating the conditional Bernoulli variance is conservative for that target. The node-level moment of Theorem~\ref{thm:onegraph} enters only if she wants inference about the population projection under resampling of the nodes, a superpopulation claim; with fitted agents that target is itself defined by the limit of the fitting procedure and is treated by cross-fitting in the companion work. The condition that the mask be independent of the outcomes is why Section~\ref{sec:real} adds a dyad-masking arm beside the positive-edge thinning protocol: the theorem covers the former and not the latter.

\subsection{Calibration of the synthesis posterior}
\label{sec:calibration}

Whether the posterior spread matches the estimation error is a separate question from the rate, because the synthesis function is a working model. The posterior covariance is $\nu G^{-1}/m$ while the sampling covariance of $\bhat$ is the sandwich $G^{-1}MG^{-1}/m$ with $M=\int\sigma^2\Fv\Fv^{\top}$, so intervals are conservative or anti-conservative according to whether $\nu$ dominates the residual second moment in the operator order. Proposition~S2 makes this precise: $\nu\ge\bar\sigma^2$ is sufficient for conservative intervals, the safe choice $\nu=1/4$ is conservative in the simulations, and the plug-in of Section~\ref{sec:real} matches the spread to the risk only when $\hat\nu$ estimates $c(\wstar)/d$, which the validation Brier score of the selected agent does not. Under the one-graph design the intervals are valid when $\kappa_1=0$ and fail in loaded directions otherwise, by Theorem~\ref{thm:onegraph}(iii). Guarantees transfer to smooth functionals of the graphon by the Lipschitz argument of Proposition~S1, for risk and posterior contraction, not for interval coverage.

\section{Simulations}
\label{sec:sims}

The simulations are designed to falsify the identities of Section~\ref{sec:theory}: the exact constants, the rate switch, and the coverage statements are computed and compared against their predicted values, and each experiment states the claim it would break. Risks are $L^2$ distances to $\wstar$ computed by quadrature on a $1000\times1000$ grid to relative error below $10^{-3}$. Section~S2 reports them in full, with four figures and six tables. Four findings carry into the rest of the paper. The superposition geometry of Theorem~\ref{thm:superposition} holds to three digits: the approximation error of the span is $\rho^4\dist^2(e_2,\calH)$ with the ratio constant to $3.00\times10^{-3}$ across $\rho$ from $0.05$ to $0.4$, while the hull floors at $\rho^2$, a risk ratio rising from $13$ to $73$. Under the one-graph design the Monte Carlo mean of the score quadratic matches the identity of Theorem~\ref{thm:onegraph}(i) at every $n$ within two standard errors, and the estimation error leaves the $n^{-2}$ line for $n^{-1}$ once the node-level moment dominates. Coverage of the working posterior falls to $0.73$--$0.77$ at $n=1600$ under a superposition truth and stays near nominal in the span, while the dyadic-robust intervals sit at $0.88$ to $1.00$ from $n=200$ on and undercover below $n\approx100$, where the working posterior is the better choice. With agents fitted to the graph the weights no longer sit near one per mechanism, because each fitted agent already absorbs part of the union, which is the regime Section~\ref{sec:real} operates in.

\section{Real networks: union of layers, pooled fits, and an evaluation artifact}
\label{sec:real}

\paragraph{Data and protocol.}
Two kinds of data are used: three multiplex networks whose observed graph is the union of separately recorded layers, the regime of Theorems~1, $1'$ and $1''$, and six single-layer networks, the standard benchmarking setting in which the mechanisms are not observed separately. The multiplexes are described with their results below and in Section~S3. The six public networks are treated as simple undirected graphs: polblogs, email-Eu-core, Cora, ca-GrQc, wiki-Vote and soc-Epinions1, with sizes and densities in Table~S13, sources and checksums in the archive and preprocessing conventions in Section~S3. Edges split at random into training ($70\%$), validation ($10\%$) and test ($20\%$); ten uniformly sampled non-edges per held-out edge, weighted so every criterion is an unbiased estimate of its uniform-dyad value. Five agents are fitted to the training graph with hyperparameters fixed in advance: Erd\H{o}s--R\'enyi \citep{ErdosRenyi1960}, Chung--Lu \citep{ChungLu2002}, a ten-bin degree-histogram blockmodel, a stochastic blockmodel \citep{Holland1983} with ten spectral blocks \citep{Sussman2012, Abbe2018}, and a rank-eight random dot product graph \citep{YoungScheinerman2007, Athreya2018}. The training graph is a Bernoulli thinning of the observed graph with retention $0.7$, so every agent fitted to it estimates $0.7$ times the observed-graph probability; we divide by the retention and clip to $[0,1]$ before any combination or evaluation. The division is exact for agents linear in edge counts and first order for the spectral agents; the device is the one \citet{LiLe2024} build into their protocol, with \citet{ElkanNoto2008} its one-dimensional ancestor. Because only positive edges are held out while sampled non-edges stay visible as zeros, the design is positive-edge thinning, not missing-at-random dyad masking; treating held-out dyads as missing, as in \citet{LiLevinaZhu2020}, would remove the attenuation at the source. Weights are learned on the validation dyads, the constant agent absorbed by the intercept, with $\tau^2=100$ and $\nu$ the validation Brier score of the selected agent. Comparators: selection, the calibrated selected agent (affine and Platt), hull synthesis under both scores, the nonnegative cone, unpenalized least squares, logistic stacking, and the rank-sixteen random dot product graph, the enlarge-instead-of-combine alternative \citep{Hoff2002}. Scores are the Brier score over the density, the logarithmic score over the entropy of the constant predictor, and the AUC. Paired comparisons use $t$ intervals on the ten split differences, beside the dependence-adjusted intervals of \citet{NadeauBengio2003}.

\paragraph{Union of layers.}
The regime of Theorems~$1'$ and $1''$, a graph that is the union of separately observed mechanisms, is instantiated by multiplex networks: the observed graph is an edge wherever at least one layer connects the pair, and an agent fitted to one layer carries that layer's information only. We use three public multiplexes (Table~\ref{tab:multiplex}, sources and preprocessing in Section~S3): CS-Aarhus, $61$ people and five relations; Amazon, $10{,}099$ products and two relations; and YouTube, $2{,}000$ users and five relations, whose union has density $0.42$. The protocol is that of the six-network study on the union edges, every layer restricted to the training union edges and de-thinned by its own retention. Libraries are the five union agents; one agent per layer, its class chosen on the validation dyads; the same with a rank-eight spectral agent on every layer, which is the bagging control; and the union agents, the rank-sixteen union model and the layer agents together. Beside the rules of Table~\ref{tab:real} we fit the structural rule of Theorem~$1''$, the noisy-OR synthesis with $\gamma\ge0$ by weighted maximum likelihood on the validation dyads. The single models fitted to the union are those a statistician would reach for first, each tuned on the same held-out dyads as the synthesis: singular value thresholding \citep{Chatterjee2015} with the rank chosen on the validation dyads (Section~\ref{sec:artifact}), neighborhood smoothing \citep{ZhangLevinaZhu2017} where $n\le8{,}000$, a degree-corrected blockmodel \citep{KarrerNewman2011} with $K$ by edge cross-validation \citep{LiLevinaZhu2020}, the multilayer degree-corrected blockmodel \citep{HanXuAiroldi2015, PaulChen2020}, and the calibrated heuristics, topological logistic and random forest of Table~S8. Every single model is Platt-calibrated on the validation dyads and also enters the library.

Three things happen, in order of how much they favor the synthesis. First, the predictions of Section~\ref{sec:theory} hold on all three networks. The span weights on the layer agents are $0.97$ and $1.03$ on Amazon, near one per mechanism, and $0.84$, $0.71$, $0.12$, $0.15$, $0.54$ on CS-Aarhus, unequal, with a sum of $2.4$ against the value one that a sum-to-one rule would impose; the noisy-OR strengths are $1.00$ and $0.98$ on Amazon and $0.99$, $1.02$, $1.00$, $0.82$, $0.87$ on YouTube; the span beats the hull on the layer library by $3.6$ $[3.4,3.9]$, $9.0$ $[6.4,11.5]$ and $16.4$ $[16.4,16.4]\%$ of the density on Amazon, CS-Aarhus and YouTube, the largest span--hull gaps in the paper; and the per-layer span beats the pooled calibrated agent by $2.6$ $[2.4,2.9]$ on Amazon and $7.7$ $[3.3,12.1]$ on CS-Aarhus. On YouTube (density $0.42$) the affine span of the layer agents loses to the rank-sixteen model by $3.6$ while the noisy-OR synthesis beats it by $1.1$ $[1.1,1.1]$: at layer densities of $0.08$ to $0.25$ the interactions that Theorem~$1'$ treats as second order are not small and the exact link of Theorem~$1''$ is needed. The bagging control is decisive only on CS-Aarhus, where the diverse library beats the homogeneous one by $3.1$ $[0.4,5.9]$; on Amazon and YouTube validation selects the spectral class on every layer and the two libraries coincide, so the gain there is layer information, not class diversity. Second, against single models of fixed size the synthesis wins: on Amazon the full library beats the rank-sixteen union model by $1.05$ $[0.87,1.22]$ ($[0.72,1.37]$ dependence-adjusted), on YouTube its noisy-OR form beats it by $1.3$, and on CS-Aarhus the full-library syntheses beat it by $19$ to $26$ and the layer hull by $8$. Third, and this is the comparison that decides whether combining is worth anything, the single models are tuned as carefully as the synthesis is. The rank of the thresholded spectral estimator is chosen on the validation dyads instead of by \citeauthor{Chatterjee2015}'s threshold, for the reason given in Section~\ref{sec:artifact}; the blockmodel's $K$ is chosen by edge cross-validation; every heuristic and learner is calibrated on the same dyads. Against that comparator set the synthesis still wins. On Amazon the best single model is the spectral estimator at $0.709$ of the density, ahead of the multilayer blockmodel and Adamic--Adar (both $0.858$), the blockmodel ($0.878$) and the rank-sixteen model ($0.834$); the library containing all of them synthesizes to $0.698$ under the span and $0.702$ under the cone, gains of $1.06$ $[0.56,1.57]$ and $0.67$ $[0.20,1.14]$, the first still positive under the node bootstrap $[0.12,2.01]$. On CS-Aarhus the best single model is the seven-feature topological logistic at $0.491$, and the same library gives $0.474$ under the cone and $0.482$ under the span, gains of $1.68$ $[-2.22,5.58]$ and $0.89$ $[-3.81,5.59]$. Both cover zero under the paired intervals and under the node bootstrap, so CS-Aarhus supports no claim of improvement; at $61$ nodes it could not. On YouTube the spectral estimator at rank $30$ is best at $0.180$, the rank-sixteen model reaches $0.192$, neighborhood smoothing $0.241$ and the multilayer blockmodel $0.304$; the noisy-OR synthesis of the library containing them scores $0.173$, a gain of $0.70$ $[0.68,0.73]$, $[0.65,0.76]$ under the node bootstrap, while the affine span of the same library only edges it, $0.09$ $[0.07,0.11]$. The affine rule is not enough at density $0.42$ and the exact link is; that is the empirical content of Theorem~$1''$. The noisy-OR rule is beaten by the affine span on Amazon, by $2.56$ $[2.27,2.85]$ of the density, because the two layers there are nearly disjoint and the union link buys nothing, so the ordering of the rules itself depends on the regime, which is why Section~\ref{sec:fitted} treats the choice among them as the estimand and Proposition~S3 bounds its cost.

\begin{table}[htbp]
\centering
\caption{Union-of-layers multiplex networks: Brier score over the union density, mean (s.e.) over ten splits, best per column in bold. Union agents: the five agents of Table~\ref{tab:real} fitted to the union; rank-16: the Platt-calibrated rank-sixteen spectral model of the union; layers: one agent per layer, either a rank-eight spectral agent on every layer (the bagging control) or the class selected on the validation dyads; the rows that follow combine the union agents, the rank-sixteen union model and the layer agents, then list the strong single models fitted to the union (thresholded spectral estimation with the rank chosen on the validation dyads, neighborhood smoothing where $n\le8{,}000$, DC-SBM with $K$ by edge cross-validation, the multilayer DC-SBM, the calibrated heuristics and the nested-calibrated forest), every one Platt-calibrated, and the synthesis of the full library augmented with them. Ten splits throughout. All rules and criteria are in \texttt{results/multiplex\_summary.txt} of the archive.}
\label{tab:multiplex}
\footnotesize\setlength{\tabcolsep}{4pt}
\begin{tabular}{lccc}
\toprule
 & CS-Aarhus & Amazon & YouTube \\
 & $n=61$, $L=5$ & $n=10,099$, $L=2$ & $n=2,000$, $L=5$ \\
 & density $0.193$ & density $0.003$ & density $0.418$ \\
\midrule
union agents: selected, calibrated & 0.650 (0.016) & 0.878 (0.001) & 0.224 (0.000) \\
union agents: span & 0.566 (0.014) & 0.877 (0.001) & 0.224 (0.000) \\
union rank-16, calibrated & 0.744 (0.009) & 0.834 (0.001) & 0.192 (0.000) \\
layers, rank-8 each: span (bagging control) & 0.604 (0.010) & 0.851 (0.002) & 0.228 (0.000) \\
layers, selected classes: hull & 0.663 (0.015) & 0.888 (0.001) & 0.392 (0.000) \\
layers, selected classes: span & 0.573 (0.012) & 0.851 (0.002) & 0.228 (0.000) \\
layers, selected classes: noisy-OR & 0.579 (0.011) & 0.851 (0.002) & 0.181 (0.000) \\
union $+$ layers: span & 0.518 (0.010) & 0.848 (0.001) & 0.210 (0.000) \\
union $+$ rank-16 $+$ layers: span & 0.518 (0.010) & 0.823 (0.001) & 0.191 (0.000) \\
union $+$ rank-16 $+$ layers: cone & 0.512 (0.011) & 0.825 (0.001) & 0.191 (0.000) \\
union $+$ rank-16 $+$ layers: noisy-OR & 0.552 (0.012) & 0.839 (0.003) & 0.179 (0.000) \\
spectral, rank by validation, calibrated & 0.620 (0.020) & 0.709 (0.002) & 0.180 (0.000) \\
neighborhood smoothing, calibrated & 0.674 (0.013) & -- & 0.241 (0.000) \\
DC-SBM ($K$ by edge CV), calibrated & 0.593 (0.017) & 0.878 (0.002) & 0.337 (0.000) \\
multilayer DC-SBM, calibrated & 0.569 (0.016) & 0.858 (0.002) & 0.304 (0.001) \\
topological logistic (7 features) & 0.491 (0.012) & 0.871 (0.003) & 0.235 (0.000) \\
random forest, nested-calibrated & 0.549 (0.008) & 0.751 (0.002) & 0.204 (0.000) \\
full $+$ strong models: span & 0.490 (0.018) & 0.699 (0.003) & 0.179 (0.000) \\
full $+$ strong models: noisy-OR & 0.507 (0.018) & 0.736 (0.001) & \textbf{0.173 (0.000)} \\
full $+$ strong $+$ multilayer: span & 0.482 (0.019) & \textbf{0.698 (0.003)} & 0.179 (0.000) \\
full $+$ strong $+$ multilayer: cone & \textbf{0.474 (0.018)} & 0.702 (0.003) & 0.179 (0.000) \\
full $+$ strong $+$ multilayer: noisy-OR & 0.497 (0.019) & 0.734 (0.002) & \textbf{0.173 (0.000)} \\
\bottomrule
\end{tabular}

\end{table}

\begin{figure}[!b]
\centering
\includegraphics[width=\textwidth]{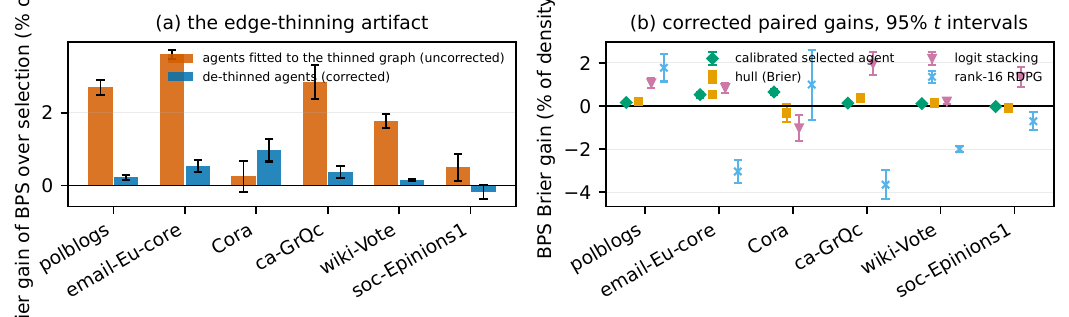}
\caption{Six networks. (a) Brier gain of span synthesis over selection with agents fitted to the thinned graph and with de-thinned agents: the uncorrected protocol manufactures most of the gain. (b) Corrected paired gains of span synthesis over the calibrated selected agent, hull synthesis, logistic stacking and the rank-sixteen random dot product graph; $95\%$ intervals over ten splits.}
\label{fig:real}
\end{figure}

\paragraph{Pooled fits on six networks: when the premise fails.}
Table~\ref{tab:real} reports all rules on de-thinned agents. Under the Brier score the span synthesis beats the affine-calibrated selected agent on five networks, by $0.16$ $[0.10,0.22]$, $0.52$ $[0.38,0.67]$, $0.64$ $[0.52,0.77]$, $0.13$ $[0.11,0.16]$ and $0.11$ $[0.09,0.13]$ percent of the density on polblogs, email-Eu-core, Cora, ca-GrQc and wiki-Vote, and loses nothing on soc-Epinions1 ($-0.03$ $[-0.07,0.01]$); intervals are paired $t$ over ten splits (Figure~\ref{fig:real}(b)). Ten splits of one graph are dependent, so we also report dependence-adjusted intervals \citep{NadeauBengio2003}, which widen every margin and still exclude zero on the five (Section~S3). Under Platt calibration the selected agent edges the synthesis on Cora ($0.995$ against $0.997$ of density), so the five-network claim is specific to the affine calibration. The margins are small because, after correction, the projection sits close to the selected agent, so the span--hull gap of Theorem~\ref{thm:superposition}(i), the most synthesis can gain, is small (Table~S7). Enlarging beats combining on four networks: the rank-sixteen random dot product graph is best under the Brier score on email-Eu-core, ca-GrQc, wiki-Vote and soc-Epinions1, by $3.1$, $3.7$, $2.0$ and $0.7\%$ of the density, while combining wins on polblogs and Cora. On Cora, where the weights are genuinely multi-agent and unstable across splits, the hull beats the span ($-0.34$ $[-0.70,0.03]$, covering zero) and logistic stacking beats both; the hypothesis, that weak collinear agents make the span's extra freedom cost more than it buys, is testable in Figure~S1(b). Under the logarithmic score the span synthesis is never best: the winner is logistic stacking on email-Eu-core, Cora and wiki-Vote, the hull refitted under that score on polblogs and soc-Epinions1, and the cone on ca-GrQc. Selection fails outright on ca-GrQc and Cora (log scores $1.17$ and $1.35$ times the constant predictor) because the selected agent puts probabilities below $10^{-3}$ on most held-out edges; the intercept repairs this ($0.73$ and $0.92$), the linear-pool behavior \citet{RanjanGneiting2010} describe. Under AUC logistic stacking is best on all six (Table~S12). Across the twelve Brier and log-score cells of Table~\ref{tab:real} the synthesis is best in one, polblogs under Brier, where it is the same fit as unpenalized least squares to six decimals. A single larger model, the rank-sixteen spectral agent, is best in four of the six Brier cells, and logistic stacking is best in four of the twelve cells overall. In none of the six AUC cells of Table~S12 is the synthesis best. The least-squares column equals the synthesis in every cell at the displayed precision, since at $\tau^2=100$ the ridge $\kappa=\nu/(m\tau^2)$ is of order $10^{-7}$ and the prior is numerically inert; those two columns are the nonnegative and OLS mixers of \citet{LiLe2024} on this library, and their exponential mixer coincides with selection to three decimals except on Cora, where it gains $0.73\%$. Prior sensitivity and clip rates are in Section~S3.

\begin{table}[htbp]
\centering
\caption{Test scores on de-thinned agents, mean (s.e.) over ten splits; best per row in bold: Brier over density and log score over the entropy of the constant predictor; AUC is in Table~S12.}
\label{tab:real}
\resizebox{\textwidth}{!}{\begin{tabular}{llcccccccccc}
\toprule
Network & score & Select & Select+affine & Hull (Brier) & Hull (log) & Cone & Least squares & Exp.\ mix & Logit stacking & BPS (span) & RDPG-16 \\
\midrule
polblogs & Brier$/e(\mathbf{A})$ & 0.702 (0.002) & 0.701 (0.002) & 0.702 (0.002) & 0.703 (0.002) & \textbf{0.700 (0.002)} & \textbf{0.700 (0.002)} & 0.702 (0.002) & 0.710 (0.002) & \textbf{0.700 (0.002)} & 0.717 (0.004) \\
 & log score$/H_0$ & 0.662 (0.003) & 0.670 (0.004) & 0.609 (0.007) & \textbf{0.587 (0.002)} & 0.634 (0.004) & 0.638 (0.004) & 0.662 (0.003) & 0.590 (0.002) & 0.638 (0.004) & 0.746 (0.004) \\
email-Eu-core & Brier$/e(\mathbf{A})$ & 0.673 (0.002) & 0.673 (0.002) & 0.673 (0.002) & 0.675 (0.002) & 0.668 (0.002) & 0.668 (0.002) & 0.673 (0.002) & 0.676 (0.002) & 0.668 (0.002) & \textbf{0.637 (0.003)} \\
 & log score$/H_0$ & 0.659 (0.004) & 0.666 (0.005) & 0.596 (0.003) & 0.589 (0.002) & 0.625 (0.004) & 0.625 (0.004) & 0.659 (0.004) & \textbf{0.586 (0.002)} & 0.625 (0.004) & 0.666 (0.006) \\
Cora & Brier$/e(\mathbf{A})$ & 1.007 (0.004) & 1.004 (0.004) & 0.994 (0.002) & 0.992 (0.000) & 0.997 (0.003) & 0.997 (0.003) & 1.000 (0.003) & \textbf{0.987 (0.002)} & 0.997 (0.003) & 1.007 (0.008) \\
 & log score$/H_0$ & 1.353 (0.066) & 0.966 (0.008) & 0.972 (0.008) & 0.893 (0.001) & 0.997 (0.020) & 0.917 (0.006) & 1.065 (0.045) & \textbf{0.884 (0.002)} & 0.917 (0.006) & 1.329 (0.007) \\
ca-GrQc & Brier$/e(\mathbf{A})$ & 0.772 (0.003) & 0.770 (0.003) & 0.772 (0.003) & 0.827 (0.001) & 0.769 (0.003) & 0.768 (0.003) & 0.772 (0.003) & 0.788 (0.003) & 0.768 (0.003) & \textbf{0.732 (0.004)} \\
 & log score$/H_0$ & 1.169 (0.004) & 0.765 (0.002) & 1.169 (0.004) & 0.742 (0.002) & \textbf{0.720 (0.002)} & 0.726 (0.004) & 1.169 (0.004) & 0.722 (0.002) & 0.726 (0.004) & 0.960 (0.010) \\
wiki-Vote & Brier$/e(\mathbf{A})$ & 0.802 (0.002) & 0.802 (0.002) & 0.802 (0.002) & 0.804 (0.001) & 0.802 (0.002) & 0.801 (0.002) & 0.802 (0.002) & 0.803 (0.002) & 0.801 (0.002) & \textbf{0.781 (0.002)} \\
 & log score$/H_0$ & 0.595 (0.002) & 0.594 (0.001) & 0.575 (0.004) & 0.551 (0.001) & 0.592 (0.002) & 0.596 (0.002) & 0.595 (0.002) & \textbf{0.549 (0.001)} & 0.596 (0.002) & 0.598 (0.002) \\
soc-Epinions1 & Brier$/e(\mathbf{A})$ & 0.925 (0.002) & 0.926 (0.003) & 0.926 (0.002) & 0.932 (0.001) & 0.926 (0.003) & 0.926 (0.003) & 0.925 (0.002) & 0.940 (0.002) & 0.926 (0.003) & \textbf{0.919 (0.002)} \\
 & log score$/H_0$ & 0.702 (0.001) & 0.702 (0.005) & 0.681 (0.006) & \textbf{0.640 (0.001)} & 0.685 (0.004) & 0.701 (0.005) & 0.702 (0.001) & 0.645 (0.001) & 0.701 (0.005) & 0.716 (0.002) \\
\bottomrule
\end{tabular}
}
\end{table}

A related experiment in which the five agents hold disjoint private shares of one graph is in Section~S3: its own bagging control attributes most of the gain to pooling of the shares, so we do not count it as evidence for mechanism diversity.

The corrected weights, posterior spreads and the recovery of triangle densities and clustering coefficients are reported in Section~S3 (Tables~S3 and S7): the dominant spectral agent receives $0.90$ to $1.10$ on five networks and $0.65$ on Cora, the intercepts are zero to two decimals, and the weight posterior is up to fifteen times wider on Cora, where two agents are nearly collinear.

\subsection{Two ways a benchmark can mislead}
\label{sec:artifact}

\paragraph{The edge-thinning artifact: a methodological result.}
Figure~\ref{fig:real}(a) and Table~S4 show the same experiment with and without the retention correction. Without it, the span synthesis beats selection by $2.7$, $3.6$, $2.8$, $1.8$ and $0.5\%$ of the density on the five non-Cora networks, the weight on the dominant random dot product agent is $1.25$ to $1.52$, near $1/0.7=1.43$, and the weight sums are $1.47$ to $2.07$: the fitted rule is mostly multiplying the attenuated best agent back up, which selection and the hull cannot do. With the correction the dominant weight returns to $0.90$--$1.10$ on those five networks ($0.65$ on Cora), the recalibration slope of the selected agent is $0.86$--$1.11$, and the gains shrink by factors of $6.8$ to $12.5$ on four networks, while Cora's grows from $0.25$ to $0.97$ and soc-Epinions1's changes sign, from $0.50$ to $-0.18$. An uncorrected edge-holdout benchmark therefore gives any free-weight combiner an advantage unrelated to combining, and the triangle densities it appears to recover are the cube of the same rescaling. Four rules follow for anyone benchmarking network model combination: divide by the retention before learning free weights, score gains against a calibrated selected agent, attach dependence-aware uncertainty, and report the clip rate (Section~S3).

\paragraph{Tuning the comparator: the threshold that keeps everything.} The second way cuts against the paper's own thesis. A combination rule is interesting only if the single models it is compared against are tuned with the same care, and the obvious way to fit the strongest single model is the estimator's published default. For singular value thresholding that default keeps the components above $2.01\sqrt{n\hat\rho}$ \citep{Chatterjee2015}, a threshold derived under a dense graphon with $\hat\rho$ the fraction of observed entries; our implementation uses the edge density, which is the reading appropriate to a fully observed graph and is not a uniform variance bound for a heterogeneous network. Read that way it does not bind on sparse graphs: it retains all $256$ leading components on Amazon, wiki-Vote and soc-Epinions1, and $231$ to $240$ on ca-GrQc. On wiki-Vote the resulting estimator scores $0.967$ of the density against $0.996$ for the constant predictor, recovering almost none of the available structure, while a rank chosen on the validation dyads keeps $20$ to $30$ components and scores $0.777$; on Amazon the switch moves the score from $0.726$ to $0.709$ and on ca-GrQc from $0.820$ to $0.794$. Every rank of a \emph{tuned comparator} in this paper is therefore chosen on held-out dyads, as the blockmodel's number of communities already was, while the rank-eight and rank-sixteen library agents stay fixed by construction. The effect runs against us: it strengthens the competitor and shrinks the gains of Section~\ref{sec:real} by a factor of two to three. A combination result stated against an untuned single model measures the handicap, in the same way and for the same reason that a gain measured on a thinned holdout does.

\paragraph{External comparators.}
Section~S3 and Table~S8 place six external link predictors beside the synthesis rules on the identical held-out dyads: calibrated heuristics, a seven-feature topological logistic and a random forest, raw and nested-calibrated. Under the uniform-dyad Brier score the graphon rules stay best on three networks; the nested-calibrated forest wins Cora, the raw forest wins email-Eu-core, and on ca-GrQc common neighbors, Adamic--Adar and the calibrated topological learners beat every model-based rule, the triadic information lying outside the agent span. The forest is the best ranker on four networks, and its uncalibrated probabilities score seven times the density under the Brier score on soc-Epinions1, a calibration failure that nested calibration repairs at unchanged AUC. Section~S3 adds a learned embedding predictor, node2vec with a logistic link on the Hadamard product of the endpoint vectors, on all nine networks: it is worse than every model-based rule under the Brier score on all nine, and best of all single models under AUC on two, which is the same separation between ranking and probability that the forest shows.

\section{Discussion}
\label{sec:discussion}

The positive claims are three. When mechanisms superpose, the simplex is the wrong constraint set at first order: floors of exact orders $\rho^2$ against $\rho^4$, with constant the distance from the pairwise interaction to the span, and the nonnegative cone on the fast side of that divide. From one graph, the estimation error obeys an exact variance identity whose node-level moment, zero in the span and possibly zero outside it, switches the rate from $n^{-2}$ to $n^{-1}$ and breaks the working posterior in loaded directions, with a dyadic-robust adjustment proved consistent and verified in simulation. For independent dyads the synthesis is minimax rate-optimal over the span. The negative claim is empirical: once the edge-thinning artifact is corrected and the selected agent is calibrated, the residual value of combining five fitted models is $0.11$ to $0.64\%$ of the density under the Brier score, and a larger single model beats the combination on four of the six networks; the uncorrected protocol overstates the gain by an order of magnitude on four, understates it on Cora, and flips its sign on soc-Epinions1. The synthesis pays when the span--hull gap is material: mechanisms of comparable strength that superpose, not one dominant model with satellites; with disjoint data the gain is $1$ to $6\%$ of the density on five of six (Table~S10), though a bagging control attributes most of it to pooling of the shares and not to mechanism diversity (Table~S11). The union-of-layers experiment is where the superposition model is applied to real data: the weights recover the union on the sparse multiplex, the span beats the fitted sum-to-one rules by an order of magnitude more than the pooled gains, the noisy-OR form of Theorem~$1''$ is needed once the layers are dense, and single models of fixed size lose to the synthesis. The comparison that matters is against single models tuned as carefully as the synthesis: with the spectral rank chosen on held-out dyads, the blockmodel's $K$ by edge cross-validation and every learner calibrated, the synthesis still wins on all three multiplexes, by $1.06$ and $0.70\%$ of the density on Amazon and YouTube, with no improvement established on CS-Aarhus. That the margins survive a fairly tuned comparator is what makes numbers of this size worth reporting; the pooled six-network study, where the mechanisms are not recorded separately, is where they do not.

Four limitations remain. The theory takes the agents as known while the six-network study fits them to the same graph; \if0\blind cross-fitting is treated in \citet{PapamichalisRuane2026minimax} and dynamics in \citet{PapamichalisRuaneEtAl2026dynamic}.\else cross-fitting and dynamics are treated in two companion papers withheld for review.\fi{} The theory is Brier-score theory, and the synthesis is never best under the log score; a logit-scale synthesis with the analogue of Theorem~\ref{thm:onegraph} is the natural next step. Agents are assumed jointly aligned; combining independently estimated graphons needs alignment over measure-preserving maps. And the theory is dense-graphon theory, with sparse joint asymptotics open.

\setlength{\bibsep}{0pt plus 0.3ex}
\bibliographystyle{agsm-jasa}
\bibliography{refs}

\bigskip
\if0\blind

\clearpage
\setcounter{section}{0}\setcounter{table}{0}\setcounter{figure}{0}\setcounter{equation}{0}
\renewcommand{\thesection}{S\arabic{section}}
\renewcommand{\thetable}{S\arabic{table}}
\renewcommand{\thefigure}{S\arabic{figure}}
\renewcommand{\theequation}{S\arabic{equation}}
\begin{center}
{\LARGE\bf Supplement to ``Graphon Synthesis of Network Mechanisms''}
\end{center}
\medskip
\spacingset{1.5}

This supplement contains the proofs of Theorems 1--4 and 6, Theorems $1'$ and $1''$, Propositions $2'$, $2''$ and 5, and Corollary 1 of the main text (Section~S1), the calibration proposition and the Lipschitz transfer to network functionals (Propositions~S1 and S2 in Section~S1), and additional tables for the experiments (Section~S2). Numbering without the prefix S refers to the main text; all notation is that of Section~2 there. In the main text, display (1) defines the node-level residual moments, (2) is the independent-dyad design, (3) the one-graph design, (5) the posterior and (6) the dyadic-robust covariance.

\section{Simulation detail}
\label{sec:Ssims}

\paragraph{Superposition.}
The geometry of Theorem~1 of the main text is verified in Figure~S2 with $J=2$, a two-block assortative graphon and a rank-one degree graphon: the span floor divided by $\rho^4$ equals $\dist^2(\tilde w_1\tilde w_2,\calH)=3.00\times10^{-3}$ at every $\rho$, the exact identity of part (iii); the hull floor tracks slope two between its bounds; the best hull element at $\rho=0.4$ drops the degree mechanism entirely, with error about $180$ times the span floor; and at $m=5000$ the hull-to-span risk ratio grows from $13$ at $\rho=0.05$ to $73$ at $\rho=0.4$, the crossover of Section~3.1 of the main text read at fixed $m$.

\paragraph{One graph.} Figure~S4 runs the one-graph design with the $d=3$ library at $n=50$ to $3200$, $R=60$ replications, two truths: one in the span, one the superposition of two mechanisms. The Monte Carlo mean of the score quadratic matches the identity of Theorem~2 of the main text(i) at every $n$ within two standard errors. In the span the estimation error follows $n^{-2}$; under the superposition truth it leaves that line for $n^{-1}$ once the node-level term dominates, which at these sizes happens near the right edge of the plotted range, so the plot shows the crossover and not a long $n^{-1}$ stretch. Coverage of the working posterior at three probed link probabilities, the target being $w^{\circ}$, falls from nominal to $0.73$--$0.77$ at $n=1600$ for the superposition truth and stays near nominal in the span; the dyadic-robust intervals sit at $0.88$ to $1.00$ for both truths from $n=200$ on, with Monte Carlo standard error $0.03$ at nominal and $0.06$ at $0.73$, and undercover at $n\le100$, where the working posterior is the better choice. Computational limits are stated in Section~S3.

\paragraph{Truth in the span.} Three libraries of size $J=2,4,8$ ($d=3,5,9$), built from block, degree and smooth kernels with the truth inside the span, test the rate of Theorem~S1 and its minimax part (Figure~S3, Table~S6). The estimation error tracks $c(\wstar)/m$ with no free constant across $m=50$ to $25{,}600$, the product $m\times\text{risk}$ is flat in $m$ at $0.15$ to $0.21$ times $d$, and the span overtakes the hull at the first grid point where the hull's floor exceeds $c(\wstar)/m$. The full bound of Theorem~S1(i) sits a factor of eight above the realized risk at the largest $m$, and the lower bound is below it by the stated factor of about $400$, so the bound is informative about the rate and not about the constant.

\paragraph{Misspecification, collinearity and the fitted-agent pipeline.}
Three further experiments are in Figure~S1 and Section~S3: a truth outside the span, where the synthesis converges to the approximation floor while the hull stops a factor $1.5$ and selection a factor $19$ above it; two nearly collinear agents, where least squares reaches coefficient error $51$ while the posterior mean at $\tau^2=1$ has coefficient error $0.02$ and posterior standard deviation near one, so unidentified weights are reported as unidentified without harming the synthesized graphon; and the full fitted-agent pipeline on single graphs from a three-mechanism superposition, where combining beats calibrating the best agent by about a fifth of the remaining risk and the rank-sixteen model does not substitute for it. There the fitted agents already absorb part of the union, so the weights sum toward one instead of toward $J$: the near-unit-per-mechanism regime of Theorem~1 of the main text concerns population mechanisms.

\section{The noisy-OR fit}
\label{sec:Shazard}

\begin{theoremS}[Least-squares fit of the noisy-OR family]
\label{thm:hazardS}
In the independent-dyad design the independent-dyad design (2) of the main text, let $\sup_j\supn{w_j}\le1-\epsilon$ for some $\epsilon>0$, let $x_j=-\log(1-w_j)$, let $(1,x_1,\dots,x_J)$ be linearly independent in $L^2([0,1]^2)$, and let $\Gamma_B=\{\gamma\in\R^{d}:\gamma\ge0,\ \norm{\gamma}_1\le B\}$ with $B\ge J$. Let $\hat\gamma$ minimize $m^{-1}\sum_s\{Y_s-w_\gamma(U_s,V_s)\}^2$ over $\Gamma_B$.
\begin{enumerate}
\item[(i)] A minimizer exists, and $w_\gamma=w_{\gamma'}$ almost everywhere implies $\gamma=\gamma'$. Under the superposition model of Theorem~1 of the main text the population risk over $\Gamma_B$ is uniquely minimized at $\gamma^\star=(0,1,\dots,1)$, where it equals the Bayes risk.
\item[(ii)] Under the superposition model, for every $m\ge3$,
\[
\E\norm{w_{\hat\gamma}-\wstar}^2\le C_{B,\epsilon}\,\frac{d\log m}{m},\qquad C_{B,\epsilon}=80\,\{1+\log(3BL_\epsilon)\},\quad L_\epsilon=1+\log(1/\epsilon).
\]
\item[(iii)] Without the superposition model, for every $m\ge5$, $\E\norm{w_{\hat\gamma}-\wstar}^2\le\dist^2(\wstar,\{w_\gamma:\gamma\in\Gamma_B\})+4\sqrt{1+\log(3BL_\epsilon)}\,\sqrt{d\log m/m}$.
\end{enumerate}
\end{theoremS}

\section{Rate and minimaxity under the independent-dyad design}
\label{sec:S0}

Under the independent-dyad design the synthesis is a ridge-regularized least-squares fit of outcomes on agent probabilities, and its behavior is that of the projection $\wproj$. Let $\calW_{\calH}=\{w\in\calH:0\le w\le1\}$, let $\calH_0=\{h-\int h:h\in\calH\}$, and let $B=\inf\sum_{k=1}^{d-1}\supn{\psi_k}$ over orthonormal bases $(\psi_k)$ of $\calH_0$.

\begin{theoremS}[Rate and minimaxity for independent dyads]
\label{thm:minimaxS}
Let Assumption~1 of the main text hold and let the data follow its independent-dyad design (2). For every graphon $\wstar$, every $m\ge1$ and every $\nu,\tau^2>0$, with $\kappa=\nu/(m\tau^2)$:
\begin{enumerate}
\item[(i)]
\begin{equation}
\label{eq:rate}
\E\norm{\whatp-\wstar}^2\;\le\;\norm{\wstar-\wproj}^2+\frac{8\,c(\wstar)}{m}+\frac{8\kappa^2\norm{\bproj}^2}{\lambda}+d\exp\Big(-\frac{m}{7L}\Big),
\end{equation}
and the posterior contracts at the same rate: for every $r>0$, $\E\,\Pi_m(\{w:\norm{w-\wproj}^2>r\})\le r^{-1}\{(8c(\wstar)+2\nu d)/m+8\kappa^2\norm{\bproj}^2/\lambda\}+d\exp\{-m/(7L)\}$.
\item[(ii)] For every $m\ge3B^2/8$,
\[
\frac{3(d-1)}{512\,m}\le\inf_{\tilde w}\sup_{\wstar\in\calW_{\calH}}\E\norm{\tilde w-\wstar}^2\le\sup_{\wstar\in\calW_{\calH}}\E\norm{\whatp-\wstar}^2\le\frac{2d}{m}+\frac{8\nu^2}{m^2\tau^4\lambda^2}+d\exp\Big(-\frac{m}{7L}\Big),
\]
the infimum over all estimators measurable in $(U_s,V_s,Y_s)_{s\le m}$.
\end{enumerate}
\end{theoremS}

The Gram-matrix term is the probability that the empirical Gram matrix drops below $G/2$, and controlling it by the leverage $L$ instead of by $\lambda$ turns the bound from symbolic into numerical: it is below $10^{-3}$ once $m\ge7L\log(10^{3}d)$, which is $m\ge1468$, $1605$ and $3847$ for the three simulated libraries, whereas the entrywise Hoeffding route needs $m\gtrsim d^2\lambda^{-2}\log(10^{3}d)$, above $10^{5}$ already at $d=3$. At $m=25{,}600$ the full bound is $1.6\times10^{-4}$, $3.0\times10^{-4}$ and $4.9\times10^{-4}$ for $d=3,5,9$, a factor of eight above the realized risk, while $c(\wstar)/m$ matches the realized risk with no free constant. Part (ii) says no procedure improves the order $d/m$ over the span, including ones that see the latent positions; its two constants differ by a factor of about $400$. The rate $d/m$ over a fixed span is the classical linear-aggregation rate \citep{Tsybakov2003, BuneaTsybakovWegkamp2007}; the content here is the graphon instance with explicit constants. The lower bound concerns the worst case over $\calW_{\calH}$ and is displayed only where $m\ge3B^2/8$; orthonormalized bases give $B\le6.7$, $13.6$ and $30.3$, so it holds from $m\ge17$, $69$ and $345$.

\section{Proofs}
\label{sec:proofs}

Throughout, Assumption~1 holds. Under the independent-dyad design (2) write
\[
\xi_s=Y_s-\wproj(U_s,V_s),\qquad S=\frac1m\sum_{s=1}^m\Fv_s\xi_s,\qquad \sigma^2(u,v)=\wstar(1-\wstar)(u,v)+(\wstar-\wproj)^2(u,v);
\]
under the one-graph design (3) write $\xi_{ij}=A_{ij}-\wproj(U_i,U_j)$, $\delta=\wstar-\wproj$, $h(u,v)=\delta(u,v)\Fv(u,v)$ and $S=N^{-1}\sum_{i<j}\Fv_{ij}\xi_{ij}$. In both cases $A=\Ghat+\kappa I_d$ with $\kappa=\nu/(m\tau^2)$ or $\nu/(N\tau^2)$, and the Gram event is $\mathcal{E}=\{\Ghat\succeq G/2\}$, on which $A\succeq G/2$, so $G\preceq2A$ and $A^{-1}\preceq2G^{-1}$.

\begin{lemma}[Projection]
\label{lem:proj}
For every graphon $\wstar$, the projection $\wproj=(\bproj)^{\top}\Fv$ with $\bproj=G^{-1}\int\Fv\wstar$ is the unique minimizer of $\norm{w-\wstar}$ over $w\in\calH$; it satisfies $\int(\wstar-\wproj)\Fv=0$ and $\norm{\wstar-\wproj}\le\min_j\norm{\wstar-w_j}$.
\end{lemma}

\begin{proof}
For $\beta\in\R^d$, $\norm{\beta^{\top}\Fv-\wstar}^2=\beta^{\top}G\beta-2\beta^{\top}\int\Fv\wstar+\norm{\wstar}^2$ is a strictly convex quadratic because $G\succ0$, with unique minimizer $\bproj$; the first-order condition $G\bproj=\int\Fv\wstar$ is $\int(\wstar-\wproj)\Fv=0$. Each $w_j\in\calH$, so the minimum is at most $\norm{\wstar-w_j}$.
\end{proof}

\begin{lemma}[Matrix Chernoff for averages]
\label{lem:chernoff}
Let $Z_1,\dots,Z_k$ be independent symmetric positive semidefinite $d\times d$ random matrices with $\E Z_t=I_d$ and $\norm{Z_t}_{\mathrm{op}}\le L$ almost surely, $L\ge1$. Then
\[
\Prob\Big(\lambda_{\min}\Big(\frac1k\sum_{t=1}^{k}Z_t\Big)\le\frac12\Big)\;\le\;d\exp\Big(-\frac{k}{7L}\Big).
\]
\end{lemma}

\begin{proof}
Fix $\theta>0$. Since $\exp\{-\theta\lambda_{\min}(A)\}=\lambda_{\max}(e^{-\theta A})\le\tr e^{-\theta A}$ for symmetric $A$, Markov's inequality gives
\[
\Prob\Big(\lambda_{\min}\Big(\sum_tZ_t\Big)\le\frac k2\Big)\le e^{\theta k/2}\,\E\,\tr\exp\Big(-\theta\sum_tZ_t\Big)\le e^{\theta k/2}\,\tr\exp\Big(\sum_t\log\E e^{-\theta Z_t}\Big),
\]
the last step by the master bound for sums of independent random matrices, which rests on Lieb's concavity theorem \citep{Tropp2012}. On $[0,L]$ convexity gives $e^{-\theta x}\le1+(e^{-\theta L}-1)x/L$, and applying this spectrally to $0\preceq Z_t\preceq LI_d$ and taking expectations with $\E Z_t=I_d$ yields $\E e^{-\theta Z_t}\preceq\{1+(e^{-\theta L}-1)/L\}I_d$, hence $\log\E e^{-\theta Z_t}\preceq\{(e^{-\theta L}-1)/L\}I_d$ by $\log(1+x)\le x$. Therefore the bound is $d\exp\{\theta k/2+k(e^{-\theta L}-1)/L\}$. The choice $\theta=(\log2)/L$ gives exponent $-k(1-\log2)/(2L)\le-k/(7L)$ because $(1-\log2)/2=0.1534>1/7$.
\end{proof}

\begin{lemma}[Gram concentration]
\label{lem:gram}
Let $L=\supn{\Fv^{\top}G^{-1}\Fv}$. Then
\[
\Prob(\mathcal{E}^c)\le d\exp\{-m/(7L)\}\ \text{under (2)},\qquad
\Prob(\mathcal{E}^c)\le d\exp\{-\lfloor n/2\rfloor/(7L)\}\ \text{under (3)}.
\]
\end{lemma}

\begin{proof}
Let $Z_s=G^{-1/2}\Fv_s\Fv_s^{\top}G^{-1/2}$, so $\E Z_s=I_d$, $\norm{Z_s}_{\mathrm{op}}=\Fv_s^{\top}G^{-1}\Fv_s\le L$, and $\Ghat\succeq G/2$ if and only if $\lambda_{\min}(m^{-1}\sum_sZ_s)\ge1/2$. The independent-dyad case is Lemma~\ref{lem:chernoff} with $k=m$.

For one graph, let $k=\lfloor n/2\rfloor$ and, for a permutation $\pi$ of $\{1,\dots,n\}$, let
\[
W_\pi=k^{-1}\sum_{t=1}^{k}Z_{\pi(2t-1)\pi(2t)},
\]
an average over a perfect matching of $2k$ nodes, hence over $k$ dyads with pairwise disjoint node pairs, so the summands of $W_\pi$ are independent with mean $I_d$ and norm at most $L$. By symmetry every unordered dyad $\{i,j\}$ appears in the same number, $n!\,k/N$, of the pairs $(\pi,t)$, so
\[
\frac1N\sum_{i<j}Z_{ij}=\frac{1}{n!}\sum_{\pi}W_\pi
\]
\citep{Hoeffding1963}. Concavity of $\lambda_{\min}$ gives $\lambda_{\min}(N^{-1}\sum Z_{ij})\ge(n!)^{-1}\sum_\pi\lambda_{\min}(W_\pi)$, and convexity of the exponential gives, for $\theta>0$,
\[
\begin{aligned}
\Prob(\mathcal{E}^c)&\le e^{\theta/2}\,\E\exp\Big(-\frac{\theta}{n!}\sum_\pi\lambda_{\min}(W_\pi)\Big)\\
&\le e^{\theta/2}\,\frac{1}{n!}\sum_\pi\E\exp\{-\theta\lambda_{\min}(W_\pi)\}=e^{\theta/2}\,\E\exp\{-\theta\lambda_{\min}(W)\},
\end{aligned}
\]
where $W$ is distributed as any $W_\pi$. Writing $\theta=k\theta'$ and bounding $\exp\{-\theta\lambda_{\min}(W)\}\le\tr\exp\{-\theta'\sum_{t\le k}Z_t\}$ as in Lemma~\ref{lem:chernoff} yields $\Prob(\mathcal{E}^c)\le d\exp\{k\theta'/2+k(e^{-\theta'L}-1)/L\}$, and the same choice of $\theta'$ completes the proof.
\end{proof}

\begin{lemma}[Error decomposition]
\label{lem:decomp}
In either design, $\bhat-\bproj=A^{-1}(S-\kappa\bproj)$ and, for the unclipped $\what=\bhat^{\top}\Fv$,
\[
\norm{\what-\wstar}^2=\norm{\wstar-\wproj}^2+(\bhat-\bproj)^{\top}G(\bhat-\bproj),
\]
and on $\mathcal{E}$, $(\bhat-\bproj)^{\top}G(\bhat-\bproj)\le8S^{\top}G^{-1}S+8\kappa^2\norm{\bproj}^2/\lambda$. Under (2), $\E[S^{\top}G^{-1}S]=c(\wstar)/m$.
\end{lemma}

\begin{proof}
The outcome decomposition $Y_s=\Fv_s^{\top}\bproj+\xi_s$ gives $m^{-1}\sum\Fv_sY_s=\Ghat\bproj+S$, hence $\bhat-\bproj=A^{-1}(S-\kappa\bproj)$ since $A^{-1}\Ghat-I_d=-\kappa A^{-1}$; the same algebra holds with $N$ under (3). The Pythagorean identity follows from $\what-\wproj\in\calH$ and Lemma~\ref{lem:proj}. On $\mathcal{E}$, with $\zeta=S-\kappa\bproj$,
\[
\zeta^{\top}A^{-1}GA^{-1}\zeta\le2\zeta^{\top}A^{-1}\zeta\le4S^{\top}A^{-1}S+4\kappa^2(\bproj)^{\top}A^{-1}\bproj\le8S^{\top}G^{-1}S+8\kappa^2\norm{\bproj}^2/\lambda,
\]
using $G\preceq2A$, $A^{-1}\preceq2G^{-1}$ and $G^{-1}\preceq\lambda^{-1}I_d$. Under (2),
\[
\E[S^{\top}G^{-1}S]=m^{-2}\textstyle\sum_{s,t}\E[\xi_s\xi_t\Fv_s^{\top}G^{-1}\Fv_t];
\]
the off-diagonal terms vanish by independence and $\E[\xi\Fv]=\int\delta\Fv=0$, and the diagonal gives $m^{-1}\E[\sigma^2\Fv^{\top}G^{-1}\Fv]=c(\wstar)/m$.
\end{proof}

\begin{lemma}[One-graph score moments]
\label{lem:onegraphscore}
Under the one-graph design (3), with $g$, $\kappa_1$ and $\Sigma_g$ as in (1) of the main text:
\begin{enumerate}
\item[(a)] $\displaystyle\E[S^{\top}G^{-1}S]=\frac{c(\wstar)}{N}+\frac{4(n-2)}{n(n-1)}\,\kappa_1(\wstar)$, exactly, for every $n\ge2$.
\item[(b)] Write $S=S_1+S_2$ with $S_1=N^{-1}\sum_{i<j}h(U_i,U_j)$ and $S_2=N^{-1}\sum_{i<j}\Fv_{ij}\{A_{ij}-\wstar(U_i,U_j)\}$, and let $\hat S_1=(2/n)\sum_{i=1}^{n}g(U_i)$. Then $\E\norm{S_2}^2\le d/(4N)$ and $\E\norm{S_1-\hat S_1}^2\le2\,\E\norm{h(U,V)}^2/\{n(n-1)\}$.
\end{enumerate}
\end{lemma}

\begin{proof}
(a) Expand $\E[S^{\top}G^{-1}S]=N^{-2}\sum\E[\xi_{ij}\xi_{kl}\Fv_{ij}^{\top}G^{-1}\Fv_{kl}]$ over ordered pairs of dyads. Given $U$, the $\xi_{ij}$ are independent with $\E[\xi_{ij}\mid U]=\delta(U_i,U_j)$ and $\E[\xi_{ij}^2\mid U]=\sigma^2(U_i,U_j)$. The $N$ equal-dyad terms each contribute $\E[\sigma^2\Fv^{\top}G^{-1}\Fv]=c(\wstar)$. For dyads with disjoint node sets, conditioning on $U$ and using independence of the two position pairs, the term is $\E[\delta_{ij}\Fv_{ij}]^{\top}G^{-1}\E[\delta_{kl}\Fv_{kl}]=0$ by Lemma~\ref{lem:proj}. For the $2N(n-2)$ ordered pairs sharing exactly one node, say $(i,j)$ and $(i,l)$ with $j\neq l$, conditioning first on $U$ and then on $U_i$, and using that $(U_j,A_{ij})$ and $(U_l,A_{il})$ are conditionally independent given $U_i$,
\[
\E[\xi_{ij}\xi_{il}\Fv_{ij}^{\top}G^{-1}\Fv_{il}]=\E\big[\E[\delta_{ij}\Fv_{ij}\mid U_i]^{\top}G^{-1}\E[\delta_{il}\Fv_{il}\mid U_i]\big]=\E[g(U_i)^{\top}G^{-1}g(U_i)]=\kappa_1(\wstar).
\]
Summing, $\E[S^{\top}G^{-1}S]=N^{-2}\{Nc(\wstar)+2N(n-2)\kappa_1(\wstar)\}$, which is (a) since $2(n-2)/N=4(n-2)/\{n(n-1)\}$.

(b) Given $U$, the summands of $S_2$ are independent and mean zero, so $\E[\norm{S_2}^2\mid U]=N^{-2}\sum_{i<j}\wstar(1-\wstar)(U_i,U_j)\norm{\Fv_{ij}}^2\le d/(4N)$ using $\wstar(1-\wstar)\le1/4$ and $\norm{\Fv}^2\le d$. For $S_1$: the kernel $h$ is symmetric with $\int\!\!\int h=0$ and $\E[h(u,V)]=g(u)$, so $\E[S_1\mid U_i]=(2/n)g(U_i)$ and $\hat S_1$ is the H\'ajek projection. The variance count of (a), applied to the kernel $h$, gives $\E\norm{S_1}^2=\E\norm{h}^2/N+\{4(n-2)/(n(n-1))\}\tr\Sigma_g$, while $\E\norm{\hat S_1}^2=(4/n)\tr\Sigma_g$ and $\E[S_1^{\top}\hat S_1]=(2/n)\sum_i\E[\E[S_1\mid U_i]^{\top}g(U_i)]=(4/n)\tr\Sigma_g$. Hence
\[
\E\norm{S_1-\hat S_1}^2=\E\norm{S_1}^2-\E\norm{\hat S_1}^2=\frac{\E\norm{h}^2}{N}-\frac{4\tr\Sigma_g}{n(n-1)}\le\frac{2\,\E\norm{h}^2}{n(n-1)}. \qedhere
\]
\end{proof}

\subsection*{Proof of Theorem 1}

(i) $\what_{\mathrm{sel}}$ is one of the $w_j$, so $\norm{\what_{\mathrm{sel}}-\wstar}^2\ge\min_j\norm{w_j-\wstar}^2$; each $w_j\in\calC$, so the minimum is at least $\dist^2(\wstar,\calC)$; and $\calC\subset\calH$ gives $\dist(\wstar,\calC)\ge\dist(\wstar,\calH)=\norm{\wstar-\wproj}$, likewise $\norm{\what_{\calC}-\wstar}\ge\dist(\wstar,\calC)$. If $\wproj\in\calC$, equality holds. Conversely, $\calC$ is the continuous image of the compact simplex, so the distance to $\calC$ is attained at some $h\in\calC\subset\calH$; if $\norm{\wstar-h}=\norm{\wstar-\wproj}$ then $h$ minimizes over $\calH$, and uniqueness in Lemma~\ref{lem:proj} gives $h=\wproj\in\calC$.

(ii) Write $e_l$ for the elementary symmetric function of degree $l$ in $(\tilde w_1,\dots,\tilde w_J)$, so that pointwise
\[
\wstar=1-\prod_{j=1}^{J}(1-\rho\tilde w_j)=\sum_{l=1}^{J}(-1)^{l+1}\rho^{\,l}e_l=\rho\,e_1-\rho^2e_2+R_\rho,\qquad R_\rho=\sum_{l=3}^{J}(-1)^{l+1}\rho^{\,l}e_l .
\]
Each $\tilde w_j=w_j/\rho\in\calH$, so $\rho e_1\in\calH$ and $\dist(\wstar,\calH)=\dist(-\rho^2e_2+R_\rho,\calH)$, because the distance to a linear subspace is invariant under translation by its elements. Since $0\le e_l\le\binom Jl$ pointwise and $\rho\le1$, $\norm{R_\rho}\le\supn{R_\rho}\le\sum_{l\ge3}\rho^{\,l}\binom Jl\le c_J\rho^3$ with $c_J=\sum_{l=3}^{J}\binom Jl=2^J-1-J-J(J-1)/2$. The distance to a subspace satisfies $|\dist(x+y,\calH)-\dist(x,\calH)|\le\norm y$ and $\dist(-\rho^2e_2,\calH)=\rho^2\dist(e_2,\calH)$, giving the two-sided bound of (ii). For the crude upper bound, the first two Bonferroni inequalities for the union of independent events give $0\le\sum_jw_j-\wstar\le\sum_{j<k}w_jw_k$ pointwise, so with $\sum_jw_j\in\calH$,
\[
\norm{\wstar-\wproj}\le\Big\|\sum_jw_j-\wstar\Big\|\le\sum_{j<k}\norm{w_jw_k}\le\rho^2\sum_{j<k}\min(\norm{\tilde w_j},\norm{\tilde w_k})=\rho^2a,
\]
using $\norm{w_jw_k}\le\min(\supn{w_j}\norm{w_k},\supn{w_k}\norm{w_j})$ and $\supn{w_j}\le\rho$.

For the hull, take $\pi\in\Delta_J$ and write $\wstar-\sum_j\pi_jw_j=\sum_j(1-\pi_j)w_j-\Gamma$ with $\Gamma=\sum_jw_j-\wstar$, which satisfies $0\le\Gamma$ and $\norm\Gamma\le\rho^2a$ as above. The agents are nonnegative, so the cross terms in $\norm{\sum_j(1-\pi_j)w_j}^2$ are nonnegative and $\norm{\sum_j(1-\pi_j)w_j}^2\ge\sum_j(1-\pi_j)^2\norm{w_j}^2\ge\rho^2b^2\sum_j(1-\pi_j)^2\ge\rho^2b^2(J-1)^2/J$, the last step by Cauchy--Schwarz from $\sum_j(1-\pi_j)=J-1$. The triangle inequality then gives $\norm{\wstar-\sum_j\pi_jw_j}\ge\rho\{(J-1)J^{-1/2}b-\rho a\}$, and taking the infimum over $\pi$ and the positive part proves the lower bound; the upper bound is $\dist(\wstar,\calC)\le\norm{\wstar-w_1}=\norm{\sum_{j\ge2}w_j-\Gamma}\le\rho\sum_{j\ge2}\norm{\tilde w_j}+\rho^2a$.

(iii) If $e_2\notin\calH$ then $\dist(e_2,\calH)>0$ and (ii) gives $\norm{\wstar-\wproj}=\rho^2\dist(e_2,\calH)\{1+O(\rho)\}$ with the $O(\rho)$ term bounded in absolute value by $c_J\rho/\dist(e_2,\calH)$; squaring gives the claim. For $J=2$, $c_2=0$, the expansion terminates at $l=2$, and $\dist(\wstar,\calH)=\rho^2\dist(\tilde w_1\tilde w_2,\calH)$ exactly for every $\rho$. For the hull, if $\rho\le(J-1)b/(2\sqrt Ja)$ the bracket in the lower bound is at least $(J-1)b/(2\sqrt J)$, giving $\dist^2(\wstar,\calC)\ge\rho^2(J-1)^2b^2/(4J)$, and the upper bound of (ii) is $O(\rho^2)$. \qed

\subsection*{Proof of Theorem $1'$}
Write $\wstar=\rho e_1-\rho^2e_2+R_\rho$ with $\norm{R_\rho}\le c_J\rho^3$ as in the proof of Theorem~1, and $\tilde F=(1,\tilde w_1,\dots,\tilde w_J)$ with Gram matrix $\tilde G$, invertible by Assumption~1.

\emph{(i) Sum-to-one with a free intercept.} For $c\in\R$ and $\pi\in\Delta_J$, $\wstar-c-\sum_j\pi_jw_j=\rho\{\sum_j(1-\pi_j)\tilde w_j-c/\rho\}-\rho^2e_2+R_\rho$, so $\dist(\wstar,\calC_c)\ge\rho\,b_c-\rho^2\norm{e_2}-c_J\rho^3$ with $b_c=\min_{\pi\in\Delta_J}\dist\big(\sum_j(1-\pi_j)\tilde w_j,\ \mathrm{span}(1)\big)$. The map $\pi\mapsto\dist(\sum_j(1-\pi_j)\tilde w_j,\mathrm{span}(1))$ is continuous on the compact simplex and never zero: $\sum_j(1-\pi_j)=J-1\ge1$, so $\sum_j(1-\pi_j)\tilde w_j$ is a nonzero element of $\mathrm{span}(\tilde w_1,\dots,\tilde w_J)$, which meets $\mathrm{span}(1)$ only at zero by Assumption~1. Hence $b_c>0$ and $\dist^2(\wstar,\calC_c)\ge\rho^2b_c^2/4$ for small $\rho$. The upper bound is attained by $c=0$, $\pi=e_1$: $\norm{\wstar-w_1}\le\rho\norm{e_1-\tilde w_1}+\rho^2\norm{e_2}+c_J\rho^3=O(\rho)$. Since $\calC\subset\calC_c$, the same bounds hold for $\calC$, and $\dist(\wstar,\calC)\le\norm{\wstar-w_1}$ gives the matching upper bound.

\emph{(ii) Cones and spans.} With intercept: the projection of $\wstar$ onto $\calH=\mathrm{span}(\tilde F)$ has coefficients $\tilde\beta(\rho)=\tilde G^{-1}\int\tilde F\wstar=\rho\tilde G^{-1}\int\tilde Fe_1+O(\rho^2)$, and $\int\tilde F e_1=\int\tilde F\sum_j\tilde w_j=\tilde G(0,1,\dots,1)^{\top}$, hence $\tilde\beta_j(\rho)=\rho\{1+O(\rho)\}$ for every template coordinate $j\ge1$, with the $O(\rho)$ constant depending only on $\tilde G^{-1}$, $\norm{e_2}$ and $c_J$; the coefficient on the agent $w_j=\rho\tilde w_j$ itself is $\beta^{\circ}_j=\tilde\beta_j(\rho)/\rho=1+O(\rho)$. For $\rho\le\rho_0$ these are positive, so the projection $\wproj$ lies in the relative interior of $\calK$ within $\calH$; it minimizes $\norm{w-\wstar}$ over $\calH\supset\calK$ and belongs to $\calK$, hence it is the minimizer over $\calK$ and $\dist(\wstar,\calK)=\dist(\wstar,\calH)$. Each $\tilde w_j=w_j/\rho\in\calH$, so $\rho e_1\in\calH$ and $\dist(\wstar,\calH)=\dist(-\rho^2e_2+R_\rho,\calH)$, because the distance to a linear subspace is invariant to translation along it; then $\rho^2\dist(e_2,\calH)-c_J\rho^3\le\dist(\wstar,\calH)\le\rho^2\dist(e_2,\calH)+c_J\rho^3$, which is the fourth-order statement when $e_2\notin\calH$. Without intercept the same computation applies to $\tilde F_0=(\tilde w_1,\dots,\tilde w_J)$ with Gram matrix $\tilde G_0$, since $\int\tilde F_0e_1=\tilde G_0(1,\dots,1)^{\top}$: the coefficients are $1+O(\rho)$, the projection onto $\calH_0$ lies in $\calK_0$ for small $\rho$, $\rho e_1\in\calH_0$, and $\dist(\wstar,\calH_0)=\rho^2\dist(e_2,\tilde{\calH}_0)+O(\rho^3)$, with $\tilde{\calH}_0=\calH_0$ as sets. \qed

\subsection*{Proof of Theorem $1''$}
Since $0\le w_j\le\rho<1$, each $x_j=-\log(1-w_j)$ is finite and nonnegative, and $1-\prod_{j}(1-w_j)=1-\exp\{\sum_j\log(1-w_j)\}=1-\exp(-\sum_jx_j)=w_{(0,1,\dots,1)}$, which lies in $\calN$; hence $\dist(\wstar,\calN)=0$. For the expansion, $x_j=w_j+r_j$ with $0\le r_j\le w_j^2/(1-w_j)\le\rho^2/(1-\rho)$, so for $\gamma$ in a bounded set $t:=\sum_j\gamma_jx_j=\sum_j\gamma_jw_j+O(\rho^2)$ with $t=O(\rho)$, and $\exp(-\gamma_0-t)=e^{-\gamma_0}\{1-t+O(t^2)\}$ gives $w_\gamma=(1-e^{-\gamma_0})+e^{-\gamma_0}\sum_j\gamma_jw_j+O(\rho^2)$ uniformly. The leading term has intercept $1-e^{-\gamma_0}$ and nonnegative slopes $e^{-\gamma_0}\gamma_j$, so it belongs to $\calK$. \qed

\subsection*{Proof of Theorem 2}

(i) is Lemma~\ref{lem:onegraphscore}(a).

(ii) On $\mathcal{E}$, Lemma~\ref{lem:decomp} and the contraction property of clipping give $\norm{\whatp-\wstar}^2\le\norm{\wstar-\wproj}^2+8S^{\top}G^{-1}S+8\kappa^2\norm{\bproj}^2/\lambda$; on $\mathcal{E}^c$, $\norm{\whatp-\wstar}^2\le1$ because both functions take values in $[0,1]$. Taking expectations, bounding $\E[S^{\top}G^{-1}S\ind_{\mathcal{E}}]$ by the identity of (i), and applying the one-graph part of Lemma~\ref{lem:gram} completes the proof.

(iii) By Lemma~\ref{lem:decomp}, $\bhat-\bproj=A^{-1}(S-\kappa\bproj)$. Each entry of $\Ghat$ is a $U$-statistic of the $U_i$ with a kernel bounded by one, so by the variance count in the proof of Lemma~\ref{lem:onegraphscore}, $\operatorname{Var}(\Ghat_{ab})\le4/n+2/\{n(n-1)\}\to0$ and $\Ghat\to G$ in probability; with $\kappa\to0$, $A^{-1}\to G^{-1}$ in probability.

Suppose $\kappa_1(\wstar)>0$, equivalently $\Sigma_g\neq0$. Lemma~\ref{lem:onegraphscore}(b) gives $n\,\E\norm{S_2}^2\le nd/(4N)\to0$ and $n\,\E\norm{S_1-\hat S_1}^2\le2n\,\E\norm h^2/\{n(n-1)\}\to0$, so $\sqrt n\,S=\sqrt n\,\hat S_1+o_{L_2}(1)$. The vectors $g(U_i)$ are independent, bounded and mean zero (Lemma~\ref{lem:proj} gives $\int g=0$) with covariance $\Sigma_g$, so the multivariate central limit theorem gives $\sqrt n\,\hat S_1=(2/\sqrt n)\sum_ig(U_i)\Rightarrow N(0,4\Sigma_g)$. Since $\sqrt n\,\kappa\norm{\bproj}\to0$, Slutsky's lemma yields $\sqrt n(\bhat-\bproj)\Rightarrow N(0,4G^{-1}\Sigma_gG^{-1})$.

If $\wstar\in\calH$ then $\delta\equiv0$, $S=S_2$, and given $U$ the summands $\Fv_{ij}\{A_{ij}-\wstar(U_i,U_j)\}$ are independent, mean zero and bounded by $\sqrt d$, with conditional covariance $N^{-1}\sum_{i<j}\wstar(1-\wstar)(U_i,U_j)\Fv_{ij}\Fv_{ij}^{\top}\to M$ in probability, entrywise, by the same $U$-statistic variance bound. The Lindeberg condition holds conditionally because the normalized summands are bounded by $\sqrt{d/N}$, so the conditional characteristic function of $\sqrt N\,S_2$ converges in probability to $\exp(-t^{\top}Mt/2)$ for every $t$, and dominated convergence gives $\sqrt N\,S_2\Rightarrow N(0,M)$; Slutsky's lemma as before gives $\sqrt N(\bhat-\bproj)\Rightarrow N(0,G^{-1}MG^{-1})$.

In both cases $N\Sighat=\nu A^{-1}\to\nu G^{-1}$ in probability. For the coverage claim, let $s^2=4\ell^{\top}G^{-1}\Sigma_gG^{-1}\ell>0$: the interval is $\ell^{\top}\bhat\pm z\{\ell^{\top}\Sighat\ell\}^{1/2}$ with $z=z_{1-\alpha/2}$, and $n\,\ell^{\top}\Sighat\ell=(n/N)\cdot N\ell^{\top}\Sighat\ell\to0$ in probability. For any $\epsilon>0$,
\[
\Prob\big(|\ell^{\top}(\bhat-\bproj)|\le z\{\ell^{\top}\Sighat\ell\}^{1/2}\big)\le\Prob\big(|\sqrt n\,\ell^{\top}(\bhat-\bproj)|\le\epsilon\big)+\Prob\big(z\{n\,\ell^{\top}\Sighat\ell\}^{1/2}>\epsilon\big)\to2\Phi(\epsilon/s)-1,
\]
and letting $\epsilon\downarrow0$ shows the coverage tends to zero. If $\wstar\in\calH$, the coverage calculation in the proof of Proposition~S2, with $\sqrt N$ in place of $\sqrt m$ and the limits just established, gives the stated limit. \qed

\subsection*{Proposition $2''$, parts (iv) and (v)}
The statement continues from the main text.
\begin{enumerate}
\item[(iv)] With the plug-in working variance $\hat\nu\asymp\rho_n$ of Section~5 of the main text, the ridge term of Theorem~2 of the main text(ii) is $O(\norm{\bar\beta^{\circ}}^2/(N^2\tau^4))$, negligible against the dyad term whenever $N\tau^4\rho_n\to\infty$ (a sufficient condition; it is also necessary when the projected slopes $\bar\beta^{\circ}_j$, $j\ge1$, are not all zero).
\item[(v)] Let $D_n=\mathrm{diag}(1,\rho_n,\dots,\rho_n)$ and $\kappa=\hat\nu/(N\tau^2)$ with $\hat\nu\asymp\rho_n$. If $\bar\kappa_1>0$ and $n\rho_n\to\infty$, then $(\sqrt n/\rho_n)\,D_n(\bhat-\beta^{\circ}_n)\Rightarrow N(0,4\tilde G^{-1}\bar\Sigma_g\tilde G^{-1})$ with $\bar\Sigma_g=\int\bar g\bar g^{\top}$; if $\bar w\in\tilde{\calH}$ and $N\rho_n\to\infty$, then $\sqrt{N/\rho_n}\,D_n(\bhat-\beta^{\circ}_n)\Rightarrow N(0,\tilde G^{-1}\bar M\tilde G^{-1})$ with $\bar M=\int\bar w\,\tilde F\tilde F^{\top}$. In the original coordinates the slopes converge at rate $\sqrt n$ and $\sqrt{N\rho_n}$ and the intercept at $\sqrt n/\rho_n$ and $\sqrt{N/\rho_n}$.
\end{enumerate}

\subsection*{Proof of Proposition $2''$}
Write $D_n=\mathrm{diag}(1,\rho_n,\dots,\rho_n)$, so that $F_n=D_n\tilde F$, $G_n=\int F_nF_n^{\top}=D_n\tilde GD_n$ and $G_n^{-1}=D_n^{-1}\tilde G^{-1}D_n^{-1}$.

(i) The span of $(1,\rho_n\tilde w_1,\dots,\rho_n\tilde w_J)$ is that of $(1,\tilde w_1,\dots,\tilde w_J)$ because $\rho_n>0$. The leverage is invariant, $F_n^{\top}G_n^{-1}F_n=\tilde F^{\top}D_nD_n^{-1}\tilde G^{-1}D_n^{-1}D_n\tilde F=\tilde\ell$. The empirical Gram matrix is $\hat G_n=N^{-1}\sum_{i<j}F_{n,ij}F_{n,ij}^{\top}=D_n\hat{\tilde G}D_n$, and for a positive definite congruence $\hat G_n-G_n/2=D_n(\hat{\tilde G}-\tilde G/2)D_n$ is positive semidefinite if and only if $\hat{\tilde G}-\tilde G/2$ is. For the eigenvalue bounds, $e_j^{\top}G_ne_j=\rho_n^2\norm{\tilde w_j}^2$ for $j\ge1$ gives the upper bound, and for a unit vector $v$, $v^{\top}G_nv=(D_nv)^{\top}\tilde G(D_nv)\ge\lambda_{\min}(\tilde G)\norm{D_nv}^2\ge\lambda_{\min}(\tilde G)\rho_n^2$ because every diagonal entry of $D_n$ is at least $\rho_n$ when $\rho_n\le1$.

(ii) Projection onto $\tilde{\calH}$ is linear, so $w^{\circ}_n=\Pi(\rho_n\bar w)=\rho_n\bar w^{\circ}$ and $\delta_n=\rho_n\bar\delta$. Writing $\bar w^{\circ}=\bar\beta^{\circ\top}\tilde F$ gives $w^{\circ}_n=\rho_n\bar\beta^{\circ}_0+\sum_{j\ge1}\bar\beta^{\circ}_j\,(\rho_n\tilde w_j)$, which identifies the coefficients on $F_n$.

(iii) The residual second moment is $\sigma_n^2=\wstar_n(1-\wstar_n)+\delta_n^2=\rho_n\bar w-\rho_n^2\bar w^2+\rho_n^2\bar\delta^2$, and $c(\wstar_n)=\int\sigma_n^2F_n^{\top}G_n^{-1}F_n=\int\sigma_n^2\tilde\ell$ by (i), which is the stated expression. Next $g_n(u)=\int\delta_n(u,v)F_n(u,v)\,dv=\rho_nD_n\bar g(u)$, so $\kappa_1(\wstar_n)=\int g_n^{\top}G_n^{-1}g_n=\rho_n^2\int\bar g^{\top}D_nD_n^{-1}\tilde G^{-1}D_n^{-1}D_n\bar g=\rho_n^2\bar\kappa_1$. Substituting into Theorem~2(i) gives the display, and the ratio of the second term to the first is
\[
\frac{4(n-2)\rho_n^2\bar\kappa_1/\{n(n-1)\}}{\{\rho_n\int\bar w\tilde\ell+\rho_n^2\int(\bar\delta^2-\bar w^2)\tilde\ell\}/N}=\frac{2(n-2)\rho_n\bar\kappa_1}{\int\bar w\tilde\ell+\rho_n\int(\bar\delta^2-\bar w^2)\tilde\ell},
\]
which is of exact order $n\rho_n$ since the denominator tends to $\int\bar w\tilde\ell>0$.

(iv) With $\kappa=\hat\nu/(N\tau^2)$ and $\hat\nu\asymp\rho_n$, the ridge term $8\kappa^2\norm{\beta^{\circ}_n}^2/\lambda_{\min}(G_n)$ is at most of order $\rho_n^2\norm{\bar\beta^{\circ}}^2/(N^2\tau^4\rho_n^2)=\norm{\bar\beta^{\circ}}^2/(N^2\tau^4)$ by (i) and (ii), and of exactly that order when some slope $\bar\beta^{\circ}_j$, $j\ge1$, is nonzero (if all slopes vanish, $\norm{\beta^{\circ}_n}=\rho_n|\bar\beta^{\circ}_0|$ and the term is smaller); its ratio to the dyad term $\rho_n\int\bar w\tilde\ell/N$ is of order $1/(N\tau^4\rho_n)$.

(v) Write $\xi_{ij}=A_{ij}-w^{\circ}_n(U_i,U_j)=\{A_{ij}-\wstar_n(U_i,U_j)\}+\rho_n\bar\delta(U_i,U_j)$ and $S_n=N^{-1}\sum_{i<j}F_{n,ij}\xi_{ij}$. Since $D_n^{-1}F_{n,ij}=\tilde F_{ij}$,
\[
D_n^{-1}S_n=E_n+\rho_n\bar T,\qquad E_n=\frac1N\sum_{i<j}\tilde F_{ij}\{A_{ij}-\wstar_n(U_i,U_j)\},\qquad \bar T=\frac1N\sum_{i<j}\tilde F_{ij}\bar\delta(U_i,U_j),
\]
where $\bar T$ is a $U$-statistic whose kernel does not depend on $n$, with $\E\bar T=0$ by the normal equations of $\bar w^{\circ}$ and Hájek projection $(2/n)\sum_i\bar g(U_i)$. The ridge fit satisfies $\bhat-\beta^{\circ}_n=(\hat G_n+\kappa I)^{-1}S_n-\kappa(\hat G_n+\kappa I)^{-1}\beta^{\circ}_n$, and multiplying by $D_n$, with $\hat G_n=D_n\hat{\tilde G}D_n$,
\[
D_n(\bhat-\beta^{\circ}_n)=(\hat{\tilde G}+\kappa D_n^{-2})^{-1}D_n^{-1}S_n-\kappa(\hat{\tilde G}+\kappa D_n^{-2})^{-1}D_n^{-1}\beta^{\circ}_n .
\]
Here $\hat{\tilde G}\to\tilde G$ in probability (a $U$-statistic law of large numbers), $\kappa\norm{D_n^{-2}}=\kappa\rho_n^{-2}\asymp1/(N\rho_n)\to0$, and $\norm{D_n^{-1}\beta^{\circ}_n}\asymp\rho_n^{-1}$, so the ridge term is $O_p(\kappa/\rho_n)=O_p(1/N)$. Conditionally on $U$ the summands of $E_n$ are independent, centered and bounded by $\norm{\tilde F}$, with $\E(\norm{E_n}^2\mid U)\le C\rho_n/N$.

\emph{Loaded misspecification.} If $\bar\kappa_1>0$, the classical $U$-statistic central limit theorem gives $\sqrt n\,\bar T\Rightarrow N(0,4\bar\Sigma_g)$, while $\E\norm{(\sqrt n/\rho_n)E_n}^2\le Cn/(N\rho_n)\asymp1/(n\rho_n)\to0$ and $(\sqrt n/\rho_n)\cdot O_p(1/N)\to0$ when $n\rho_n\to\infty$. Hence $(\sqrt n/\rho_n)D_n(\bhat-\beta^{\circ}_n)=\tilde G^{-1}\sqrt n\,\bar T+o_p(1)\Rightarrow N(0,4\tilde G^{-1}\bar\Sigma_g\tilde G^{-1})$ by Slutsky.

\emph{Correct specification.} If $\bar w\in\tilde{\calH}$ then $\bar\delta=0$ and $D_n^{-1}S_n=E_n$. Conditionally on $U$, $\sqrt{N/\rho_n}\,E_n$ is a sum of $N$ independent centered vectors bounded by $\norm{\tilde F}/\sqrt{N\rho_n}$, with conditional covariance $N^{-1}\sum_{i<j}\tilde F_{ij}\tilde F_{ij}^{\top}\bar w(U_i,U_j)\{1-\rho_n\bar w(U_i,U_j)\}\to\bar M$ in probability; the Lindeberg condition holds because $\norm{\tilde F}^2/(N\rho_n)\to0$ when $N\rho_n\to\infty$, so the conditional Lindeberg--Feller theorem and the convergence of the conditional covariance to a constant give $\sqrt{N/\rho_n}\,E_n\Rightarrow N(0,\bar M)$ unconditionally. The ridge term contributes $\sqrt{N/\rho_n}\cdot O_p(1/N)=O_p\{(N\rho_n)^{-1/2}\}\to0$, and Slutsky gives $\sqrt{N/\rho_n}\,D_n(\bhat-\beta^{\circ}_n)\Rightarrow N(0,\tilde G^{-1}\bar M\tilde G^{-1})$. Reading off the coordinates of $D_n$ gives the rates stated for slopes and intercept. \qed

\subsection*{Proof of Corollary 1}
Under the superposition model $\wstar=\rho e_1-\rho^2e_2+R_\rho$ with $\rho e_1\in\calH$, so $\delta=\wstar-\wproj=-\rho^2(e_2-\Pi_{\calH}e_2)+(R_\rho-\Pi_{\calH}R_\rho)=-\rho^2r+O(\rho^3)$ in $L^2$ and in $L^\infty$ (the projection onto a finite-dimensional space of bounded functions is bounded in sup norm). With $F=D_\rho\tilde F$, $g(u)=\int\delta(u,v)F(u,v)\,dv=-\rho^2D_\rho\bar g_r(u)+O(\rho^3)$, and $G^{-1}=D_\rho^{-1}\tilde G^{-1}D_\rho^{-1}$ gives $\kappa_1(\wstar)=\int g^{\top}G^{-1}g=\rho^4\int\bar g_r^{\top}\tilde G^{-1}\bar g_r+O(\rho^5)=\rho^4K+O(\rho^5)$. For $c(\wstar)=\int\sigma^2F^{\top}G^{-1}F=\int\sigma^2\tilde\ell$ (leverage invariance, Proposition~$2''$(i)) with $\sigma^2=\wstar(1-\wstar)+\delta^2=\rho e_1+O(\rho^2)$, so $c(\wstar)=\rho\int e_1\tilde\ell+O(\rho^2)$. The ratio of the two terms of Theorem~2(i) is $\{4(n-2)\rho^4K/(n(n-1))\}/\{\rho\int e_1\tilde\ell/N\}\asymp n\rho^3$. \qed

\subsection*{A finite-class inequality for the Brier score}
\begin{lemma}
\label{lem:finite}
Let $(F_s,Y_s)$, $s=1,\dots,m$, be independent with $Y_s\in\{0,1\}$ and $\E(Y_s\mid F_s)=\wstar_s$, let $f_1,\dots,f_K$ be fixed predictors with values in $[0,1]$, $R(f)=m^{-1}\sum_s\E(Y_s-f_s)^2$, $R^\star=R(\wstar)$, $\hat R(f)=m^{-1}\sum_s(Y_s-f_s)^2$, and let $\hat k$ satisfy $\hat R(f_{\hat k})\le\min_k\hat R(f_k)+\varsigma$ for some $\varsigma\ge0$. For $\varepsilon\in(0,1)$, $\delta\in(0,1)$ and $u=\log(2K/\delta)/m$, with probability at least $1-\delta$,
\[
R(f_{\hat k})-R^\star\le\frac{1+\varepsilon}{1-\varepsilon}\min_k\{R(f_k)-R^\star\}+\frac{2c_\varepsilon u+\varsigma}{1-\varepsilon},\qquad c_\varepsilon=\frac2\varepsilon+\frac43 .
\]
\end{lemma}
\begin{proof}
Fix $k$ and set $Z_s=(Y_s-f_{k,s})^2-(Y_s-\wstar_s)^2=(\wstar_s-f_{k,s})(2Y_s-f_{k,s}-\wstar_s)$, so $|Z_s|\le1$, $\E Z_s=\E(f_{k,s}-\wstar_s)^2$ and $\E Z_s^2\le4\E(f_{k,s}-\wstar_s)^2$; with $A_k=R(f_k)-R^\star=m^{-1}\sum_s\E Z_s\ge0$, the average $\bar Z_k=\hat R(f_k)-\hat R(\wstar)$ has mean $A_k$ and Bernstein's inequality gives $P(|\bar Z_k-A_k|\ge t)\le2\exp\{-mt^2/(8A_k+4t/3)\}$. Solving for $t$ at level $\delta/K$ and using $\sqrt{8A_ku}\le\varepsilon A_k+2u/\varepsilon$ yields, simultaneously for all $k$ with probability $1-\delta$, $|\bar Z_k-A_k|\le\varepsilon A_k+c_\varepsilon u$. On this event, with $k^\star$ a population minimizer, $(1-\varepsilon)A_{\hat k}-c_\varepsilon u\le\bar Z_{\hat k}\le\bar Z_{k^\star}+\varsigma\le(1+\varepsilon)A_{k^\star}+c_\varepsilon u+\varsigma$.
\end{proof}

\subsection*{Proof of Theorem 4}
Write $\tilde x=(1,x_1,\dots,x_J)$, $\eta_\gamma=\gamma^{\top}\tilde x$ and $L=L_\epsilon=1+\log(1/\epsilon)$, so that $0\le\tilde x_k\le L$ for every $k$ and $|w_\gamma-w_{\gamma'}|=|e^{-\eta_{\gamma'}}-e^{-\eta_\gamma}|\le|\eta_\gamma-\eta_{\gamma'}|\le L\norm{\gamma-\gamma'}_1$ pointwise.

\emph{(i).} The empirical risk is continuous in $\gamma$ on the compact set $\Gamma_B$, so a minimizer exists. If $w_\gamma=w_{\gamma'}$ a.e. then $\eta_\gamma=\eta_{\gamma'}$ a.e. because $t\mapsto1-e^{-t}$ is injective, hence $(\gamma-\gamma')^{\top}\tilde x=0$ a.e. and $\gamma=\gamma'$ by linear independence. Under the superposition model $\wstar=w_{\gamma^\star}$ (Theorem~$1''$) with $\gamma^\star\in\Gamma_B$ since $\norm{\gamma^\star}_1=J\le B$, and $R(\gamma)-R^\star=\norm{w_\gamma-\wstar}^2$ vanishes only at $\gamma^\star$.

\emph{(ii).} Let $\calG\subset\Gamma_B$ be an $\ell_1$-net of $\Gamma_B$ at scale $\delta_g=1/(mL)$, of cardinality $K\le(3B/\delta_g)^d=(3BLm)^d$. Let $\tilde\gamma\in\calG$ be a nearest point to $\hat\gamma$; since $|(y-p)^2-(y-p')^2|\le2|p-p'|$, $\hat R(\tilde\gamma)\le\hat R(\hat\gamma)+2L\delta_g\le\min_{\calG}\hat R+2/m$, so $\tilde\gamma$ is a selector over $\calG$ with slack $\varsigma=2/m$, and the net contains a point $\gamma_g$ with $\norm{w_{\gamma_g}-\wstar}^2\le(L\delta_g)^2=m^{-2}$. Lemma~\ref{lem:finite} with $\varepsilon=1/2$ gives, with probability $1-\delta$, $\norm{w_{\tilde\gamma}-\wstar}^2\le3m^{-2}+(64/3)u+4/m$ with $u=\{d\log(3BLm)+\log(2/\delta)\}/m$, and $\norm{w_{\hat\gamma}-\wstar}^2\le2\norm{w_{\hat\gamma}-w_{\tilde\gamma}}^2+2\norm{w_{\tilde\gamma}-\wstar}^2\le2m^{-2}+2\norm{w_{\tilde\gamma}-\wstar}^2$. Integrating the tail in $\delta$ (the bound is $a+b\log(2/\delta)$ with $P(\text{excess}>a+bt)\le2e^{-t}$, so the expectation is at most $a+b(1+\log2)$) and using $m\ge3$, $B\ge2$, $L\ge1$, so that $\log(3BLm)\ge2$ and $\log(3BLm)\le\{1+\log(3BL)\}\log m$,
\[
\E\norm{w_{\hat\gamma}-\wstar}^2\le\frac{8}{m^2}+\frac8m+\frac{128}{3}\,\frac{d\log(3BLm)+1.7}{m}\le80\,\{1+\log(3BL)\}\,\frac{d\log m}{m}.
\]

\emph{(iii).} On the net, Hoeffding's inequality and a union bound give $\max_{\calG}|\hat R-R|\le\sqrt{\log(2K/\delta)/(2m)}$ with probability $1-\delta$, and the discretization adds at most $2L\delta_g=2/m$ to each of $\hat R$ and $R$. With $\gamma^{\circ}_B$ a population minimizer, $R(\hat\gamma)\le\hat R(\hat\gamma)+\text{dev}\le\hat R(\gamma^{\circ}_B)+\text{dev}\le R(\gamma^{\circ}_B)+2\,\text{dev}$ with $\text{dev}\le\sqrt{\log(2K/\delta)/(2m)}+2/m$; taking expectations, using $\E\sqrt X\le\sqrt{\E X}$ and $m\ge5$, $R(\hat\gamma)-R^\star\le R(\gamma^{\circ}_B)-R^\star+4\sqrt{1+\log(3BL)}\sqrt{d\log m/m}$, and $R(\gamma^{\circ}_B)-R^\star=\dist^2(\wstar,\{w_\gamma:\gamma\in\Gamma_B\})$. \qed

\subsection*{Proposition S3 and its proof}
\begin{propositionS}[Choosing a rule on held-out dyads]
\label{prop:selectS}
Let $r_1,\dots,r_K$ be predictors with values in $[0,1]$, fitted on data independent of $m_2$ further independent dyads $(F_s,Y_s)$ from the independent-dyad design (2) of the main text, and let $\hat k$ minimize the empirical Brier score of $r_k$ on those dyads. Write $R(r)=\E\{Y-r\}^2$ and $R^\star=\E\{Y-\wstar\}^2$. For every $\varepsilon\in(0,1)$ and $\delta\in(0,1)$, with probability at least $1-\delta$ conditionally on the fitted rules,
\[
R(r_{\hat k})-R^\star\le\frac{1+\varepsilon}{1-\varepsilon}\min_k\{R(r_k)-R^\star\}+\frac{2(2/\varepsilon+4/3)}{1-\varepsilon}\,\frac{\log(2K/\delta)}{m_2}.
\]
\end{propositionS}

Conditionally on the fitted rules, the $m_2$ held-out dyads are independent draws from the independent-dyad design (2) of the main text, so Lemma~\ref{lem:finite} applies with $\varsigma=0$. \qed

\subsection*{Theorem 6, parts (ii) and (iv)}
The statement continues from the main text.
\begin{enumerate}
\item[(ii)] Conditionally on $(U,A_{\calT})$, $\E\norm{\bhat-\beta^{\circ}_{c}}^2_{\hat G}\le2\hat c/m+2\kappa^2\norm{\beta^{\circ}_{c}}^2/\lambda_{\min}(\hat G)$, where $\norm{\beta}^2_{\hat G}=\beta^{\top}\hat G\beta$ is the squared error of the synthesis averaged over the evaluation dyads.
\item[(iv)] When the agents are known, the unconditional identity of Theorem~2 of the main text(i) is the law of total variance for the conditional one: the dyad term $c(\wstar)/N$ is the average of the conditional sampling variance, and the node-level term $4(n-2)\kappa_1(\wstar)/\{n(n-1)\}$ is the H\'ajek part of the variance, under resampling of the nodes, of the conditional projection $\beta^{\circ}_{c}(U)$ around the population projection $\bproj$, the remainder of that variance being the $\int\delta^2\,\Fv^{\top}G^{-1}\Fv/N$ share of $c(\wstar)/N$.
\end{enumerate}

\subsection*{Proof of Theorem 6}
\emph{(i).} Because $\calT$ and $\calE$ are chosen independently of $A$ given $U$ and the $A_{ij}$ are conditionally independent given $U$, the outcomes on $\calE$ are, conditionally on $(U,A_{\calT})$ and the masks, independent Bernoulli$\{\wstar(U_i,U_j)\}$, and $\hat F$ is a fixed function of $(i,j)$. The definition of $\beta^{\circ}_{c}$ is the normal equation $m^{-1}\sum_{\calE}\hat F_{ij}\{\wstar(U_i,U_j)-\hat w^{\circ}_{c}(i,j)\}=0$, so the approximation residual cancels from the score,
\[
S=\frac1m\sum_{\calE}\hat F_{ij}\{A_{ij}-\hat w^{\circ}_{c}\}=\frac1m\sum_{\calE}\hat F_{ij}\{A_{ij}-\wstar(U_i,U_j)\},
\]
which has conditional mean zero. Expanding $S^{\top}\hat G^{-1}S$, the cross terms between distinct dyads vanish by conditional independence and centering, and the diagonal terms give $m^{-2}\sum_{\calE}\wstar(1-\wstar)(U_i,U_j)\hat F_{ij}^{\top}\hat G^{-1}\hat F_{ij}=\hat b/m$. The residual $\wstar-\hat w^{\circ}_{c}$ is a fixed vector orthogonal to the columns of $\hat F$ in the empirical inner product; it biases $\beta^{\circ}_{c}$ as a target for $\wstar$ and contributes nothing to the conditional variance of the score.

\emph{A two-mechanism library with $K=0$, for Corollary~1.} Partition $[0,1]$ into four equal cells carrying independent balanced signs $h,k\in\{-1,1\}$, and set $x=h(u)h(v)$, $y=k(u)k(v)$, $\tilde w_1=(1+x)/2$, $\tilde w_2=(1+y)/2$. Assumption~1 holds, $e_2=\tilde w_1\tilde w_2=(1+x+y+xy)/4$ and $r=e_2-\Pi_{\calH}e_2=xy/4\neq0$, so $e_2\notin\calH$. But $\int x\,y\,dv=0$ and $\int xy\cdot x\,dv=\int xy\cdot y\,dv=0$, so $\bar g_r\equiv0$ and $K=0$: then $\kappa_1(\wstar)=0$ for every $\rho$ and the estimation error stays on the $N^{-1}$ line however large $n\rho^3$ becomes. This is why Corollary~1 requires $K>0$.

\emph{(ii).} Since $m^{-1}\sum_{\calE}\hat F_{ij}A_{ij}=\hat G\beta^{\circ}_{c}+S$, the ridge fit satisfies $\bhat-\beta^{\circ}_{c}=(\hat G+\kappa I)^{-1}S-\kappa(\hat G+\kappa I)^{-1}\beta^{\circ}_{c}$. Using $\norm{a+b}^2_{\hat G}\le2\norm a^2_{\hat G}+2\norm b^2_{\hat G}$, $(\hat G+\kappa I)^{-1}\hat G(\hat G+\kappa I)^{-1}\preceq\hat G^{-1}$ and $\kappa^2(\hat G+\kappa I)^{-1}\hat G(\hat G+\kappa I)^{-1}\preceq\kappa^2\lambda_{\min}(\hat G)^{-1}I$ (eigenvalues $\kappa^2\lambda/(\lambda+\kappa)^2\le\kappa^2/\lambda$), then (i).

\emph{(iii).} Conditionally, $\sqrt mS=m^{-1/2}\sum_{\calE}\hat F_{ij}\xi_{ij}$ is a sum of independent centered vectors with covariance $\hat B\to B_\infty$ and summands bounded by $\norm{\hat F_{ij}}\le\{\lambda_{\max}(\hat G)\hat F_{ij}^{\top}\hat G^{-1}\hat F_{ij}\}^{1/2}$, so the Lindeberg condition follows from $\max_{\calE}\hat F_{ij}^{\top}\hat G^{-1}\hat F_{ij}/m\to0$ and the Lindeberg--Feller theorem gives $\sqrt mS\Rightarrow N(0,B_\infty)$. The ridge term is $\sqrt m\kappa\norm{(\hat G+\kappa I)^{-1}\beta^{\circ}_{c}}\le\sqrt m\kappa\norm{\beta^{\circ}_{c}}/\lambda_{\min}(\hat G)\to0$ since $m\kappa\to0$, and $(\hat G+\kappa I)^{-1}\to G_\infty^{-1}$; Slutsky gives the limit. The working posterior covariance is $\nu(\hat G+\kappa I)^{-1}/m$, so the proof of Proposition~S2 applies with $\hat G,\hat M$ in place of $G,M$; its trace is $\nu d/m+o(1/m)$ against the conditional risk $\hat c/m+o(1/m)$ from (i)--(ii), which agree if and only if $\nu=\hat c/d$.

\emph{(iv).} With known agents take $\calE$ to be all dyads, $m=N$, and write $S=N^{-1}\sum F_{ij}(A_{ij}-\wproj)=E+T$ with $E=N^{-1}\sum F_{ij}\{A_{ij}-\wstar(U_i,U_j)\}$ and $T=N^{-1}\sum F_{ij}\delta(U_i,U_j)$. Then $\E[E\mid U]=0$, so $\E[S^{\top}G^{-1}S]=\E\{\E[E^{\top}G^{-1}E\mid U]\}+\E[T^{\top}G^{-1}T]$; the first term is $N^{-1}\int\wstar(1-\wstar)\Fv^{\top}G^{-1}\Fv$, and Hoeffding's decomposition of the $U$-statistic $T$ gives $\E[T^{\top}G^{-1}T]=N^{-1}\int\delta^2\Fv^{\top}G^{-1}\Fv+4(n-2)\kappa_1(\wstar)/\{n(n-1)\}$, the second part being the variance of its Hájek projection $(2/n)\sum_ig(U_i)$. Finally $T=\hat G_U(\beta^{\circ}_{c}(U)-\bproj)$ with $\hat G_U$ the empirical Gram matrix over all dyads, so $T$ is the conditional projection's deviation from the population projection, and the node-level term is the Hájek part of its variance. \qed

\subsection*{Proof of Proposition $2'$}

Throughout, residuals at $\bproj$ are $\xi_{ij}$ with $\E[\xi_{ij}\mid U]=\delta(U_i,U_j)$, $\Var(\xi_{ij}\mid U)=\wstar(1-\wstar)(U_i,U_j)$, and the $\xi_{ij}$ independent given $U$; features and residuals are bounded, $|\xi_{ij}|\le1+\norm{\bproj}_1$ and each entry of $\Fv$ in $[0,1]$.

\emph{(i).} $\E S=\int\delta\Fv=0$ by the normal equations, so $\Var(S)=\E[SS^{\top}]$. Expanding $N^2\E[SS^{\top}]$ over ordered pairs of dyads: the $N$ diagonal terms give $\E[\xi^2\Fv\Fv^{\top}]=M$ each; the $2N(n-2)$ ordered pairs sharing exactly one node give, by conditional independence, $\E[\delta(U,V)\delta(U,V')\Fv(U,V)\Fv(U,V')^{\top}]=\int g\,g^{\top}=\Sigma_g$ each; pairs sharing no node factor into two zero means. Dividing by $N^2$ gives the display, and taking the trace against $G^{-1}$ recovers Theorem~2(i).

\emph{(ii).} Write $\hat V(\beta)$ for the estimator with residuals at $\beta$, $T_1=\sum_ir_ir_i^{\top}$ and $T_2=\sum_{i<j}\xi_{ij}^2\Fv_{ij}\Fv_{ij}^{\top}$, so $N^2\hat V=T_1-T_2$. At $\bproj$: $\E[T_1]=\sum_i\sum_{j,k\neq i}\E[\xi_{ij}\xi_{ik}\Fv_{ij}\Fv_{ik}^{\top}]=n(n-1)M+n(n-1)(n-2)\Sigma_g$ and $\E[T_2]=NM$, hence $\E[\hat V(\bproj)]=M/N+\{4(n-2)/(n(n-1))\}\Sigma_g$, the display of (i). For consistency, $n\,T_2/N^2\le Cn\cdot N/N^2=O(1/n)$ deterministically by boundedness. Each entry of $T_1$ is a sum over $O(n^3)$ ordered triples $(i;j,k)$ of terms bounded by $C$; two such terms are correlated only if their index sets intersect, at most $O(n^5)$ ordered pairs of triples, so $\Var(T_{1,ab})\le Cn^5$ and $\Var(n\,T_{1,ab}/N^2)\le Cn^7/N^4=O(1/n)$. Chebyshev and the mean computation give $n\hat V(\bproj)\to M\cdot0+4\Sigma_g$ elementwise in probability, since $nM/N\to0$. For the plug-in, write $h=\bhat-\bproj$, $\norm h=O_p(n^{-1/2})$ by Theorem~2(iii) in either regime, $r_i(\bhat)=r_i(\bproj)-Q_ih$ with $Q_i=\sum_{j\neq i}\Fv_{ij}\Fv_{ij}^{\top}$, $\norm{Q_i}\le C n$. Then
\[
N^2\{\hat V(\bhat)-\hat V(\bproj)\}=-\sum_i\big(r_ih^{\top}Q_i^{\top}+Q_ihr_i^{\top}\big)+\sum_iQ_ihh^{\top}Q_i^{\top}-\{T_2(\bhat)-T_2(\bproj)\}.
\]
Since $\E\norm{r_i}^2\le Cn^2$, the first sum is $O_p(n\cdot n\cdot n\cdot n^{-1/2})=O_p(n^{5/2})$, the second $O_p(n\cdot n^2\cdot n^{-1})=O_p(n^{2})$, and the third $O(N\norm h)=O_p(n^{3/2})$ by boundedness; multiplying by $n/N^2=O(n^{-3})$ sends every term to zero in probability. Hence $n\hat V(\bhat)\to4\Sigma_g$ elementwise.

\emph{(iii).} On the Gram event of Theorem~2(ii), $A\to G$ in probability, so $n\,\ell^{\top}\Sighat_{\mathrm{dy}}\ell=n\,\ell^{\top}A^{-1}\hat VA^{-1}\ell\to4\,\ell^{\top}G^{-1}\Sigma_gG^{-1}\ell>0$. Theorem~2(iii) gives $\sqrt n\,\ell^{\top}(\bhat-\bproj)\Rightarrow N(0,\,4\ell^{\top}G^{-1}\Sigma_gG^{-1}\ell)$; Slutsky's theorem yields the studentized limit $N(0,1)$. Flooring the directional variance at zero changes nothing in the limit because the limiting variance is positive. \hfill$\square$

\subsection*{Proof of Theorem 3}

(i) On $\mathcal{E}$, Lemma~\ref{lem:decomp} and clipping give $\norm{\whatp-\wstar}^2\le\norm{\wstar-\wproj}^2+8S^{\top}G^{-1}S+8\kappa^2\norm{\bproj}^2/\lambda$; on $\mathcal{E}^c$ the clipped error is at most one. Taking expectations with $\E[S^{\top}G^{-1}S]=c(\wstar)/m$ and Lemma~\ref{lem:gram} proves the risk bound. For contraction, given the data, $w=\beta^{\top}\Fv$ with $\beta\sim N(\bhat,\Sighat)$ under $\Pi_m$, so $\E_{\Pi_m}\norm{w-\wproj}^2=(\bhat-\bproj)^{\top}G(\bhat-\bproj)+\tr(G\Sighat)$, and on $\mathcal{E}$, $\tr(G\Sighat)=(\nu/m)\tr(GA^{-1})\le2\nu d/m$ by $A^{-1}\preceq2G^{-1}$. Markov's inequality on $\mathcal{E}$, the bound $\Pi_m(\cdot)\le1$ on $\mathcal{E}^c$, and Lemmas~\ref{lem:decomp} and~\ref{lem:gram} complete the argument.

(ii) \emph{Upper bound.} For $\wstar\in\calW_{\calH}$, $\wproj=\wstar$, $\sigma^2=\wstar(1-\wstar)\le1/4$, so $c(\wstar)\le d/4$ and $8c(\wstar)/m\le2d/m$; and $\norm{\wstar}^2=(\bproj)^{\top}G\bproj\ge\lambda\norm{\bproj}^2$ with $\norm{\wstar}^2\le1$ gives $\norm{\bproj}^2\le1/\lambda$. Substituting into (i) with $\kappa^2=\nu^2/(m^2\tau^4)$ gives the bound.

\emph{Lower bound.} Write $D=d-1$. The set of orthonormal bases of $\calH_0$ is the image of one basis under the orthogonal group of $\R^D$, a compact set, and the map from a basis to $\sum_k\supn{\psi_k}$ is continuous because every element of $\calH_0$ is bounded; the infimum $B$ is attained by some basis $(\psi_k)_{k\le D}$, which we fix. Let $\eta=\{3/(128m)\}^{1/2}$ and, for $\theta\in\{-1,1\}^D$, $w_\theta=\tfrac12+\eta\sum_k\theta_k\psi_k$. Since $m\ge3B^2/8$ is equivalent to $\eta B\le1/4$, every $w_\theta$ takes values in $[1/4,3/4]$ and lies in $\calW_{\calH}$, and by orthonormality $\norm{w_\theta-w_{\theta'}}^2=4\eta^2\,\mathrm{ham}(\theta,\theta')$.

Let $\tilde w$ be any estimator and $\hat\theta$ a measurable minimizer of $\norm{\tilde w-w_\theta}$ over $\theta$; the triangle inequality gives $\norm{w_{\hat\theta}-w_\theta}\le2\norm{\tilde w-w_\theta}$, hence $\norm{\tilde w-w_\theta}^2\ge\eta^2\,\mathrm{ham}(\hat\theta,\theta)$, so
\[
\sup_{\theta}\E_\theta\norm{\tilde w-w_\theta}^2\ge\eta^2\sum_{k=1}^{D}2^{-D}\sum_\theta\Prob_\theta(\hat\theta_k\neq\theta_k)\ge\eta^2\sum_{k=1}^{D}\tfrac12\big(1-\max_\theta\TV(\Prob_\theta,\Prob_{\theta^{(k)}})\big),
\]
where $\theta^{(k)}$ flips the $k$th coordinate, by pairing $\theta$ with $\theta^{(k)}$ and using $P_1(\mathcal{V})+P_2(\mathcal{V}^c)\ge1-\TV(P_1,P_2)$ with $\mathcal{V}=\{\hat\theta_k=-1\}$. The two laws share the distribution of the positions, so the Kullback--Leibler divergence is $m$ times the expected divergence of one Bernoulli outcome: with $p,p'\in[1/4,3/4]$ and $|p-p'|=2\eta|\psi_k|$, $\KL(p\Vert p')\le(p-p')^2/\{p'(1-p')\}\le\frac{16}{3}(p-p')^2$, so $\KL(\Prob_\theta\Vert\Prob_{\theta^{(k)}})\le\frac{16}{3}\,m\cdot4\eta^2\int\psi_k^2=\frac{64}{3}m\eta^2$ and Pinsker's inequality gives $\TV\le\{\frac{32}{3}m\eta^2\}^{1/2}=\frac12$. Hence $\sup_\theta\E_\theta\norm{\tilde w-w_\theta}^2\ge\eta^2D/4=3(d-1)/(512m)$. Since the estimator was arbitrary and $\{w_\theta\}\subset\calW_{\calH}$, the lower bound follows; the middle inequality of the display in the theorem is immediate. \qed

\begin{propositionS}[Lipschitz transfer]
\label{prop:functionals}
For graphons $w,w'$ with values in $[0,1]$,
\begin{align*}
|e(w)-e(w')|&\le\norm{w-w'}, & \norm{d_w-d_{w'}}_{L_2[0,1]}&\le\norm{w-w'}, & |t(w)-t(w')|&\le3\norm{w-w'},\\
|\varsigma(w)-\varsigma(w')|&\le2\norm{w-w'}, & \big|\norm{w}_{\mathrm{op}}-\norm{w'}_{\mathrm{op}}\big|&\le\norm{w-w'}, & |C(w)-C(w')|&\le \frac{5\norm{w-w'}}{\varsigma(w)},
\end{align*}
the last when $\varsigma(w)>0$, with the convention $C(w')=0$ when $\varsigma(w')=0$. Consequently the plug-in estimator of any functional $\phi$ among these satisfies $\E|\phi(\whatp)-\phi(\wstar)|^2\le L_\phi^2\,\E\norm{\whatp-\wstar}^2$ with the stated Lipschitz constant ($5/\varsigma(\wstar)$ for the clustering coefficient), and the posterior of $\phi(w^{+})$ under $\Pi_m$, with $w^{+}$ the clipped graphon, concentrates around $\phi((\wproj)^{+})$: for every $r>0$, $\E\,\Pi_m(|\phi(w^{+})-\phi((\wproj)^{+})|>r)$ is bounded by the contraction bound of Theorem~3(i) evaluated at $r^2/L_\phi^2$, with $L_\phi=5/\varsigma((\wproj)^{+})$ for the clustering coefficient, which requires $\varsigma((\wproj)^{+})>0$.
\end{propositionS}

The clustering ratio is not Lipschitz, and not continuous, at vanishing wedge density, so the transfer for $C$ is local to $\{\varsigma\ge\varsigma_0>0\}$; the other functionals transfer globally.

\begin{proof}
All integrals are over unit cubes with Lebesgue measure, a probability measure, so $\int|h|\le(\int h^2)^{1/2}$. Edge density: $|e(w)-e(w')|\le\int|w-w'|\le\norm{w-w'}$. Degrees: by Jensen's inequality $\{d_w(x)-d_{w'}(x)\}^2\le\int(w-w')^2(x,y)\,dy$, and integrating over $x$ gives the bound. Triangles: the telescoping identity
\[
w_{xy}w_{yz}w_{xz}-w'_{xy}w'_{yz}w'_{xz}=(w-w')_{xy}w_{yz}w_{xz}+w'_{xy}(w-w')_{yz}w_{xz}+w'_{xy}w'_{yz}(w-w')_{xz}
\]
bounds each term in absolute value by $|w-w'|$ at one coordinate pair, and integrating gives $3\norm{w-w'}$. Wedges: $|d_w^2-d_{w'}^2|=|d_w-d_{w'}|(d_w+d_{w'})\le2|d_w-d_{w'}|$, so $|\varsigma(w)-\varsigma(w')|\le2\norm{d_w-d_{w'}}_{L_2[0,1]}\le2\norm{w-w'}$. Operator norm: the integral operators are self-adjoint and Hilbert--Schmidt, and $|\norm w_{\mathrm{op}}-\norm{w'}_{\mathrm{op}}|\le\norm{\mathcal{T}_{w-w'}}_{\mathrm{op}}\le\norm{\mathcal{T}_{w-w'}}_{\mathrm{HS}}=\norm{w-w'}$. Clustering: $t(w')\le\int w'_{xy}w'_{xz}=\varsigma(w')$ gives $C(w')\le1$, and
\[
|C(w)-C(w')|\le\frac{|t(w)-t(w')|+C(w')|\varsigma(w)-\varsigma(w')|}{\varsigma(w)}\le\frac{5\norm{w-w'}}{\varsigma(w)};
\]
if $\varsigma(w')=0$ then $w'=0$ almost everywhere, $C(w')=0$ by convention, and $|C(w)-C(w')|=t(w)/\varsigma(w)\le3\norm{w}/\varsigma(w)=3\norm{w-w'}/\varsigma(w)$. The plug-in bound applies the Lipschitz property with $w=\wstar$, $w'=\whatp$. For the posterior statement, clipping is a pointwise contraction, so $|\phi(w^{+})-\phi((\wproj)^{+})|\le L_\phi\norm{w^{+}-(\wproj)^{+}}\le L_\phi\norm{w-\wproj}$, where for the clustering coefficient the Lipschitz property is applied with first argument $(\wproj)^{+}$, and the event $\{|\phi(w^{+})-\phi((\wproj)^{+})|>r\}$ is contained in $\{\norm{w-\wproj}^2>r^2/L_\phi^2\}$, to which the contraction bound of Theorem~3(i) applies.
\end{proof}

\begin{propositionS}[Asymptotic calibration, independent dyads]
Let Assumption~1 hold and let the data follow (2). As $m\to\infty$,
\[
\sqrt m\,(\bhat-\bproj)\Rightarrow N\big(0,\,G^{-1}MG^{-1}\big),\qquad m\Sighat\to\nu G^{-1}\ \text{a.s.},\qquad m\,\E\norm{\what-\wproj}^2\to c(\wstar).
\]
For every linear functional $\ell^{\top}\beta$, and in particular for every link probability $\beta^{\top}\Fv(u,v)$, the frequentist coverage of the posterior $(1-\alpha)$ credible interval converges to a limit that is at least $1-\alpha$ for all $\ell$ if and only if $M\preceq \nu G$; the condition $\nu\ge\bar\sigma^2$ is sufficient, and $\nu=1/4$ suffices whenever $\wstar\in\calH$.
\end{propositionS}

The proposition is a Bernstein--von Mises statement for a misspecified working likelihood: the posterior covariance is $\nu G^{-1}/m$ while the sampling covariance of $\bhat$ is the sandwich $G^{-1}MG^{-1}/m$, so credible intervals are conservative or anti-conservative according to whether $\nu$ dominates the residual second moment in the operator order. All three readings are tested in Section~4: the safe choice $\nu=1/4$ is conservative; the plug-in choice used in Section~5 is calibrated in trace but can be anti-conservative for individual dyads, with coverage down to $0.88$; and under the one-graph design the intervals are valid only when $\kappa_1=0$, by Theorem~2(iii). Intervals here are for functionals of $\wproj$, and their coverage of $\wstar$ tends to zero whenever $\wstar\neq\wproj$. Guarantees transfer to smooth functionals such as triangle density and clustering through explicit Lipschitz constants, recorded as Proposition~S1.

\subsection*{Proof of Proposition S2}

By Lemma~\ref{lem:decomp}, $\sqrt m(\bhat-\bproj)=A^{-1}\{\sqrt m\,S-\sqrt m\,\kappa\bproj\}$. The strong law gives $\Ghat\to G$ almost surely and $A^{-1}\to G^{-1}$; $\sqrt m\,\kappa=\nu/(\sqrt m\,\tau^2)\to0$; and $\sqrt m\,S$ is a normalized sum of independent, bounded, mean-zero vectors with covariance $\E[\xi^2\Fv\Fv^{\top}]=M$, so the multivariate central limit theorem and Slutsky's lemma give $\sqrt m(\bhat-\bproj)\Rightarrow N(0,G^{-1}MG^{-1})$, and $m\Sighat=\nu A^{-1}\to\nu G^{-1}$ almost surely.

For the risk limit let $K_m=m\norm{\what-\wproj}^2=\{\sqrt m(\bhat-\bproj)\}^{\top}G\{\sqrt m(\bhat-\bproj)\}$, so $K_m\Rightarrow Z^{\top}GZ$ with $Z\sim N(0,G^{-1}MG^{-1})$ and $\E[Z^{\top}GZ]=\tr(G^{-1}M)=c(\wstar)$; it suffices that $(K_m)$ is uniformly integrable. On $\mathcal{E}$, $K_m\le8m\,S^{\top}G^{-1}S+8m\kappa^2\norm{\bproj}^2/\lambda$; the second term vanishes, and $m\,S^{\top}G^{-1}S\le\lambda^{-1}\norm{\sqrt m\,S}^2$ has bounded second moment because the summands $X_s=\xi_s\Fv_s$ are independent, mean zero and bounded by $c_0=(1+\supn{\wproj})\sqrt d$, so that only paired index patterns contribute to $\E\norm{\sum_sX_s}^4\le3m^2c_0^4$. On $\mathcal{E}^c$, $A^{-1}\preceq\kappa^{-1}I_d$ and $\norm{m^{-1}\sum\Fv_sY_s}\le\sqrt d$ give $\norm{\bhat}\le\sqrt d\,m\tau^2/\nu$, so $K_m$ is bounded by a polynomial in $m$, and Lemma~\ref{lem:gram} makes $\E[K_m\ind_{\mathcal{E}^c}]\to0$. Hence $\E K_m\to c(\wstar)$.

For coverage, the posterior $(1-\alpha)$ interval for $\ell^{\top}\beta$ is $\ell^{\top}\bhat\pm z_{1-\alpha/2}\{\ell^{\top}\Sighat\ell\}^{1/2}$, and by the two limits and Slutsky's lemma its coverage converges to
\[
2\Phi\Big(z_{1-\alpha/2}\Big\{\frac{\nu\,\ell^{\top}G^{-1}\ell}{\ell^{\top}G^{-1}MG^{-1}\ell}\Big\}^{1/2}\Big)-1,
\]
with the convention that the limit is one when the denominator vanishes. The limit is at least $1-\alpha$ for all $\ell$ if and only if $\nu G^{-1}\succeq G^{-1}MG^{-1}$, which after multiplication by $G$ on both sides is $\nu G\succeq M$; and $M=\int\sigma^2\Fv\Fv^{\top}\preceq\bar\sigma^2G$ makes $\nu\ge\bar\sigma^2$ sufficient, with $\bar\sigma^2\le1/4$ when $\wstar\in\calH$. \qed

\subsection*{Transfer to network functionals}

For a graphon $w$ write $e(w)=\int w$ for the edge density, $d_w(x)=\int w(x,y)\,dy$ for the degree function, $t(w)=\int w(x,y)w(y,z)w(x,z)\,dx\,dy\,dz$ for the triangle density, $\varsigma(w)=\int d_w(x)^2dx$ for the wedge density, $C(w)=t(w)/\varsigma(w)$ for the clustering coefficient when $\varsigma(w)>0$, and $\norm{w}_{\mathrm{op}}$ for the operator norm of the integral operator with kernel $w$.

\section{Additional tables}
\label{sec:tables}

Table~\ref{tab:agents} reports the single de-thinned agents beside the calibrated selection, least-squares and enlarged-model comparators on all three criteria. Table~\ref{tab:single} is the fitted-agent single-graph pipeline of Figure~S1(c). Table~\ref{tab:functionals} reports implied triangle densities and clustering coefficients; on de-thinned agents no combination dominates the single agents, and the apparent recovery reported by the uncorrected protocol is absent. Table~\ref{tab:thinning} is the artifact comparison of Section~5. Table~\ref{tab:onegraphtab} reports the one-graph experiment behind Figure~S4: the Monte Carlo mean of $S^{\top}G^{-1}S$ beside the exact identity of Theorem~2(i), and the coverage of both interval constructions. Table~\ref{tab:coverage} is the coverage study under independent dyads. Table~\ref{tab:weightsupp} reports the posterior-mean synthesis weights of Section~5. Table~\ref{tab:external} places the external topological predictors beside the synthesis rules on the identical held-out dyads. Table~\ref{tab:retention} repeats the paired span-over-selection gain with the evaluation dyads held fixed and the training edges subsampled to $30$, $50$ and $70$ per cent. Table~\ref{tab:heterosum} summarizes the heterogeneous-information experiment of Section~5, and Table~\ref{tab:heterofull} reports all its rules and criteria. Figure~S1 is the misspecification, collinearity and fitted-agent-pipeline simulation of Section~4.

\begin{table}[ht]
\centering
\caption{Six networks, de-thinned agents: single agents and additional comparators; best per row in bold. Mean (s.e.) over ten splits.}
\label{tab:agents}
\resizebox{\textwidth}{!}{\begin{tabular}{llcccccccc}
\toprule
Network & score & ER & Chung--Lu & Degree blocks & SBM & RDPG & Select+Platt & Least squares & RDPG-16 \\
\midrule
polblogs & Brier$/e(\mathbf{A})$ & 0.978 & 0.861 & 0.869 & 0.945 & 0.702 & 0.708 & \textbf{0.700} & 0.717 \\
 & log score$/H_0$ & 1.000 & 0.755 & 0.737 & 0.884 & 0.662 & \textbf{0.613} & 0.638 & 0.746 \\
 & AUC & 0.500 & 0.899 & 0.892 & 0.783 & 0.921 & 0.921 & \textbf{0.930} & 0.899 \\
email-Eu-core & Brier$/e(\mathbf{A})$ & 0.967 & 0.833 & 0.855 & 0.904 & 0.673 & 0.680 & 0.668 & \textbf{0.637} \\
 & log score$/H_0$ & 1.000 & 0.780 & 0.780 & 0.860 & 0.659 & \textbf{0.616} & 0.625 & 0.666 \\
 & AUC & 0.500 & 0.859 & 0.854 & 0.785 & 0.919 & 0.920 & \textbf{0.931} & 0.919 \\
Cora & Brier$/e(\mathbf{A})$ & 0.999 & 0.997 & 0.999 & 0.996 & 1.006 & \textbf{0.995} & 0.997 & 1.007 \\
 & log score$/H_0$ & 1.000 & 1.119 & 1.095 & 0.979 & 1.480 & 0.960 & \textbf{0.917} & 1.329 \\
 & AUC & 0.500 & 0.622 & 0.608 & \textbf{0.759} & 0.646 & 0.631 & 0.754 & 0.685 \\
ca-GrQc & Brier$/e(\mathbf{A})$ & 0.999 & 0.980 & 0.990 & 0.989 & 0.772 & 0.777 & 0.768 & \textbf{0.732} \\
 & log score$/H_0$ & 1.000 & 0.975 & 0.965 & 0.926 & 1.169 & 0.767 & \textbf{0.726} & 0.960 \\
 & AUC & 0.500 & 0.737 & 0.740 & 0.807 & 0.733 & 0.683 & \textbf{0.809} & 0.789 \\
wiki-Vote & Brier$/e(\mathbf{A})$ & 0.996 & 0.914 & 0.944 & 0.990 & 0.802 & 0.803 & 0.801 & \textbf{0.781} \\
 & log score$/H_0$ & 1.000 & 0.690 & 0.715 & 0.917 & 0.595 & \textbf{0.557} & 0.596 & 0.598 \\
 & AUC & 0.500 & 0.933 & 0.928 & 0.775 & 0.947 & 0.949 & \textbf{0.949} & 0.942 \\
soc-Epinions1 & Brier$/e(\mathbf{A})$ & 1.000 & 0.971 & 0.995 & 1.000 & 0.925 & 0.926 & 0.926 & \textbf{0.919} \\
 & log score$/H_0$ & 1.000 & 0.707 & 0.778 & 0.965 & 0.702 & \textbf{0.649} & 0.701 & 0.716 \\
 & AUC & 0.500 & \textbf{0.886} & 0.881 & 0.759 & 0.878 & 0.876 & 0.877 & 0.868 \\
\bottomrule
\end{tabular}
}
\end{table}

\begin{table}[ht]
\centering
\caption{Single graphs with fitted, de-thinned agents (Figure~S1(c)): exact risk over density, mean (s.e.) over $R=20$ graphs; best per row in bold.}
\label{tab:single}
\resizebox{\textwidth}{!}{\begin{tabular}{lrcccccccc}
\toprule
$n$ & density & Select & Select+affine & Hull (Brier) & Hull (log) & Cone & Logit stacking & BPS (span) & RDPG-16 \\
\midrule
400 & 0.095 & 0.0591 (0.0018) & 0.0572 (0.0016) & 0.0457 (0.0009) & 0.0461 (0.0008) & 0.0458 (0.0008) & 0.0489 (0.0010) & \textbf{0.0456 (0.0009)} & 0.1102 (0.0010) \\
800 & 0.095 & 0.0361 (0.0005) & 0.0357 (0.0005) & 0.0266 (0.0003) & 0.0273 (0.0002) & 0.0267 (0.0003) & 0.0305 (0.0003) & \textbf{0.0265 (0.0003)} & 0.0645 (0.0006) \\
1600 & 0.095 & 0.0231 (0.0001) & 0.0229 (0.0001) & 0.0185 (0.0001) & 0.0193 (0.0001) & 0.0185 (0.0001) & 0.0231 (0.0003) & \textbf{0.0184 (0.0001)} & 0.0352 (0.0003) \\
3200 & 0.095 & 0.0172 (0.0000) & 0.0172 (0.0000) & 0.0150 (0.0000) & 0.0158 (0.0001) & 0.0149 (0.0000) & 0.0198 (0.0001) & \textbf{0.0148 (0.0000)} & 0.0181 (0.0001) \\
\bottomrule
\end{tabular}
}
\end{table}

\begin{table}[ht]
\centering
\caption{Implied triangle density and clustering coefficient of the fitted rules on de-thinned agents, against the observed values; closest per row in bold.}
\label{tab:functionals}
\resizebox{\textwidth}{!}{\begin{tabular}{llcccccccccc}
\toprule
Network & functional & observed & Chung--Lu & Degree blocks & SBM & RDPG & Select+affine & RDPG-16 & Hull (Brier) & Logit stacking & BPS \\
\midrule
polblogs & $t$ (triangles) & 3.31e-04 & 2.48e-04 & 1.71e-04 & 3.56e-05 & 3.27e-04 & 3.19e-04 & 3.34e-04 & 3.10e-04 & 3.17e-04 & \textbf{3.31e-04} \\
 & $C$ (clustering) & 0.226 & 0.173 & 0.131 & 0.061 & 0.194 & 0.197 & 0.184 & 0.188 & \textbf{0.217} & 0.197 \\
email-Eu-core & $t$ (triangles) & 6.62e-04 & 4.03e-04 & 3.53e-04 & 8.49e-05 & 6.14e-04 & 5.98e-04 & \textbf{6.51e-04} & 5.37e-04 & 5.77e-04 & 6.17e-04 \\
 & $C$ (clustering) & 0.267 & 0.164 & 0.153 & 0.073 & 0.219 & 0.220 & 0.215 & 0.203 & \textbf{0.225} & 0.217 \\
Cora & $t$ (triangles) & 4.93e-07 & 6.38e-08 & 1.48e-08 & 2.70e-08 & 1.14e-07 & 9.22e-08 & \textbf{1.97e-07} & 3.80e-08 & 3.90e-08 & 1.31e-07 \\
 & $C$ (clustering) & 0.093 & 0.011 & 0.004 & 0.012 & 0.045 & 0.021 & \textbf{0.047} & 0.017 & 0.011 & 0.018 \\
ca-GrQc & $t$ (triangles) & 2.01e-06 & 3.61e-08 & 4.14e-08 & 3.71e-08 & 1.03e-06 & 1.22e-06 & 1.11e-06 & 1.03e-06 & 1.20e-06 & \textbf{1.23e-06} \\
 & $C$ (clustering) & 0.630 & 0.010 & 0.015 & 0.026 & \textbf{0.499} & 0.342 & 0.376 & 0.499 & 0.274 & 0.323 \\
wiki-Vote & $t$ (triangles) & 1.01e-05 & 7.59e-06 & 3.57e-06 & 1.53e-07 & \textbf{1.01e-05} & 9.78e-06 & 1.03e-05 & 9.67e-06 & 8.96e-06 & 9.63e-06 \\
 & $C$ (clustering) & 0.125 & 0.095 & 0.058 & 0.009 & 0.104 & 0.105 & 0.098 & 0.102 & \textbf{0.116} & 0.105 \\
soc-Epinions1 & $t$ (triangles) & 2.23e-08 & 1.08e-08 & 6.37e-10 & 1.34e-11 & 1.17e-08 & 1.08e-08 & \textbf{1.38e-08} & 1.05e-08 & 6.70e-09 & 1.20e-08 \\
 & $C$ (clustering) & 0.066 & \textbf{0.032} & 0.006 & 0.001 & 0.028 & 0.028 & 0.030 & 0.027 & 0.031 & 0.029 \\
\bottomrule
\end{tabular}
}
\end{table}

\begin{table}[ht]
\centering
\caption{The edge-thinning artifact: dominant weight, weight sum and Brier gain of span synthesis over selection, with agents fitted to the thinned graph and with de-thinned agents; mean (s.e.) over ten splits.}
\label{tab:thinning}
\resizebox{\textwidth}{!}{\begin{tabular}{lcccccc}
\toprule
Network & \multicolumn{3}{c}{agents fitted to the thinned graph} & \multicolumn{3}{c}{de-thinned agents} \\
 & RDPG weight & $\sum_j\beta_j$ & BPS gain over Select (\%) & RDPG weight & $\sum_j\beta_j$ & BPS gain over Select (\%) \\
\midrule
polblogs & 1.34 & 1.59 & 2.70 (0.09) & 0.96 & 1.12 & 0.22 (0.03) \\
email-Eu-core & 1.25 & 1.85 & 3.61 (0.06) & 0.90 & 1.28 & 0.53 (0.08) \\
Cora & 0.92 & 1.87 & 0.25 (0.19) & 0.65 & 1.33 & 0.97 (0.14) \\
ca-GrQc & 1.52 & 2.07 & 2.84 (0.20) & 1.10 & 1.53 & 0.36 (0.07) \\
wiki-Vote & 1.42 & 1.50 & 1.78 (0.08) & 1.00 & 1.05 & 0.14 (0.01) \\
soc-Epinions1 & 1.35 & 1.47 & 0.50 (0.17) & 0.96 & 1.00 & -0.18 (0.09) \\
\bottomrule
\end{tabular}
}
\end{table}

\begin{table}[ht]
\centering
\caption{One-graph experiment (Figure~S4): Monte Carlo mean of the score quadratic against the exact identity of Theorem~2(i), estimation and clipped risks, and per-link coverage of nominal $95\%$ intervals at the three probed dyads; $R=60$ (Monte Carlo s.e.\ about $0.03$ per coverage entry). The third block is the off-span truth with $g\equiv0$ of Section~4, for which the directional-variance floor of (6) never bound.}
\label{tab:onegraphtab}
\resizebox{\textwidth}{!}{\begin{tabular}{lrcccccc}
\toprule
truth & $n$ & $\mathbb{E}\,S^{\top}G^{-1}S$ (MC) & identity of Theorem~2(i) & $\mathbb{E}\|\hat w-w^\circ\|^2$ & $\mathbb{E}\|\hat w^{+}-w^\star\|^2$ & coverage, posterior & coverage, dyadic \\
\midrule
in the span & 50 & 4.45e-04 (5.9e-05) & 4.17e-04 & 4.38e-04 & 4.37e-04 & 0.97, 0.98, 0.95 & 0.92, 0.90, 0.83 \\
 & 100 & 1.16e-04 (1.3e-05) & 1.03e-04 & 1.22e-04 & 1.22e-04 & 0.97, 0.98, 0.93 & 0.93, 0.93, 0.82 \\
 & 200 & 2.80e-05 (2.8e-06) & 2.57e-05 & 2.73e-05 & 2.73e-05 & 0.97, 0.97, 0.97 & 0.90, 0.92, 0.92 \\
 & 400 & 6.40e-06 (8.3e-07) & 6.40e-06 & 6.88e-06 & 6.88e-06 & 0.97, 1.00, 0.97 & 0.93, 0.93, 0.97 \\
 & 800 & 1.37e-06 (1.6e-07) & 1.60e-06 & 1.41e-06 & 1.41e-06 & 1.00, 1.00, 0.95 & 1.00, 0.93, 0.93 \\
 & 1600 & 4.94e-07 (6.1e-08) & 3.99e-07 & 4.92e-07 & 4.92e-07 & 0.90, 0.98, 0.98 & 0.90, 0.93, 0.95 \\
superposition & 50 & 5.02e-04 (5.3e-05) & 4.98e-04 & 5.35e-04 & 1.35e-03 & 0.98, 0.95, 0.98 & 0.88, 0.87, 0.83 \\
 & 100 & 1.31e-04 (1.3e-05) & 1.35e-04 & 1.42e-04 & 9.76e-04 & 0.93, 0.95, 0.95 & 0.93, 0.92, 0.92 \\
 & 200 & 4.73e-05 (5.9e-06) & 3.96e-05 & 4.77e-05 & 8.87e-04 & 0.92, 0.92, 0.88 & 0.93, 0.92, 0.90 \\
 & 400 & 1.33e-05 (1.5e-06) & 1.28e-05 & 1.39e-05 & 8.54e-04 & 0.90, 0.90, 0.88 & 0.95, 0.95, 0.88 \\
 & 800 & 3.72e-06 (4.4e-07) & 4.68e-06 & 3.82e-06 & 8.44e-04 & 0.85, 0.83, 0.85 & 0.98, 0.93, 0.98 \\
 & 1600 & 2.10e-06 (2.8e-07) & 1.91e-06 & 2.11e-06 & 8.43e-04 & 0.77, 0.73, 0.73 & 0.95, 0.93, 0.90 \\
off-span, $g=0$ & 100 & 1.25e-04 (1.6e-05) & 1.49e-04 & 1.22e-04 & 1.02e-03 & 0.98, 0.95, 0.98 & 0.98, 0.92, 0.98 \\
 & 200 & 3.70e-05 (4.1e-06) & 3.70e-05 & 3.96e-05 & 9.40e-04 & 0.92, 0.98, 0.90 & 0.90, 0.95, 0.88 \\
 & 400 & 8.20e-06 (1.0e-06) & 9.23e-06 & 8.17e-06 & 9.08e-04 & 0.98, 0.90, 0.98 & 0.98, 0.90, 0.98 \\
 & 800 & 2.23e-06 (3.1e-07) & 2.30e-06 & 2.20e-06 & 9.02e-04 & 0.97, 0.97, 0.97 & 0.93, 0.98, 0.97 \\
 & 1600 & 5.69e-07 (5.8e-08) & 5.76e-07 & 5.72e-07 & 9.01e-04 & 0.95, 0.97, 0.93 & 0.93, 0.97, 0.93 \\
\bottomrule
\end{tabular}
}
\end{table}

\begin{table}[ht]
\centering
\caption{Coverage of nominal $95\%$ link-probability intervals under independent dyads, for a truth in the span and the superposition truth, with $\nu=1/4$ and with the plug-in $\nu$ of Section~5, beside the pooled posterior spread; $R=400$.}
\label{tab:coverage}
\begin{tabular}{lrcccc}
\toprule
truth & $m$ & coverage, $\nu=1/4$ & spread & coverage, plug-in $\nu$ & spread \\
\midrule
in the span & 400 & 0.96, 1.00, 0.96 & 0.002 & 0.88, 0.97, 0.90 & 0.001 \\
 & 1600 & 0.96, 0.99, 0.96 & 0.000 & 0.88, 0.97, 0.90 & 0.000 \\
 & 6400 & 0.96, 0.99, 0.97 & 0.000 & 0.91, 0.97, 0.89 & 0.000 \\
superposition & 400 & 0.97, 0.96, 0.98 & 0.002 & 0.96, 0.94, 0.97 & 0.002 \\
 & 1600 & 0.95, 0.98, 0.98 & 0.000 & 0.94, 0.96, 0.96 & 0.000 \\
 & 6400 & 0.96, 0.98, 0.98 & 0.000 & 0.94, 0.97, 0.97 & 0.000 \\
\bottomrule
\end{tabular}

\end{table}

\begin{table}[ht]
\centering
\caption{Posterior-mean synthesis weights on de-thinned agents, averaged over ten splits (posterior s.d.\ in parentheses), their sum, the selected agent, and the affine recalibration slope of the selected agent.}
\label{tab:weightsupp}
\resizebox{\textwidth}{!}{\begin{tabular}{lcccccccc}
\toprule
Network & intercept & Chung--Lu & Degree blocks & SBM & RDPG & $\sum_j \beta_j$ & selected agent & slope $\hat\gamma$ \\
\midrule
polblogs & -0.00 (0.00) & 0.15 (0.03) & -0.12 (0.04) & 0.13 (0.03) & 0.96 (0.02) & 1.12 & RDPG (10/10) & 1.00 \\
email-Eu-core & -0.01 (0.00) & 0.18 (0.04) & -0.08 (0.04) & 0.28 (0.03) & 0.90 (0.02) & 1.28 & RDPG (10/10) & 1.00 \\
Cora & 0.00 (0.00) & 0.70 (0.19) & -0.51 (0.31) & 0.48 (0.13) & 0.65 (0.12) & 1.33 & RDPG (7/10) & 0.86 \\
ca-GrQc & 0.00 (0.00) & -0.15 (0.11) & 0.24 (0.10) & 0.34 (0.07) & 1.10 (0.02) & 1.53 & RDPG (10/10) & 1.11 \\
wiki-Vote & -0.00 (0.00) & 0.06 (0.02) & -0.14 (0.02) & 0.13 (0.03) & 1.00 (0.01) & 1.05 & RDPG (10/10) & 0.99 \\
soc-Epinions1 & -0.00 (0.00) & 0.09 (0.01) & -0.16 (0.02) & 0.11 (0.07) & 0.96 (0.01) & 1.00 & RDPG (10/10) & 0.99 \\
\bottomrule
\end{tabular}
}
\end{table}

\begin{table}[ht]
\centering
\caption{Six networks, external comparators on the identical held-out dyads: common neighbors, Adamic--Adar, Jaccard and preferential attachment, each Platt-calibrated on the validation dyads (the AUC equals that of the raw score), a logistic regression on seven topological features, and a random forest on the same features ($200$ trees, minimum leaf $20$, training non-edges capped at $2\times10^5$), beside four synthesis rules of Table~1; best per row in bold. The two calibrated columns refit the forest and the logistic with a nested validation split, half to train and half to Platt-adjust under the stratum weights; the raw forest's soc-Epinions1 Brier is the case-control artifact discussed in Section~5. The learned embedding predictor of the paragraph on node2vec is run on all nine networks and reported there, not in this table.}
\label{tab:external}
\resizebox{\textwidth}{!}{\begin{tabular}{llcccccccccccc}
\toprule
Network & score & CN & Adamic--Adar & Jaccard & Pref.\ attach. & Logit (heur.) & RF (topol.) & RF (cal.) & Logit (h., cal.) & Select+affine & BPS (span) & Logit stacking & RDPG-16 \\
\midrule
polblogs & Brier$/e(\mathbf{A})$ & 0.758 (0.002) & 0.765 (0.002) & 0.985 (0.003) & 0.851 (0.002) & 0.744 (0.002) & 0.745 (0.002) & 0.745 (0.002) & 0.745 (0.003) & 0.701 (0.002) & \textbf{0.700 (0.002)} & 0.710 (0.002) & 0.717 (0.004) \\
 & log score$/H_0$ & 0.645 (0.002) & 0.656 (0.002) & 0.897 (0.004) & 0.720 (0.002) & 0.619 (0.002) & 0.616 (0.002) & 0.616 (0.003) & 0.620 (0.003) & 0.670 (0.004) & 0.638 (0.004) & \textbf{0.590 (0.002)} & 0.746 (0.004) \\
 & AUC & 0.902 (0.001) & 0.906 (0.001) & 0.861 (0.001) & 0.899 (0.001) & 0.926 (0.001) & 0.927 (0.001) & 0.926 (0.001) & 0.925 (0.001) & 0.923 (0.001) & 0.930 (0.001) & \textbf{0.933 (0.001)} & 0.899 (0.001) \\
email-Eu-core & Brier$/e(\mathbf{A})$ & 0.681 (0.002) & 0.657 (0.002) & 0.746 (0.001) & 0.839 (0.001) & 0.626 (0.002) & \textbf{0.615 (0.002)} & 0.618 (0.002) & 0.629 (0.002) & 0.673 (0.002) & 0.668 (0.002) & 0.676 (0.002) & 0.637 (0.003) \\
 & log score$/H_0$ & 0.583 (0.001) & 0.570 (0.002) & 0.674 (0.001) & 0.770 (0.002) & 0.547 (0.001) & \textbf{0.537 (0.002)} & 0.538 (0.002) & 0.549 (0.001) & 0.666 (0.005) & 0.625 (0.004) & 0.586 (0.002) & 0.666 (0.006) \\
 & AUC & 0.928 (0.001) & 0.934 (0.001) & 0.916 (0.001) & 0.859 (0.001) & 0.939 (0.001) & \textbf{0.941 (0.001)} & 0.941 (0.001) & 0.939 (0.001) & 0.920 (0.001) & 0.931 (0.001) & 0.932 (0.001) & 0.919 (0.002) \\
Cora & Brier$/e(\mathbf{A})$ & 0.975 (0.013) & 0.985 (0.025) & 1.009 (0.016) & 0.997 (0.000) & 1.033 (0.038) & 1.081 (0.042) & \textbf{0.953 (0.010)} & 1.022 (0.030) & 1.004 (0.004) & 0.997 (0.003) & 0.987 (0.002) & 1.007 (0.008) \\
 & log score$/H_0$ & 0.857 (0.006) & 0.856 (0.012) & 0.935 (0.008) & 0.978 (0.001) & 0.866 (0.029) & 0.844 (0.015) & \textbf{0.807 (0.005)} & 0.865 (0.025) & 0.966 (0.008) & 0.917 (0.006) & 0.884 (0.002) & 1.329 (0.007) \\
 & AUC & 0.664 (0.003) & 0.664 (0.003) & 0.663 (0.003) & 0.622 (0.003) & 0.768 (0.004) & \textbf{0.786 (0.003)} & 0.782 (0.004) & 0.766 (0.004) & 0.655 (0.014) & 0.754 (0.003) & 0.759 (0.004) & 0.685 (0.004) \\
ca-GrQc & Brier$/e(\mathbf{A})$ & 0.653 (0.014) & \textbf{0.550 (0.014)} & 0.739 (0.008) & 0.988 (0.000) & 0.560 (0.012) & 0.726 (0.051) & 0.560 (0.007) & 0.573 (0.020) & 0.770 (0.003) & 0.768 (0.003) & 0.788 (0.003) & 0.732 (0.004) \\
 & log score$/H_0$ & 0.504 (0.012) & 0.459 (0.013) & 0.611 (0.009) & 0.900 (0.001) & 0.456 (0.012) & 0.517 (0.021) & \textbf{0.437 (0.005)} & 0.468 (0.018) & 0.765 (0.002) & 0.726 (0.004) & 0.722 (0.002) & 0.960 (0.010) \\
 & AUC & 0.863 (0.001) & 0.863 (0.001) & 0.863 (0.001) & 0.737 (0.002) & 0.886 (0.004) & \textbf{0.923 (0.001)} & 0.920 (0.002) & 0.892 (0.002) & 0.733 (0.002) & 0.809 (0.005) & 0.821 (0.002) & 0.789 (0.002) \\
wiki-Vote & Brier$/e(\mathbf{A})$ & 0.872 (0.001) & 0.896 (0.002) & 1.035 (0.005) & 0.913 (0.001) & 0.867 (0.002) & 0.937 (0.003) & 0.866 (0.002) & 0.868 (0.002) & 0.802 (0.002) & 0.801 (0.002) & 0.803 (0.002) & \textbf{0.781 (0.002)} \\
 & log score$/H_0$ & 0.643 (0.001) & 0.677 (0.001) & 0.951 (0.001) & 0.677 (0.001) & 0.616 (0.001) & 0.635 (0.002) & 0.607 (0.001) & 0.616 (0.001) & 0.594 (0.001) & 0.596 (0.002) & \textbf{0.549 (0.001)} & 0.598 (0.002) \\
 & AUC & 0.909 (0.001) & 0.910 (0.001) & 0.899 (0.000) & 0.933 (0.000) & 0.953 (0.000) & 0.955 (0.000) & 0.954 (0.000) & 0.953 (0.000) & 0.951 (0.000) & 0.949 (0.001) & \textbf{0.955 (0.000)} & 0.942 (0.000) \\
soc-Epinions1 & Brier$/e(\mathbf{A})$ & 0.959 (0.007) & 0.972 (0.009) & 1.003 (0.001) & 0.969 (0.001) & 0.932 (0.007) & 7.004 (0.067) & 0.937 (0.010) & 0.930 (0.005) & 0.926 (0.003) & 0.926 (0.003) & 0.940 (0.002) & \textbf{0.919 (0.002)} \\
 & log score$/H_0$ & 0.679 (0.002) & 0.721 (0.003) & 0.978 (0.000) & 0.693 (0.000) & 0.605 (0.003) & 2.645 (0.019) & \textbf{0.598 (0.004)} & 0.605 (0.003) & 0.702 (0.005) & 0.701 (0.005) & 0.645 (0.001) & 0.716 (0.002) \\
 & AUC & 0.838 (0.000) & 0.838 (0.000) & 0.837 (0.000) & 0.886 (0.000) & 0.951 (0.001) & \textbf{0.964 (0.000)} & 0.964 (0.000) & 0.951 (0.001) & 0.880 (0.000) & 0.877 (0.003) & 0.909 (0.000) & 0.868 (0.000) \\
\bottomrule
\end{tabular}
}
\end{table}

\begin{table}[ht]
\centering
\caption{Paired gain of the span synthesis over the calibrated selected agent, in per cent of the density, with the validation ($10\%$) and test ($20\%$) dyads held fixed at the main protocol's and the training edges subsampled to $30$, $50$ and $70\%$ of all edges; mean and paired $t$ interval over ten splits.}
\label{tab:retention}
\begin{tabular}{lccc}
\toprule
Network & training $30\%$ & training $50\%$ & training $70\%$ \\
\midrule
polblogs & 0.96 $[0.86,1.06]$ & 0.38 $[0.30,0.45]$ & 0.16 $[0.10,0.22]$ \\
email-Eu-core & 3.20 $[2.87,3.53]$ & 0.89 $[0.84,0.95]$ & 0.52 $[0.38,0.67]$ \\
Cora & 0.06 $[-0.24,0.36]$ & 0.39 $[0.06,0.73]$ & 0.64 $[0.52,0.77]$ \\
ca-GrQc & 0.17 $[0.12,0.23]$ & 0.10 $[0.08,0.13]$ & 0.13 $[0.11,0.16]$ \\
wiki-Vote & 0.20 $[0.15,0.26]$ & 0.09 $[0.06,0.11]$ & 0.11 $[0.09,0.13]$ \\
soc-Epinions1 & 0.72 $[-0.10,1.54]$ & -0.03 $[-0.07,0.02]$ & -0.03 $[-0.07,0.01]$ \\
\bottomrule
\end{tabular}

\end{table}

\begin{table}[ht]
\centering
\caption{Heterogeneous information, summary: training edges in five disjoint private shares, the class-to-share assignment rotated across seeds, each mechanism class fitted on its own share only and de-thinned by its own retention of about $0.14$; all rules fitted on the validation dyads. The selected column gives the modal selected agent (count of ten). Brier over density, mean (s.e.) over ten splits, best combination per row in bold; the single-site column is the calibrated selected agent of the pooled protocol, and the last column the paired gain of the span synthesis over the calibrated selected agent with its $95\%$ interval.}
\label{tab:heterosum}
\resizebox{\textwidth}{!}{\begin{tabular}{lccccccccc}
\toprule
Network & Select & Select+affine & Hull (Brier) & Hull (log) & BPS (span) & Logit stacking & selected & single site & BPS gain (\% of $e(\mathbf{A})$) \\
\midrule
polblogs & 0.887 (0.003) & 0.877 (0.003) & 0.820 (0.002) & 0.822 (0.002) & 0.819 (0.002) & \textbf{0.812 (0.002)} & CL (5) & 0.701 (0.002) & 5.8 $[5.1,6.5]$ \\
email-Eu-core & 0.858 (0.002) & 0.855 (0.001) & 0.821 (0.001) & 0.822 (0.001) & \textbf{0.809 (0.001)} & 0.810 (0.002) & CL (10) & 0.673 (0.002) & 4.6 $[4.4,4.9]$ \\
Cora & 1.011 (0.011) & 1.021 (0.015) & 1.003 (0.003) & \textbf{0.997 (0.000)} & 1.037 (0.016) & 1.000 (0.003) & ER (7) & 1.004 (0.004) & -1.5 $[-4.1,1.1]$ \\
ca-GrQc & 0.907 (0.009) & 0.904 (0.008) & 0.878 (0.005) & 0.911 (0.003) & 0.875 (0.006) & \textbf{0.859 (0.008)} & RDPG (10) & 0.770 (0.003) & 2.9 $[2.1,3.6]$ \\
wiki-Vote & 0.916 (0.003) & 0.899 (0.002) & 0.882 (0.002) & 0.884 (0.001) & 0.880 (0.002) & \textbf{0.872 (0.002)} & RDPG (9) & 0.802 (0.002) & 1.9 $[1.4,2.4]$ \\
soc-Epinions1 & 0.973 (0.002) & 0.972 (0.002) & 0.961 (0.002) & 0.963 (0.001) & 0.962 (0.002) & \textbf{0.957 (0.002)} & CL (9) & 0.926 (0.003) & 1.0 $[0.7,1.3]$ \\
\bottomrule
\end{tabular}
}
\end{table}

\begin{table}[ht]
\centering
\caption{Heterogeneous information, all rules: Brier over density, log score over $H_0$ and AUC, mean over ten splits; best per row in bold. The cone, least-squares and span-synthesis columns agree to the displayed precision in every cell (at $\tau^2=100$ the prior is numerically inert) and are shown once. The three bagging columns fit the same rank-eight spectral algorithm on each of the five shares and combine the copies by the same rules, the same-algorithm control of Section~5.}
\label{tab:heterofull}
\resizebox{\textwidth}{!}{\begin{tabular}{llcccccccccc}
\toprule
Network & score & Select & Select+affine & Select+Platt & Hull (Brier) & Hull (log) & Cone$=$LS$=$BPS & Logit stacking & Bag: Select & Bag: Sel+aff & Bag: BPS \\
\midrule
polblogs & Brier$/e(\mathbf{A})$ & 0.887 & 0.877 & 0.887 & 0.820 & 0.822 & 0.819 & 0.812 & 0.909 & 0.856 & \textbf{0.767} \\
 & log score$/H_0$ & 0.922 & 0.768 & 0.792 & 0.709 & 0.700 & 0.749 & \textbf{0.697} & 1.215 & 0.774 & 0.720 \\
 & AUC & 0.857 & 0.857 & 0.854 & 0.897 & 0.897 & 0.897 & 0.896 & 0.822 & 0.822 & \textbf{0.909} \\
email-Eu-core & Brier$/e(\mathbf{A})$ & 0.858 & 0.855 & 0.868 & 0.821 & 0.822 & \textbf{0.809} & 0.810 & 0.968 & 0.890 & 0.816 \\
 & log score$/H_0$ & 0.915 & 0.799 & 0.820 & 0.739 & \textbf{0.738} & 0.795 & 0.740 & 1.503 & 0.867 & 0.763 \\
 & AUC & 0.831 & 0.831 & 0.831 & 0.877 & \textbf{0.878} & 0.877 & 0.872 & 0.751 & 0.751 & 0.861 \\
Cora & Brier$/e(\mathbf{A})$ & 1.011 & 1.021 & 1.002 & 1.003 & \textbf{0.997} & 1.037 & 1.000 & 1.044 & 1.052 & 1.116 \\
 & log score$/H_0$ & 1.260 & 1.018 & 1.001 & 1.073 & \textbf{0.983} & 1.027 & 0.992 & 1.856 & 1.066 & 1.108 \\
 & AUC & 0.501 & 0.501 & 0.501 & 0.565 & \textbf{0.567} & 0.565 & 0.563 & 0.504 & 0.504 & 0.516 \\
ca-GrQc & Brier$/e(\mathbf{A})$ & 0.907 & 0.904 & 0.903 & 0.878 & 0.911 & 0.875 & 0.859 & 0.888 & 0.889 & \textbf{0.837} \\
 & log score$/H_0$ & 1.490 & 0.895 & 0.885 & 1.024 & 0.824 & 0.937 & \textbf{0.809} & 1.449 & 0.873 & 0.818 \\
 & AUC & 0.598 & 0.598 & 0.593 & 0.644 & \textbf{0.723} & 0.721 & 0.721 & 0.609 & 0.609 & 0.662 \\
wiki-Vote & Brier$/e(\mathbf{A})$ & 0.916 & 0.899 & 0.907 & 0.882 & 0.884 & 0.880 & 0.872 & 0.914 & 0.896 & \textbf{0.848} \\
 & log score$/H_0$ & 0.836 & 0.688 & 0.698 & 0.656 & 0.640 & 0.680 & \textbf{0.628} & 0.840 & 0.686 & 0.636 \\
 & AUC & 0.893 & 0.893 & 0.893 & 0.936 & 0.937 & 0.938 & 0.938 & 0.891 & 0.891 & \textbf{0.946} \\
soc-Epinions1 & Brier$/e(\mathbf{A})$ & 0.973 & 0.972 & 0.974 & 0.961 & 0.963 & 0.962 & 0.957 & 0.975 & 0.968 & \textbf{0.954} \\
 & log score$/H_0$ & 0.788 & 0.746 & 0.715 & 0.731 & 0.686 & 0.742 & \textbf{0.673} & 0.852 & 0.746 & 0.705 \\
 & AUC & 0.855 & 0.855 & 0.856 & 0.868 & 0.874 & 0.885 & \textbf{0.893} & 0.823 & 0.823 & 0.887 \\
\bottomrule
\end{tabular}
}
\end{table}

\begin{table}[t]
\centering
\caption{AUC of the synthesis rules of Table~2, mean (s.e.) over ten splits; best per row in bold. Logistic stacking is best on all six networks; the Platt-calibrated selected agent's drop on ca-GrQc is the floor artifact discussed in Section~5.}
\resizebox{\textwidth}{!}{\begin{tabular}{lcccccccccc}
\toprule
Network & Select & Select+affine & Hull (Brier) & Hull (log) & Cone & Least squares & Exp.\ mix & Logit stacking & BPS (span) & RDPG-16 \\
\midrule
polblogs & 0.921 (0.001) & 0.923 (0.001) & 0.928 (0.001) & 0.932 (0.001) & 0.931 (0.001) & 0.930 (0.001) & 0.921 (0.001) & \textbf{0.933 (0.001)} & 0.930 (0.001) & 0.899 (0.001) \\
email-Eu-core & 0.919 (0.001) & 0.920 (0.001) & 0.925 (0.001) & 0.930 (0.001) & 0.931 (0.001) & 0.931 (0.001) & 0.919 (0.001) & \textbf{0.932 (0.001)} & 0.931 (0.001) & 0.919 (0.002) \\
Cora & 0.655 (0.014) & 0.655 (0.014) & 0.746 (0.003) & 0.753 (0.003) & 0.752 (0.003) & 0.754 (0.003) & 0.747 (0.005) & \textbf{0.759 (0.004)} & 0.754 (0.003) & 0.685 (0.004) \\
ca-GrQc & 0.733 (0.002) & 0.733 (0.002) & 0.784 (0.010) & 0.818 (0.002) & 0.818 (0.002) & 0.809 (0.005) & 0.733 (0.002) & \textbf{0.821 (0.002)} & 0.809 (0.005) & 0.789 (0.002) \\
wiki-Vote & 0.947 (0.000) & 0.951 (0.000) & 0.949 (0.001) & 0.953 (0.000) & 0.951 (0.001) & 0.949 (0.001) & 0.947 (0.000) & \textbf{0.955 (0.000)} & 0.949 (0.001) & 0.942 (0.000) \\
soc-Epinions1 & 0.878 (0.000) & 0.880 (0.000) & 0.880 (0.001) & 0.885 (0.000) & 0.897 (0.002) & 0.877 (0.003) & 0.878 (0.000) & \textbf{0.909 (0.000)} & 0.877 (0.003) & 0.868 (0.000) \\
\bottomrule
\end{tabular}
}
\end{table}

\begin{table}[t]
\centering
\caption{The six networks: nodes, edges, density, triangle density and clustering of the observed graph.}
\small
\begin{tabular}{lrrrrr}
\toprule
Network & nodes $n$ & edges & density $e(\mathbf{A})$ & triangle density $t(\mathbf{A})$ & clustering $C(\mathbf{A})$ \\
\midrule
polblogs & 1,224 & 16,715 & 2.23e-02 & 3.31e-04 & 0.226 \\
email-Eu-core & 986 & 16,064 & 3.31e-02 & 6.62e-04 & 0.267 \\
Cora & 2,708 & 5,278 & 1.44e-03 & 4.93e-07 & 0.093 \\
ca-GrQc & 5,241 & 14,484 & 1.05e-03 & 2.01e-06 & 0.630 \\
wiki-Vote & 7,115 & 100,762 & 3.98e-03 & 1.01e-05 & 0.125 \\
soc-Epinions1 & 75,879 & 405,740 & 1.41e-04 & 2.23e-08 & 0.066 \\
\bottomrule
\end{tabular}

\end{table}

\begin{table}[t]
\centering
\caption{Node-holdout design (Section~S3): Brier score over the union density on the hidden halves of the test nodes' rows, mean (s.e.) over ten seeds, best per column in bold. Union agents: the five agents fitted to the partially observed union graph; layers: one validation-selected agent per layer; rank-16: the Platt-calibrated rank-sixteen union model.}
\label{tab:nodeholdout}
\small
\begin{tabular}{lcc}
\toprule
 & CS-Aarhus & Amazon \\
\midrule
union rank-16, calibrated & 0.770 (0.011) & 0.846 (0.003) \\
union agents: selected, calibrated & 0.675 (0.037) & 0.891 (0.003) \\
union agents: hull & 0.622 (0.030) & 0.891 (0.002) \\
union agents: span & 0.597 (0.032) & 0.890 (0.003) \\
union $+$ rank-16: span & 0.623 (0.026) & 0.840 (0.003) \\
layer agents: hull & 0.677 (0.016) & 0.894 (0.002) \\
layer agents: span & 0.575 (0.022) & 0.864 (0.003) \\
layer agents: noisy-OR & 0.575 (0.022) & 0.865 (0.003) \\
union $+$ rank-16 $+$ layers: span & 0.560 (0.025) & \textbf{0.837 (0.003)} \\
union $+$ rank-16 $+$ layers: cone & \textbf{0.545 (0.026)} & 0.838 (0.003) \\
union $+$ rank-16 $+$ layers: noisy-OR & 0.570 (0.021) & 0.849 (0.003) \\
\bottomrule
\end{tabular}

\end{table}

\section{Additional experimental detail}
\label{sec:expdetail}

\begin{figure}[t]
\centering
\includegraphics[width=\textwidth]{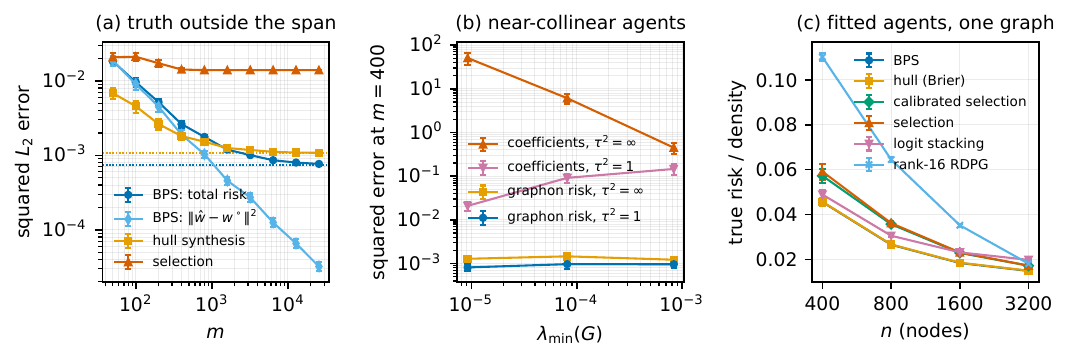}
\caption{(a) A truth outside the agent span: risk of span synthesis, hull synthesis and selection against $m$, with the approximation floors drawn. (b) Two nearly collinear agents: coefficient error and posterior spread of least squares against the posterior mean at $\tau^2=1$, as $\lambda_{\min}(G)$ falls. (c) The full fitted-agent pipeline on single graphs from a three-mechanism superposition, risk over density against $n$; $R=20$.}
\end{figure}

\begin{figure}[t]
\centering
\includegraphics[width=\textwidth]{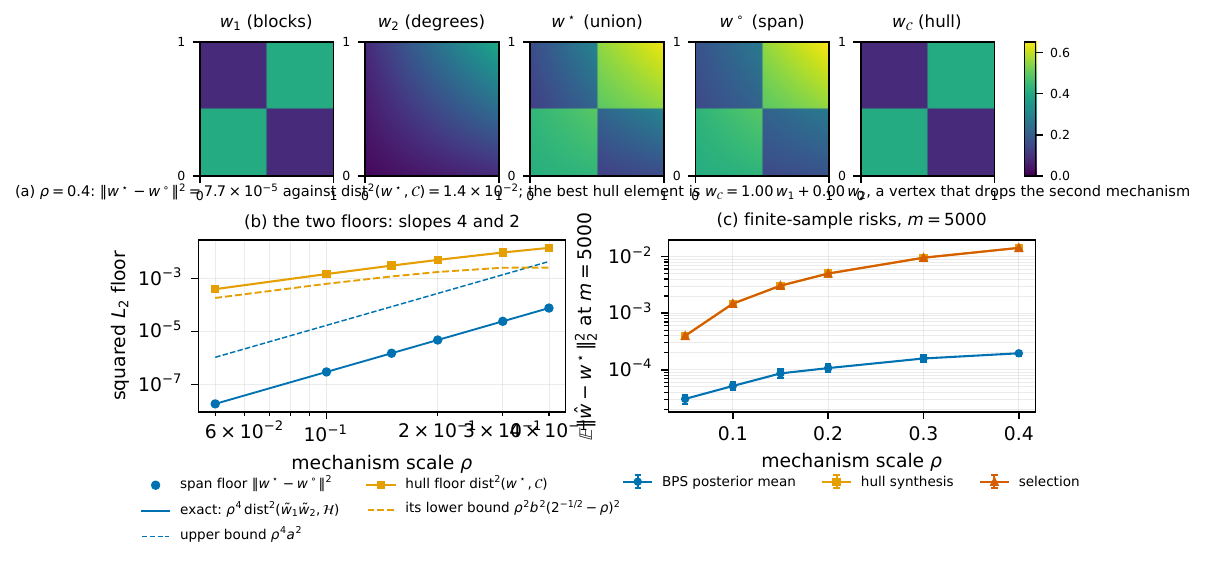}
\caption{Superposition of two mechanisms. (a) The two scaled mechanisms, their union $\wstar$, the span projection $\wproj$ and the best hull element at $\rho=0.4$. (b) The two floors against $\rho$ with the exact identity $\rho^4\dist^2(\tilde w_1\tilde w_2,\calH)$ of Theorem~1(iii) and the bounds of Theorem~1(ii). (c) Risks at $m=5000$, $R=100$; standard errors within markers.}
\end{figure}

\begin{figure}[t]
\centering
\includegraphics[width=\textwidth]{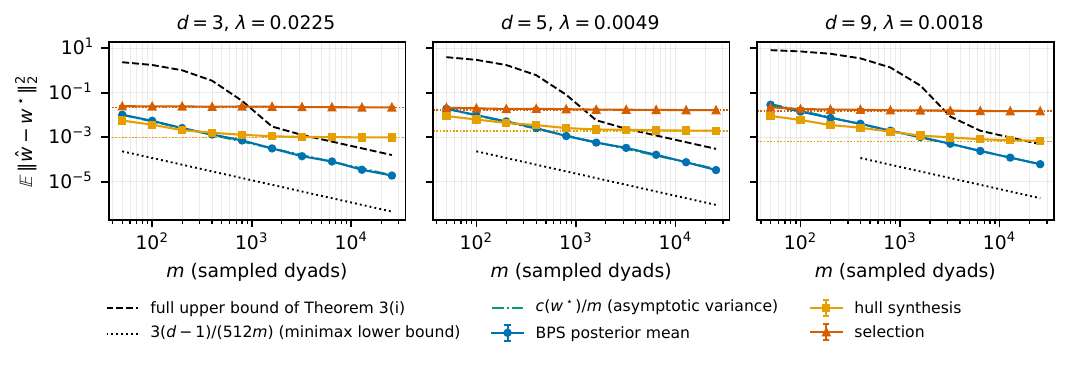}
\caption{Truth in the agent span. Risk of span synthesis, hull synthesis and selection against $m$ for $d=3,5,9$; $R=100$. Reference lines: the full upper bound of Theorem~3(i), the minimax lower bound $3(d-1)/(512m)$ drawn where $m\ge3B^2/8$, and $c(\wstar)/m$ of Proposition~S2.}
\end{figure}

\begin{figure}[t]
\centering
\includegraphics[width=\textwidth]{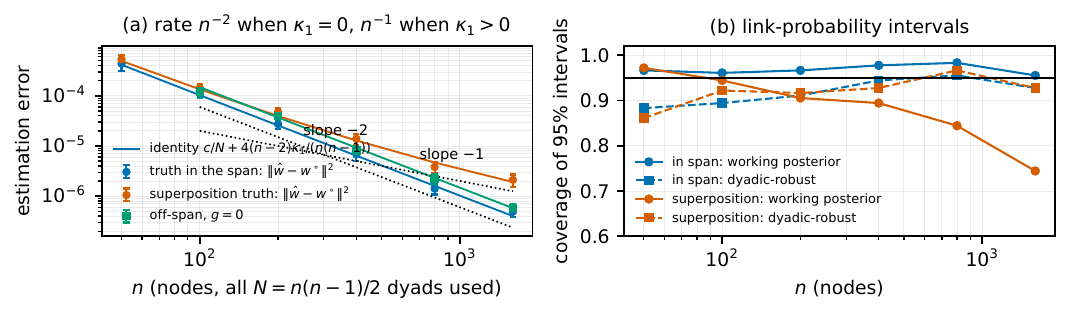}
\caption{One observed graph, all dyads, known agents. (a) Estimation error against $n$ for a truth in the span and the superposition truth, with the exact identity of Theorem~2(i) drawn through the points and reference slopes $-2$ and $-1$. (b) Coverage of nominal $95\%$ intervals for three link probabilities: working posterior against the dyadic-robust adjustment (6); $R=60$.}
\end{figure}

\paragraph{Network functionals.}
For a graphon, $e(w)=\int w$, $t(w)=\int w(u,v)w(v,z)w(z,u)\,du\,dv\,dz$ and $\varsigma(w)=\int\{\int w(u,v)\,dv\}^2du$ are the edge, triangle and wedge densities and $C(w)=t(w)/\varsigma(w)$ the clustering ratio. For a graph $A$ on $n$ nodes with degrees $d_i$, the finite versions over distinct tuples are $e(A)=\sum_id_i/\{n(n-1)\}$, $t(A)=6\,\#\{\text{triangles}\}/\{n(n-1)(n-2)\}$, $\varsigma(A)=\sum_id_i(d_i-1)/\{n(n-1)(n-2)\}$ and $C(A)=t(A)/\varsigma(A)$; observed and implied values in Table~S3 use these same conventions, so the two are directly comparable.

\paragraph{Preprocessing and scoring conventions.}
Each network is symmetrized, with self-loops, duplicate edges and isolated nodes removed and all connected components kept. Ten splits give standard errors, which describe split-to-split variability on the fixed graph, not uncertainty about a population of networks. Log scores floor predictions at $10^{-6}$ after clipping. The dependence-adjusted intervals of \citet{NadeauBengio2003} multiply the naive variance of the mean split difference by $1/J+n_2/n_1$ with $J=10$ splits and $n_2/n_1=0.2/0.8$, the scored share of the edges over the share the agents and weights are fitted on; for the corrected paired gains of Section~5 they are $[0.05,0.27]$, $[0.25,0.79]$, $[0.42,0.87]$, $[0.08,0.18]$ and $[0.07,0.14]$ percent of the density on polblogs, email-Eu-core, Cora, ca-GrQc and wiki-Vote, and $[-0.11,0.05]$ on soc-Epinions1. The more conservative ratio $0.3/0.7$ widens them further and still excludes zero on the same five networks (script \texttt{13\_dependence\_intervals.py} in the archive).

\paragraph{Prior sensitivity, the exponential mixer and clip rates.}
At $\tau^2\in\{1,0.01\}$ the synthesis Brier score is unchanged to three decimals on five networks; on Cora $\tau^2=0.01$ improves it by $0.5\%$ of the density, so shrinkage helps where the weights are unidentified. The exponential-weight mixer of \citet{LiLe2024}, run on the identical splits, coincides with selection to three decimals on five networks and beats it by $0.73\%$ of the density on Cora, where several agents nearly tie, consistent with its cross-validation limit. The de-thinning clip binds on at most $0.6\%$ of test dyads (the rank-eight agent on ca-GrQc) and below $0.3\%$ elsewhere. The Platt-calibrated selected agent's AUC on ca-GrQc, $0.683$ against $0.733$ for the raw agent, is a floor artifact: most selected probabilities sit at the logit-transform floor and tie, while the affine calibration preserves the ranking.

\paragraph{Weights, uncertainty and functionals (pooled study).}
The corrected weights (Table~S7) are predictive projection coefficients, not mechanism strengths: the dominant agent receives $0.90$ to $1.10$ on five networks and $0.65$ on Cora, the degree-histogram agent a negative weight on five, reaching $-0.51$ on Cora, and the intercepts are zero to two decimals except $-0.01$ on email-Eu-core. The posterior standard deviations come from the working model, trustworthy only when the truth is in the span; we read them as identification geometry, up to fifteen times wider on Cora, not as frequentist uncertainty. Implied triangle densities and clustering are in Table~S3; on de-thinned agents the combinations no longer dominate the single agents there, and we draw no functional-recovery claim.

\paragraph{External comparators (pooled study).}
Table~S8 places six external link predictors beside the synthesis rules on the identical held-out dyads: common neighbors, Adamic--Adar, Jaccard and preferential attachment, Platt-calibrated on the validation dyads, a logistic regression on seven topological features, and a random-forest stacker on them \citep{GhasemianEtAl2020}. The forest is the best ranker on four networks, and on ca-GrQc common neighbors, Adamic--Adar and both calibrated topological learners beat every model-based rule (Adamic--Adar at $0.55$ of the density against $0.73$ for the rank-sixteen model), while Jaccard and preferential attachment do not: the triadic information lies outside the agent span, and the constructive reading is to add such a mechanism to the library. Under the uniform-dyad Brier score the graphon rules stay best on three networks, the nested-calibrated forest wins Cora ($0.953$ against $0.987$) and the raw forest wins email-Eu-core. The forest's ranking advantage does not survive the stratified reweighting on soc-Epinions1 (Brier $7.0$ against $0.93$ times the density), but refitting it with a nested validation split, half to train and half to Platt-adjust under the stratum weights, brings that Brier to $0.94$ at unchanged AUC $0.964$: the $7.0$ was a calibration artifact, not a property of the features. The synthesis is not the best ranker; its case is calibrated probabilities, uncertainty and weights. Holding the evaluation dyads fixed and subsampling the training edges to $30$ and $50\%$ (Table~S9), the paired gain falls as training grows on polblogs ($0.96$ to $0.16$) and email-Eu-core ($3.2$ to $0.52$), is flat on ca-GrQc and wiki-Vote, and rises on Cora ($0.06$ to $0.64$): less data behind each agent means more value in combining only where one agent dominates.

\paragraph{Private shares.}
The pooled protocol gives every agent the same training graph, the regime least favorable to combination; the experiment now repeats with the training edges in five disjoint private shares, the class-to-share assignment rotated across seeds, each mechanism class fitted on its own share only and de-thinned by its own retention of about $0.14$, the decision maker holding just the validation dyads and the five fitted mechanisms: the opinion-pooling regime of \citet{West1992} and \citet{McAlinnWest2019}. Selection changes identity across splits, and the synthesis beats the calibrated selected agent by $1.0$ to $5.8\%$ of the density on five networks, from $5.8$ $[5.1,6.5]$ on polblogs to $1.0$ $[0.7,1.3]$ on soc-Epinions1, an order of magnitude above the pooled gains, with logistic stacking better still on four of the five; on Cora the span trails ($-1.5$ $[-4.1,1.1]$) and the log-score hull is the best rule. A same-algorithm control, five rank-eight spectral agents fitted one per share and combined by the same rules, separates information pooling from mechanism diversity, and it cuts against the diversity reading: the bagged copies beat the mechanism library by $3.2$ to $5.2\%$ of the density on polblogs, ca-GrQc and wiki-Vote and by $0.8$ on soc-Epinions1, while the mechanism library edges the copies on email-Eu-core ($0.7$) and beats them on Cora ($8.0$), where the spectral agent is weak. Most of the private-share gain is therefore pooling of information held in disjoint shares, which one strong algorithm fitted on every share and combined captures as well or better; mechanism diversity adds beyond bagging on two of the six networks. Fragmentation still costs: the synthesis trails the pooled calibrated agent by $3.6$ to $13.6\%$ of the density and the best combination by $3.1$ to $13.6$ on five networks, with Cora's log-score hull edging the pooled benchmark by $0.6$. Synthesis mitigates, and does not repair, the loss of pooling (Tables~S10 and S11).

\paragraph{Multiplex networks.}
The three union-of-layers networks of Table~1 are the CS-Aarhus multiplex of \citet{MagnaniEtAl2013}, $61$ members of a computer science department with five relations (lunch, Facebook friendship, coauthorship, leisure, work; $620$ layer edges, $353$ union edges), distributed with the \texttt{py3plex} library; and the Amazon and YouTube multiplexes in the form distributed by \citet{CenEtAl2019}, an Amazon co-purchase graph on $10{,}099$ products with two relation types ($130{,}899$ layer edges, $129{,}811$ union edges, so the layers are nearly disjoint) and a YouTube user graph on $2{,}000$ users with five relation types ($1{,}310{,}544$ layer edges, $835{,}330$ union edges), both taken as the union of their training, validation and test positives. Nodes isolated in the union are dropped and each layer is symmetrized with self-loops removed. The split of the union edges, the negative sampling, the stratum weights and the scoring are those of the six-network study; when fewer than ten non-edges per held-out edge exist, as on CS-Aarhus and YouTube, every non-edge of the union is used, split between validation and test in the ratio of their edge counts. Each layer graph is restricted to the training union edges, so its retention is close to $0.7$ and known exactly, and its agents are divided by it. The layers are not independent unconditionally: the density of the intersection of two layers exceeds the product of their densities by factors of $3.7$ to $9.9$ on CS-Aarhus (median $5.4$ over the ten pairs), $13$ on Amazon and $1.8$ to $8.3$ on YouTube (median $2.5$). Shared latent positions produce such marginal association even when the layer indicators are conditionally independent given the positions, so these ratios neither confirm nor refute the superposition model's independence assumption; they show that the union is far from the disjoint case and that the noisy-OR link is being used as a working model. The class of the agent on each layer is chosen on the validation dyads among Chung--Lu, degree blocks, blockmodel and rank-eight spectral; validation picks the spectral class on every layer of Amazon and YouTube and mixes classes on CS-Aarhus. The noisy-OR synthesis is fitted by L-BFGS-B on the weighted validation log-likelihood with the constraint $\gamma\ge0$ and the intercept $\gamma_0\ge0$; the logit pool is a weighted logistic regression on the agents' logits; the Bayesian model average weights the agents by $\exp$ of their weighted validation log-likelihood times the number of validation dyads, and is a selection rule at this sample size. Script \texttt{14\_multiplex.py} writes one file per split; \texttt{14b\_multiplex\_tables.py} produces Table~2 and the paired comparisons quoted in Section~5 (\texttt{results/multiplex\_summary.txt}), with paired-$t$ and dependence-adjusted intervals as in the six-network study.

\paragraph{Strong single models and the multilayer blockmodel.}
Script \texttt{16\_strong\_agents.py} regenerates the split of a seed and fits to the union training graph: universal singular value thresholding, the truncated SVD that keeps the singular values above $2.01\sqrt{n\hat\rho}$ with $\hat\rho$ the training density, computed from the sixty leading singular triplets so the rank is capped at $60$ (the cap binds on Amazon in every split; on CS-Aarhus the threshold keeps two or three); neighborhood smoothing with bandwidth $\sqrt{\log n/n}$ on the dense $n\times n$ matrix, fitted where $n\le8{,}000$; a degree-corrected blockmodel on regularized spectral clustering with $K\in\{2,4,6,8,10,12,16\}$ chosen by holding out twenty percent of the training edges and ten sampled non-edges per held-out edge, refitting for every $K$ and keeping the smallest stratified Brier score ($K=16$, the largest value offered, on Amazon in every split); and the multilayer degree-corrected blockmodel, communities from the aggregate of the training layers by the same spectral clustering with the same $K$, degree-corrected block matrices per layer divided by the layer's retention, and the union predicted by $1-\prod_l(1-p_l)$. The external baselines are those of Table~S8 computed from the union training graph, with the forest fitted twice: on all validation dyads (negatives capped at $200{,}000$), and on a random half with a Platt map fitted on the other half. Every single model enters the comparison Platt-calibrated on the validation dyads, and the strong models and the multilayer blockmodel enter the libraries \texttt{full+strong} and \texttt{full+strong+MLSBM}, whose rules are those of Section~5.

\paragraph{Node holdout.}
Script \texttt{17\_node\_holdout.py} implements the split family absent from the dyad-holdout studies. Ten percent of the nodes are test nodes and a disjoint ten percent validation nodes; every edge touching a test or validation node is hidden with probability one half (one coin per edge; an edge touching both sets counts as test); the agents are fitted to the graph with the hidden edges removed, so the held-out nodes are present but half observed; the weights are fitted on the hidden validation edges and ten sampled non-edges per hidden edge with at least one endpoint among the validation nodes, and every rule is scored on the hidden test edges and the corresponding non-edges. Since a hidden edge is observed with probability one half, agent probabilities on evaluation dyads are divided by one half, the row analogue of the de-thinning of the main protocol; stratum weights use the population density. Table~S14 reports the union library, the rank-sixteen union model, the layer agents and their combinations on CS-Aarhus and Amazon over ten seeds; the strong single models are not refitted under this design.

\paragraph{The rebuilt protocol.}
Script \texttt{19\_protocol\_v2.py} reruns the same splits with the changes the referees asked for, each recorded here so that the shipped tables and the rebuilt ones can be compared. Nested validation: the validation dyads are split in half, agent-class selection and the Platt maps of single models use one half, every combination rule is fitted on the other, and the adaptive selector of Proposition~5 fits the rules on half of the second half and chooses on the rest. Loss-matched selection: an agent selected under the log score is reported beside the Brier-selected one. Rank-preserving calibration: Platt maps are fitted to the unclipped de-thinned agent values, which removes the ties that the floor at $10^{-6}$ created (the ca-GrQc AUC drop of Table~S12 was that artifact). Structural baselines: the fixed noisy-OR $1-\prod_j(1-\hat w_j)$, the clipped sum and the max, raw and calibrated; the exact mixers of \citet{LiLe2024}, ordinary and nonnegative least squares without intercept, the latter being the no-intercept cone $\calK_0$ of Theorem~$1'$; and the intercept simplex $\calC_c$, completing the four sets. The span with $\tau^2$ chosen on held-out dyads among $\{10^{-3},\dots,10^2\}$. Strong single models with the caps lifted: universal singular value thresholding from the $200$ leading triplets and $K\in\{2,\dots,16,24,32,48\}$. And a node-level bootstrap of the key paired Brier differences, $200$ replicates with every test dyad weighted by the product of the multiplicities of its endpoints, beside the split-to-split intervals. Two points the referees raised need no rerun. The population density enters the stratum weights only through $M/N$, which the training graph determines exactly given the known retention, so a training-only estimate is the same number; and for agents that are functions of the rescaled observed adjacency, positive-edge thinning with division by the retention is the inverse-probability-weighted estimator under a uniform dyad mask, so an outcome-independent masking arm changes only which non-edges are scored and not what the agents see: the two protocols are equivalent in distribution for every rule in this paper, and differ only for methods that exploit missingness, which none of ours does.

On CS-Aarhus the rebuilt protocol (ten splits, full library of Table~1) gives, in Brier over density: selected agent $0.662$, hull $0.600$, intercept simplex $0.594$, Li--Le OLS $0.544$, Li--Le NNLS $0.537$, cone with intercept $0.518$, span $0.542$ and $0.528$ with $\tau^2$ chosen by validation ($0.1$ in nine splits), learned noisy-OR $0.568$, fixed noisy-OR $0.546$, clipped sum $0.676$, adaptive choice $0.558$ (the cone in seven splits, the hull in two); with the strong models in the library the span with validation-chosen $\tau^2$ reaches $0.499$. The four constraint sets order as Theorem~$1'$ predicts, with the two sum-to-one rules at the top and the two cones and two spans below, and the no-intercept cone $1.9$ behind the intercept cone. Two things this protocol was built to test do not generalize from this network, and we state them as they are. Learning the noisy-OR strengths instead of fixing them at one costs $2.2$ of the density here and $2.5$ on soc-Epinions1, but gains $4.3$ on polblogs, $2.4$ on wiki-Vote, $2.3$ on YouTube, $2.1$ on email-Eu-core and $0.8$ on Amazon, every one of those intervals excluding zero: learning helps on five of the nine networks and hurts on two, and there is no regime rule in our data that predicts which. Choosing $\tau^2$ on held-out dyads instead of fixing it at $100$ improves the span by $1.4$ of the density here, where the twenty-odd agents are near collinear, and changes nothing to three decimals on any other network: the prior does numerical work exactly where the Gram matrix is ill-conditioned and nowhere else. The node-bootstrap intervals are two to three times wider than the split-to-split ones on this $61$-node graph and include zero for every comparison except the span against the rank-sixteen model, which is what sixty nodes can support.

\begin{table}[htbp]
\centering
\caption{Rebuilt protocol (Section~S3): Brier over density, mean (s.e.) over ten splits, best per column in bold; rules on the base library of each network (the full library of Table~1 for the multiplexes, the five agents for the single-layer networks), the calibrated single models, and the span, noisy-OR and adaptive rules on the library augmented with the strong models.}
\label{tab:v2}
\footnotesize\setlength{\tabcolsep}{4pt}
\begin{tabular}{lc}
\toprule
 & CS-Aarhus \\
\midrule
selected agent, affine-calibrated & 0.668 (0.022) \\
selected agent, Platt (rank-preserving) & 0.682 (0.022) \\
selected under the log score & 0.674 (0.020) \\
hull & 0.604 (0.010) \\
hull with free intercept & 0.597 (0.010) \\
Li--Le OLS (no intercept) & 0.554 (0.012) \\
Li--Le NNLS (no-intercept cone) & 0.543 (0.008) \\
cone with intercept & 0.524 (0.009) \\
least squares (span) & 0.550 (0.011) \\
span, $\tau^2=100$ & 0.550 (0.011) \\
span, $\tau^2$ by validation & 0.535 (0.006) \\
logistic stacking & 0.587 (0.023) \\
noisy-OR, learned strengths & 0.567 (0.011) \\
fixed noisy-OR, calibrated & 0.542 (0.012) \\
clipped sum, calibrated & 0.671 (0.004) \\
max, calibrated & 0.596 (0.013) \\
adaptive choice (Proposition 5) & 0.565 (0.012) \\
single: RDPG-16, calibrated & 0.752 (0.008) \\
single: USVT, calibrated & 0.661 (0.017) \\
single: DCSBM, calibrated & 0.625 (0.026) \\
single: NS, calibrated & 0.695 (0.014) \\
single: best-strong, calibrated & 0.608 (0.020) \\
full+strong: BPS & 0.534 (0.013) \\
full+strong: BPS-tau & \textbf{0.512 (0.011)} \\
full+strong: NoisyOR & 0.536 (0.012) \\
full+strong: Adaptive & 0.570 (0.010) \\
\bottomrule
\end{tabular}

\end{table}

\paragraph{A learned embedding predictor.}
The comparators above are all probabilistic network models or topological learners. Script \texttt{20\_embedding\_baseline.py} adds the predictor a machine-learning reader would ask for, short of a graph neural network, which we exclude deliberately: node2vec \citep{GroverLeskovec2016} with $p=q=1$, that is DeepWalk random walks with skip-gram and negative sampling, implemented in the archive's four dependencies. Ten walks of length eighty from every node of the training graph, window five, sixty-four dimensions, five negative samples per positive, three epochs of minibatch stochastic gradient descent over two million subsampled centre--context pairs with a linearly decaying step; link prediction by weighted logistic regression on the Hadamard product of the endpoint embeddings, fitted on the validation dyads, Platt-calibrated there, and scored on the test dyads, so it enters Table~S8 on exactly the terms every other single model does. No test edge is visible to the walks. Three numerical points, since the recipe is standard but the regime is not: the embeddings are normalized to unit length and the Hadamard features standardized before the logistic fit, which is ridge regularized, because with sixty-four features and few positive validation dyads the unpenalized fit separates and its coefficients diverge; and each node's accumulated gradient within a minibatch is clipped in norm, because on a small or dense graph (CS-Aarhus has $61$ nodes) a node is touched hundreds of times per batch and the summed step overflows.

The result is the same dissociation between ranking and probability that Table~S8 reports for the random forest, and it is sharper. Under the uniform-dyad Brier score the embedding is worse than every model-based rule on all nine networks, by margins that are not close: $0.478$ of the density against $0.180$ for the thresholded spectral estimator on YouTube, $0.917$ against $0.709$ on Amazon, $0.942$ against $0.781$ on wiki-Vote, $1.000$ against $0.923$ on soc-Epinions1, and Platt calibration on the validation dyads does not move it. Under the area under the curve it is competitive, best of all single models on Cora ($0.801$ against $0.780$) and ca-GrQc ($0.875$ against $0.864$), and within a few points of the best elsewhere. A reader who asks whether a learned embedding would overturn the comparisons of Section~5 of the main text therefore has an answer: it ranks dyads about as well as the model-based rules and estimates probabilities considerably worse, so it changes nothing under the scoring rule this paper uses and would change the ordering under a rule this paper does not use. Per-split results are in \texttt{results/embed\_*.json}.

\paragraph{Benchmarking rules.}
The thinning analysis of Section~5 reduces to four rules for anyone benchmarking network model combination.
\begin{enumerate}
\item Never learn free weights on a thinned-edge holdout without dividing every fitted probability by the retention; the weights will otherwise absorb $1/\text{retention}$.
\item Score combination gains against a calibrated selected agent; affine and Platt calibration can disagree about the winner.
\item Attach dependence-aware uncertainty to every paired gain; naive standard errors over reused splits understate.
\item Report the fraction of rescaled probabilities that hit the clip, since clipping concentrates on the dyads that drive the Brier score.
\end{enumerate}

\paragraph{Misspecification, collinearity and the fitted-agent pipeline (Figure~S1).}
Panel (a) uses a truth with a diagonal band and an oscillating component outside the span of three agents, $\norm{\wstar-\wproj}^2=7.35\times10^{-4}$: the total risk of the span synthesis reaches $7.67\times10^{-4}$ at $m=25{,}600$ with estimation part $3.2\times10^{-5}$, while hull synthesis and selection stop at their floors, $1.08\times10^{-3}$ and $1.40\times10^{-2}$. Panel (b) takes two agents differing by a perturbation, with $\lambda$ from $9\times10^{-6}$ to $8\times10^{-4}$: least squares has coefficient error up to $51$ while the posterior mean with $\tau^2=1$ has coefficient error $0.02$ and posterior standard deviation near one, with stable graphon risks for both. Panel (c) runs the full pipeline on single graphs from a superposition of a four-block, a degree and a latent-distance component (density $0.095$): Chung--Lu, degree-histogram blockmodel, spectral blockmodel and rank-eight random dot product agents are fitted to $70\%$ of the edges and de-thinned as in Section~5, weights are learned on validation dyads, and the risk is the distance to the generating matrix by quadrature, $R=20$. At $n=1600$ the calibrated selected agent has risk $0.0229$ times the density, hull synthesis $0.0185$, span synthesis $0.0184$, logistic stacking $0.0231$ and the rank-sixteen random dot product graph $0.0352$, with projection weights $(-0.10,-0.01,0.39,0.64)$ summing to $0.92$.

\paragraph{Computational limits.}
The one-graph experiment was attempted at $n=3200$, where the in-span block alone requires a Monte Carlo over $5.1\times10^{6}$ dyads per replication; the runtime was prohibitive in the environment used, and the reported sizes stop at $n=1600$. The archive contains the script and the partial output, so the extension can be run where more compute is available.

\paragraph{Superposition (Figure~S2).}
Figure~S2 shows the regime of Theorem~1 with $J=2$: a two-block assortative graphon and a rank-one degree graphon, each scaled by $\rho$, and the union truth. The heatmaps make the geometry visible: the union contains both the block pattern and the degree gradient, the span projection reproduces both and misses only the interaction, and the best hull element at $\rho=0.4$ is the vertex $w_1$, dropping the degree mechanism entirely, with error $1.4\times10^{-2}$, about $180$ times the span floor $7.7\times10^{-5}$. Panel (b) is the exact test: the span floor divided by $\rho^4$ equals $\dist^2(\tilde w_1\tilde w_2,\calH)=3.00\times10^{-3}$ at every $\rho$, the identity of Theorem~1(iii), and the hull floor tracks slope two between its bounds. Panel (c) shows the finite-sample side at $m=5000$: hull and selection sit on their floors, the span sits on its estimation error, and their ratio grows from $13$ at $\rho=0.05$ to $73$ at $\rho=0.4$, the crossover of Section~3.1 read at fixed $m$.

\paragraph{Truth in the span (Figure~S3).}
Three libraries of size $J=2,4,8$ ($d=3,5,9$), built from blockmodel, rank-one, latent-distance and periodic graphons, are combined with a negative intercept and weights summing to more than one, so $\wstar\in\calW_{\calH}$ but $\wstar\notin\calC$. Figure~S3 shows the risk against $m$ with the full bound of Theorem~3(i), the minimax lower bound over its range, and $c(\wstar)/m$. The realized risk lies on $c(\wstar)/m$ from $m=50$, with $m\E\norm{\what-\wstar}^2/d$ between $0.15$ and $0.21$; the full bound is a factor of eight above it once the Gram term dies, at $m\approx1500$ to $3900$, so the bound is informative in the plotted range. Hull synthesis and selection flatten at their floors, $\dist^2(\wstar,\calC)=9.8\times10^{-4}$, $2.0\times10^{-3}$, $6.6\times10^{-4}$ and $\min_j\norm{w_j-\wstar}^2=2.2\times10^{-2}$, $1.7\times10^{-2}$, $1.5\times10^{-2}$, and the realized risks first cross at $m=400$, $400$ and $1600$; clipping matters only at $m=50$. Coverage under (2) is in Table~S6: with $\nu=1/4$ the link-probability intervals cover at $0.95$ to $1.00$ at $m\in\{400,1600,6400\}$, with posterior spread above risk by the predicted factor $\nu d/c(\wstar)$, $1.47$ and $1.36$ for the two truths; with the plug-in $\nu$ of Section~5 coverage drops to $0.88$ at some dyads, the anti-conservative case of Proposition~S2.

\end{document}